\documentclass[final]{siamltex}
\usepackage{xcolor}
\usepackage{amsmath}
\usepackage{amssymb}
\usepackage{mathrsfs}
\usepackage{graphicx}
\usepackage{bm}
\usepackage{ucs}
\usepackage[utf8x]{inputenc}
\usepackage{wrapfig}
\usepackage{color}
\usepackage{float}
\usepackage{epstopdf}
\usepackage{algorithm,algorithmic}
\usepackage[breaklinks,colorlinks=true,linkcolor=blue,citecolor=red,backref=page]{hyperref}

\def\vect#1{\mbox{\boldmath{$#1$}}}
\def\Zc{{\cal Z}}

\renewcommand{\hat}{\widehat}
\newcommand{\om}{\omega}
\newcommand{\la}{\lambda}

\newcommand{\br}{{\mathbf{r}}}
\newcommand{\vr}{{\vec{\br}}}
\newcommand{\vy}{{\vec{\mathbf{y}}}}

\newcommand{\vt}{{\vec{\bf t}}}
\newcommand{\vm}{{\vec{\bf m}}}
\newcommand{\vn}{{\vec{\bf n}}}
\newcommand{\bA}{\mathbf{A}}
\newcommand{\bd}{\mathbf{d}}

\newcommand{\cI}{{\mathcal{I}}}
\newcommand{\cN}{{\mathcal{N}}}
\newcommand{\cA}{{\mathcal{A}}}
\newcommand{\bD}{{\bf D}}
\newcommand{\cY}{{{Y}_\alpha}}
\newcommand{\cYp}{{{Y}_\alpha^\perp}}

\newcommand{\brho}{\boldsymbol{\rho}}

\begin{document}

\title{Synthetic Aperture Imaging of Direction and Frequency Dependent
  Reflectivities} \author{Liliana Borcea\footnotemark[1] \and Miguel
  Moscoso \footnotemark[2] \and George Papanicolaou\footnotemark[3]
  \and Chrysoula Tsogka\footnotemark[4]} \maketitle


\renewcommand{\thefootnote}{\fnsymbol{footnote}}
\footnotetext[1]{Department of Mathematics, University of Michigan,
  Ann Arbor, MI 48109 {\tt borcea@umich.edu}}
\footnotetext[2]{Gregorio Mill\'{a}n Institute, Universidad Carlos III
  De Madrid, Madrid 28911, Spain {\tt moscoso@math.uc3m.edu}}
\footnotetext[3]{Stanford Mathematics Department, 450 Serra Mall Bldg.
  380, Stanford CA 94305 {\tt papanico@math.stanford.edu}}
\footnotetext[4]{Mathematics and Applied Mathematics,
 University of Crete and IACM/FORTH, GR-71409 Heraklion,
 Greece {\tt tsogka@uoc.gr}}
\renewcommand{\thefootnote}{\arabic{footnote}}

\markboth{L. BORCEA, M. MOSCOSO, G. PAPANICOLAOU, AND
  C. TSOGKA}{SYNTHETIC APERTURE IMAGING}

\begin{abstract}
We introduce a synthetic aperture imaging
framework that takes into consideration directional dependence
of the reflectivity that is to be imaged, as well as its frequency
dependence. We use an $\ell_1$ minimization approach that is
coordinated with data segmentation so as to fuse information from
multiple sub-apertures and frequency sub-bands. We analyze this approach
from first principles and assess its performance with numerical
simulations in an X-band radar regime.
\end{abstract}

\section{Introduction}
{We introduce and analyze a novel algorithm for  
synthetic aperture radar (SAR) imaging, where a  moving
receive-transmit platform probes a remote region with signals $f(t)$
and records the scattered waves. The platform spans a 
large synthetic aperture so that high resolution images of the region may
be obtained by processing the recorded data. A related application is inverse 
synthetic aperture radar (ISAR), where the receive-transmit antenna is 
stationary, and the synthetic aperture is due to the motion of an unknown scatterer. 
If this motion is known or can be estimated, the problem can be restated mathematically 
as  SAR imaging of the scatterer, using the reference frame that moves with it.}

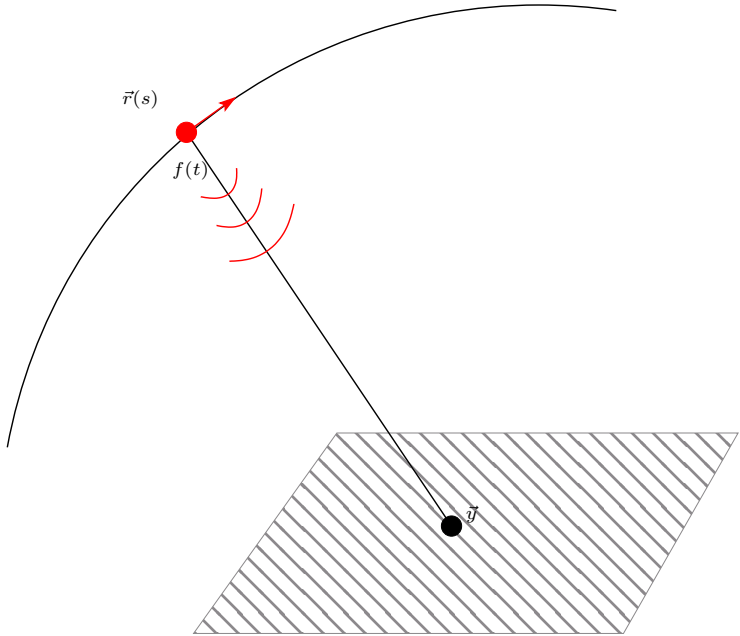
\begin{figure}[t]
\centering  
\begin{picture}(0,0)%
\includegraphics{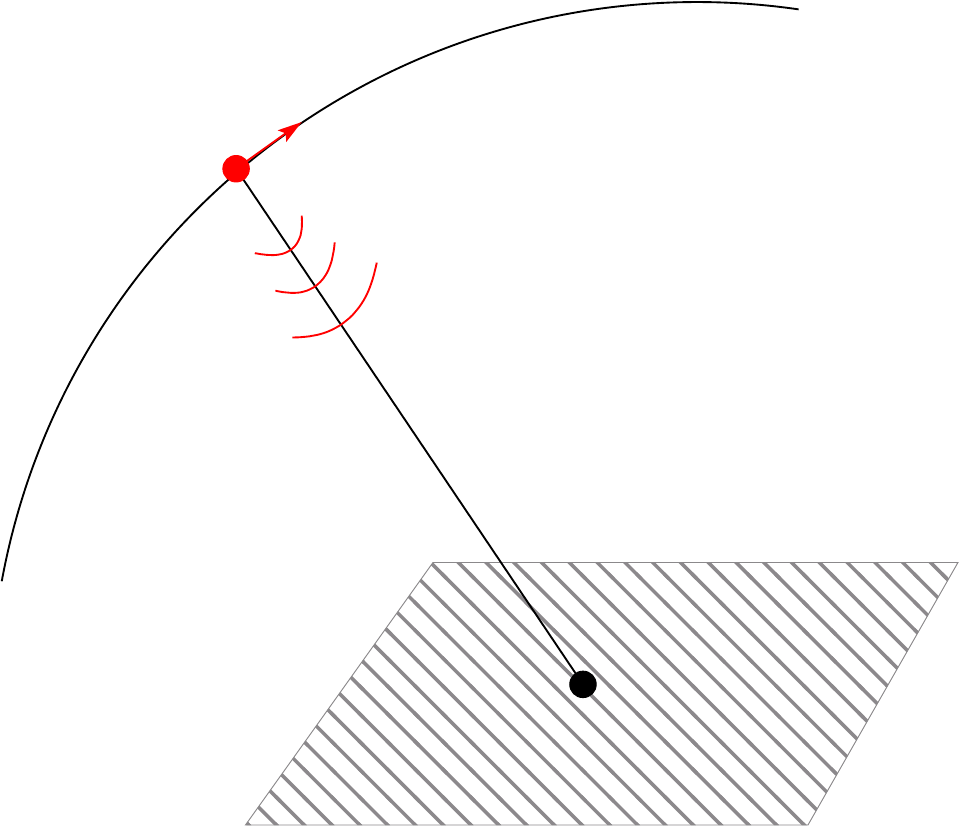}%
\end{picture}%
\setlength{\unitlength}{2368sp}%
\begingroup\makeatletter\ifx\SetFigFont\undefined%
\gdef\SetFigFont#1#2#3#4#5{%
  \reset@font\fontsize{#1}{#2pt}%
  \fontfamily{#3}\fontseries{#4}\fontshape{#5}%
  \selectfont}%
\fi\endgroup%
\begin{picture}(7677,6611)(1036,-3073)
\put(5851,-1861){\makebox(0,0)[lb]{\smash{{\SetFigFont{7}{8.4}{\familydefault}{\mddefault}{\updefault}{\color[rgb]{0,0,0}$\vec{y}$}%
}}}}
\put(2251,2489){\makebox(0,0)[lb]{\smash{{\SetFigFont{7}{8.4}{\familydefault}{\mddefault}{\updefault}{\color[rgb]{0,0,0}$\vec{r}(s)$}%
}}}}
\put(2776,1739){\makebox(0,0)[lb]{\smash{{\SetFigFont{7}{8.4}{\familydefault}{\mddefault}{\updefault}{\color[rgb]{0,0,0}$f(t)$}%
}}}}
\end{picture}%
  \caption{{Setup for imaging with a synthetic aperture.}}
  \label{fig:schematic}
\end{figure}

A schematic of the SAR imaging setup is in Figure
\ref{fig:schematic}.  The recordings $u(s,t)$ at the moving receive-transmit platform 
depend on two time variables: the slow time $s$ and the fast time
$t$. The slow time parametrizes the trajectory of the platform, and it
is discretized in uniform steps $h_s$, called the pulse repetition
rate. At time $s$ the platform is at location $\vr(s)$. It emits the
signal $f(t)$ and receives the backscattered returns $u(s,t)$. The
fast time $t$ runs between consecutive signal emissions $t \in (0,h_s)$,
and we assume a separation of time scales: The duration of
$f(t)$ is smaller than the round trip travel time of the waves between
the sensor and the imaging region, and the latter is smaller than
$h_s$.

In the usual synthetic aperture image formulation the reflectivity is
modeled as a two dimensional function of location $\vy$ on a surface
of known topography, say flat for simplicity. The assumption is that
each point on the surface reflects the waves the same way in all directions, independent
of the direction and frequency of the incident waves. This simplifies
the imaging process and makes the inverse problem formally determined: 
the data are two-dimensional and so is the unknown reflectivity function.

The reflectivity can be reconstructed by the reverse
time migration formula
\cite{skolnik1970radar,Curlander,Jakowatz,cheney}
\begin{equation}
  \cI(\vy) = \sum_{j}\int dt\,  u(s_j,t)
  \overline{f\big(t-2\tau(s_j,\vy)\big)}.
  \label{eq:migration}
\end{equation}
Here $s_j$ are the slow time emission-recording instants, spaced by $h_s$, and the
image is formed by superposing over the platform trajectory the data
$u(s_j,t)$, match-filtered with the time reversed emitted signal $f(t)$,
delayed by the roundtrip travel time $2\tau(s_j,\vy) = 2
|\vr(s_j)-\vy|/c$ between the platform location $\vr(s_j)$ and the
imaging point $\vy$. The bar denotes complex conjugate and $c$ is the
wave speed in the medium which is assumed homogeneous.

The assumption of an isotropic reflectivity may not always be
justified in applications.  Backscatter reflectivities are in general
functions of five variables: the location $\vy$ on the known (flat)
surface, the two angles of incidence and the frequency. Thus, the
inverse problem is underdetermined and we cannot expect a
reconstruction of the five dimensional reflectivity with a migration
approach.  Direct application of (\ref{eq:migration})
will produce low-resolution images of some effective, position-dependent 
reflectivity, and  there will be no information about the directivity and frequency dependence
of the actual reflectivity. 

The reconstruction of frequency dependent reflectivities with synthetic aperture radar 
has been considered in \cite{cheney2013imaging}, where Doppler 
effects are shown to be useful in inversion, and in \cite{sotirelis2012study,elachi1990radar}, where 
data are segmented over frequency sub-bands, and then images are  formed separately, for each data subset. 
Data segmentation is a natural idea, and we show here how to use it for reconstructing both frequency and direction dependent reflectivities.

The main result in this paper is the introduction
and analysis of an algorithm for imaging direction and
frequency dependent reflectivities of strong, localized
scatterers. This algorithm is based on $\ell_1$ optimization.
It reconstructs reflectivities of localized scatterers by seeking among all
those that fit the data model the ones with minimal spatial support. 
Array imaging algorithms based on  $\ell_1$ optimization are proposed and analyzed in
\cite{baraniuk2007compressive,potter2010sparsity,chai2014imaging,
  chai2013robust,
  fannjiang2010compressed,fannjiang2012coherence,borcea2015resolution}.
They consider only isotropic, frequency independent reflectivities.

A direct extension of $\ell_1$ optimization
methods to imaging direction and frequency dependent reflectivities
amounts to solving a grand optimization problem for a very long vector
$\brho$ of unknowns, the discretized reflectivity over spatial
locations on the imaging grid, the angles of incidence/backscatter and
the frequency. It has considerable computational complexity because of
the high dimension of the space in which the discretized reflectivity vector
lies. It also does not take into account the fact that many unknowns are tied to the same
spatial location points within the discretized image window.

The synthetic aperture imaging algorithm introduced in this paper is
designed to reconstruct efficiently direction and frequency dependent
reflectivities by combining two main ideas: The first is to divide the
data over carefully calibrated sub-apertures  and frequency sub-bands,
and solve an $\ell_1$ optimization problem to estimate the
reflectivity  for each data subset. {Data segmentation is useful 
assuming that  the  reflectivity
changes continuously with the direction of probing and the frequency, so that 
we can  approximate it by  a piecewise constant function, pointwise in the imaging window.  Over a sub-aperture of 
small enough linear size $a$, the platform receives scattered waves from a narrow 
cone with opening angle of the order $a/L$, where $L$ is the distance from the 
platform to the imaging window,  and we can approximate the reflectivity by that 
at the center angle. Similarly, we can approximate the reflectivity by a  constant 
over a small enough frequency sub-band. Then, we can use $\ell_1$ optimization to 
estimate the reflectivity as a function of location for each data subset.}  The size
of the sub-apertures and sub-bands determine the resolution of the
reconstruction.  The larger they are, the better the expected spatial
resolution of the reflectivity. But the resolution is worse over direction and frequency
dependence. The calibration of the data segmentation over sub-apertures and
sub-bands reflects this trade-off.  The second 
idea combines the $\ell_1$ optimizations by seeking
reflectivities that have common spatial support. Instead of a single vector
$\brho$, the unknown is a matrix with columns of spatially discretized
reflectivities. Each column corresponds to a direction of probing from
a sub-aperture and a central frequency in a sub-band. The values of
the entries in the columns are different, but they are zero
(negligible) in the same rows. Moreover, the forward model, which is
derived here from first principles, maps each column of the reflectivity
matrix to the entries in the data subsets via one common reflectivity-to-data 
model matrix. The optimization can then be carried out within the multiple
measurement vector (MMV)  formalism described in
\cite{malioutov05,cotter2005sparse,tropp2006algorithmsI,
  tropp2006algorithmsII}.

The MMV formalism is used for solving matrix-matrix equations for an unknown matrix
variable whose columns share the same support but have possibly
different nonzero values. 
We show in this paper how to reduce the synthetic aperture imaging
problem to an MMV format. The columns of the unknown matrix are
associated with the discretized spatial reflectivities for different
directions and frequencies. 
The solution of the MMV problem can
be obtained with a matrix (2,1)-norm minimization where one
seeks to minimize the $\ell_1$ norm of the vector formed by the $\ell_2$
norms of the rows of the unknown reflectivity matrix. The solutions obtained this way
preserve the common support of the columns of the unknown matrix.

This paper is organized as follows. We begin in section \ref{sect:F} with the formulation
of the imaging problem. We  derive the data model, describe the complexity of the inverse
problem, and motivate our imaging approach. The foundation of this approach is in 
section  \ref{sect:MMVRed}, where we show how to  reduce the imaging problem to an MMV format. 
The imaging algorithm is described in section \ref{sect:MMVAlg} and its performance is assessed
with numerical simulations in section \ref{sect:NUM}. The presentation in 
sections \ref{sect:F}-\ref{sect:NUM} uses the so-called start stop approximation, which neglects
the motion of the receive-transmit platform over the duration of the fast time data recording window.
This is for simplicity and also because the approximation holds in the X-band radar regime used in the 
numerical simulations. However, the imaging algorithm can include Doppler effects due to the 
motion of the receive-transmit platform, as explained in section \ref{ap:Doppler}. We end with a summary
in section \ref{sect:summary}.

\section{Formulation of the imaging problem}
\label{sect:F}
The data model is described in section \ref{sect:F1}.  Then, we review
briefly imaging of isotropic reflectivity functions via migration and
$\ell_1$ optimization in section \ref{sect:F2}. The formulation of the
problem for direction and frequency dependent reflectivities is in
section \ref{sect:F3}
\subsection{Synthetic aperture data model}
\label{sect:F1}
In synthetic aperture imaging we usually assume that the data $u(s,t)$,
depending on the slow time $s$ and the fast time $t$,
can be modeled with the single scattering approximation. For an
isotropic and frequency independent reflectivity function $\rho =
\rho(\vy)$ we have
\begin{equation}
  u(s,t) = \int \frac{d \om}{2 \pi} \hat u(s,\om) e^{-i \om t},
\label{eq:F1}
\end{equation}
with Fourier transform $\hat u(s,\om)$ given by
\begin{equation}
  \hat u(s,\om) \approx k^2 \hat f(\om) \int_{\Omega} d \vy \,
  \rho(\vy) \frac{\exp\big[2 i \om \tau(s,\vy)\big]}{(4
    \pi|\vr(s)-\vy|)^2}.
  \label{eq:F2}
\end{equation}
Here $k = \om/c$ is the wavenumber and the integral is over points
$\vy$ in $\Omega$, the support of $\rho$.  The model (\ref{eq:F2})
uses the so-called start-stop approximation, where the platform is
assumed stationary over the duration of the fast time recording
window. We use this approximation throughout most of the paper for
simplicity, and because it holds in the X-band radar regime considered in the
numerical simulations. However, the results extend to other regimes,
where Doppler effects may be important, as explained in section
\ref{ap:Doppler}.

The inverse problem is to
invert relation (\ref{eq:F2}) and thus estimate $\rho(\vy)$, given $u(s_j,t)$
at the slow time samples $s_j = (j-1)h_s$, for $j = 1, \ldots,
N_s$. Here $h_s$ is the slow time sample spacing. The inversion  is usually
done by discretizing (\ref{eq:F2}), to obtain  a linear
system of equations for the unknown vector $\brho$ of discretized
reflectivities. The support $\Omega$ in (\ref{eq:F2}) is not known, so the
inversion is done in a bounded search domain $\mathcal{Y}$ on the
imaging surface, assumed flat. We call $\mathcal{Y}$ the 
image window.  The reconstruction of $\brho$ in $\mathcal{Y}$ is a
solution of the linear system, as we review briefly in section
\ref{sect:F2}.

The discretization of $\mathcal{Y}$ is adjusted so that it is commensurate 
with the expected resolution of the image in
range and cross-range. The range
direction is the projection on the imaging plane of the unit vector pointing from the imaging location 
$\vy \in \mathcal{Y}$ to the
platform location.  The cross-range direction is orthogonal to
range. {It is well known in imaging that the range resolution is determined by the accuracy of travel time
estimation, which in turn is determined by the temporal support of $f(t)$. Thus,  it is useful 
to have a short pulse $f(t)$ whose support is of order $1/B$,
where $B$ is the bandwidth. The range resolution with such pulses is
of order $c/B$.  The cross-range resolution is proportional to the
central wavelength, which is why the emitted signals are typically modulated
by high carrier frequencies $\om_o/(2\pi)$. If $L$ is a typical
distance between the platform and the imaging window and $\cA$ is
the length of the flight path, so that the platform receives waves within 
a cone of opening angle $\cA/L$,  the
cross-range resolution is of the order $\la_o L/\cA$, where $\la_o = 2
\pi c/\om_o$ is the carrier wavelength. 
We assume that $\om_o \gg B$, which is usually the case in radar. }

In synthetic aperture imaging applications like SAR, the platform
emits relatively long signals $f(t)$ so as to carry sufficient energy to generate
strong scatter returns, and thus high signal to noise ratios. Examples of
such signals are chirps, whose frequency changes over time in an
interval centered at the carrier frequency $\om_o/(2\pi)$.  To improve
the precision of travel time estimation, and therefore range
resolution, the returns $u(s_j,t)$ are compressed in time via
match-filtering with the time reversed emitted signal
\cite{skolnik1970radar}. Moreover, to remove the large phases and
therefore avoid unnecessarily high sampling rates for the returns, the
data are migrated via travel time delays calculated with respect to a
reference point $\vy_o$ in the imaging window. The combination of
these two data pre-processing steps is called down-ramping.

For the purposes of this paper it suffices to assume that $f(t)$
is a linear chirp, in which case the Fourier transform $|\hat
f(\om)|^2$ of the compressed signal {has approximately the simple form}
\begin{equation}
  |\hat f(\om)| \approx |\hat f(\om_o)| 1_{[\om_o-\pi B , \om_o + \pi B]}(\om),
  \label{eq:L1.2}
\end{equation}
where $1_{[\om_1,\om_2]}(\om)$ denotes the indicator function of the
frequency interval $[\om_1,\om_2]$.  The down-ramped returns are 
\begin{equation}
  \int dt' \, u\Big(s,t-t' + 2 \tau(s,\vy_o)\Big) \overline{
    f(-t')} = \int \frac{d \om}{2 \pi} \, \overline{\hat
    f(\om)}\hat u(s,\om) e^{-i \om\big[t + 2\tau(s,\vy_o)\big]},
  \label{eq:L1.0}
\end{equation}
and we let $\bd$ be the vector of the samples of its Fourier transform
\begin{equation}
  \bd = \left( d(s_j,\om_l) \right)_{j = 1, \ldots N_s, l = 1, \ldots,
    N_\om}, \qquad d(s,\om) = \overline{\hat f(\om)} \hat
  u(s,\om)e^{-2 i \om \tau(s,\vy_o)} .
  \label{eq:L1.3}
\end{equation}
The size of the vector $\bd$ is $N_s N_\om$.

The linear relation between the unknown reflectivity vector $\brho$ and the
down-ramped data vector $\bd$ follows from (\ref{eq:L1.3}) and
(\ref{eq:F2}). We write it as
\begin{equation}
  \bA \brho = \bd,
  \label{eq:L1.1}
\end{equation}
where the entries in $\brho \in \mathbb{C}^{Q}$ are proportional to
$\rho(\vy_q)$, with $\vy_q$ the $Q$ discretization points of the image window
$\mathcal{Y}$, and with the constant of proportionality taken to be the area of a grid
cell. The reflectivity $\brho$ is mapped by the reflectivity-to-data matrix $\bA \in
\mathbb{C}^{N_SN_\om \times Q}$ to the data $\bd$. The assumption of frequency
independent reflectivity leads to a set of decoupled systems of
equations $\bA(\om_l) \brho = \bd(\om_l)$ indexed by the frequency
$\om_l$, where the entries of the $N_s \times Q$ matrices $\bA(\om_l)$
are 
\begin{equation}
  \bA_{j,q}(\om_l) = \frac{k_l^2|\hat f(\om_l)|^2}{ (4 \pi
    |\vr(s_j)-\vy_q|)^2} e^{2 i \om_l \big[
      \tau(s_j,\vy_q)-\tau(s_j,\vy_o)\big]}.
  \label{eq:L1.4}
\end{equation}
Here $k_l = \om_l/c$, $l = 1, \ldots, N_\om$, $j = 1, \ldots, N_s,$
and $ q = 1, \ldots, Q.$

\subsection{Imaging isotropic reflectivities}
\label{sect:F2}
Imaging of the isotropic reflectivities amounts to inverting the
linear system (\ref{eq:L1.1}). When this system is underdetermined,
there are two frequently used choices for picking a solution: either minimize the
Euclidian norm of $\brho$ or its $\ell_1$ norm. The first choice gives
\begin{equation}
  \brho = \bA^\dagger \bd,
  \label{eq:F3}
\end{equation}
where $\bA^\dagger$ is the pseudo-inverse of $\bA$. If $\bA$ is full
row rank, $ \bA^\dagger = \bA^* (\bA \bA^*)^{-1}$. The inversion
formula (\ref{eq:F3}) also applies to overdetermined problems, where
$\brho$ is the least squares solution and $\bA^\dagger = (\bA^*
\bA)^{-1} \bA^*, $ for full column rank $\bA$. The choice of the
imaging window $\mathcal{Y}$ and its discretization is an essential part of the
imaging process and, depending on the objectives and available prior
information, we may be able to control whether the system (\ref{eq:L1.1})
is overdetermined or not. We explain in Appendix
\ref{ap:discr} that by discretizing $\mathcal{Y}$ in steps commensurate
with expected resolution limits we can make the columns of $\bA$ nearly
orthogonal. This means that in the overdetermined case $\bA^* \bA$ is
close to a diagonal matrix. We also shown in Appendix
\ref{ap:discr} that in the underdetermined case, for coarse enough
sampling of the slow time $s$ and frequency $\om$, the rows of $\bA$
are nearly orthogonal, and therefore $\bA \bA^\star$ is close to a
diagonal matrix. Thus, in both cases, $\bA^\dagger$ is
approximately $\bA^\star$ up to multiplicative factors, and we can
therefore image the support of $\brho$ with $\bA^*\bd$.  This is in fact the migration formula
(\ref{eq:migration}) written in the Fourier domain, up to a
geometrical factor, since the amplitude in (\ref{eq:L1.4}) is
approximately constant for platform trajectories that are shorter than
the imaging distance and for bandwidths $B \ll \om_o$.

If we know that the imaging scene consists of a few strong, localized
scatterers, as we assume here, a better estimate of $\brho$ is given
by the optimization
\begin{equation}
  \min \|\brho\|_1 \quad \mbox{such that} ~ ~ \|\bA \brho -\bd\|_2 \le
  \epsilon.
  \label{eq:L1.5}
\end{equation}
Here $\epsilon$ is an error tolerance, commensurate with the noise
level in the data, and $\|\cdot\|_1$ and $\|\cdot \|_2$ are the $\ell_1$ and the
Euclidian norm, respectively. We refer to
\cite{baraniuk2007compressive,potter2010sparsity,chai2013robust,
  fannjiang2010compressed,fannjiang2012coherence} for studies of
imaging with $\ell_1$ optimization. The main result in this context is
that when there is no noise so that $\epsilon=0$, the reflectivities are
recovered exactly provided that the inner products of the normalized columns of $\bA$ are 
sufficiently small.   An extension of the optimization
to nonlinear data models that account for multiple scattering effects
in $\mathcal{Y}$, is considered in \cite{chai2014imaging}.  A
resolution study of imaging with $\ell_1$ optimization is in
\cite{borcea2015resolution}.

\subsection{Imaging direction and frequency dependent reflectivities}
\label{sect:F3}
In general, backscatter reflectivities are functions of five
variables: the location $\vy \in \mathcal{Y}$, the unit direction
vector $\vm$ and the frequency $\om$. Hence,
\begin{equation}
  \rho = \rho(\vy,\vm,\om).
  \label{eq:F4}
\end{equation}
This means that the down-ramped data model is more complicated than
assumed in equations (\ref{eq:F2}) and (\ref{eq:L1.0}) or,
equivalently, after discretization, in
(\ref{eq:L1.3})-(\ref{eq:L1.4}). In integral form it is given by
\begin{align}
 d(s,\om) &= \overline{\hat f(\om)} \hat u(s,\om) e^{-2 i \om
   \tau(s,\vy_o)} \nonumber \\ &= k^2 |\hat f(\om)|^2
 \int_{\Omega} d \vy \, \rho(
 \vy,\vm(s,\vy),\om) \frac{\exp \Big[ 2 i \om
     \big[ \tau(s,\vy)-\tau(s,\vy_o) \big] \Big]}{\big(4 \pi
   |\vr(s)-\vy|\big)^2},
     \label{eq:F5}
\end{align}
where $\vm(s,\vy)$ is the unit vector pointing
from the platform location $\vr(s)$ to $\vy $ in the image window $\mathcal{Y}$. In
discretized form we still have a linear system like (\ref{eq:L1.1}),
except that now $\brho$ is a vector of $Q N_s N_\om$ unknowns, the
discretized values of $\rho$ in the image window $\mathcal{Y}$.

Extending the inversion approaches described in the previous section
to this model means inverting approximately the matrix $\bA$ with a
very large number of columns. We cannot expect the migration formula
(\ref{eq:migration}) to give an accurate estimate of the reflectivity
as a function of five variables, as pointed out in the introduction. 
The $\ell_1$ optimization approach works,
but it becomes impractical for the large number $Q N_s N_\om$
of unknowns.  Moreover, it does not take into account the fact that
the entries in $\brho$ indexed by the slow time and frequency pairs
$(j,l)$, with $j = 1, \ldots, N_s$ and $l = 1, \ldots, N_\om$, refer to
the same locations $\vy_q$ on the imaging grid.

The imaging approach introduced in this paper gives an efficient way
of estimating direction and frequency dependent reflectivity functions
of strong localized scatterers in $\mathcal{Y}$. It uses an
approximation of the model (\ref{eq:F5}), motivated by the expectation
that the backscatter reflectivity should not change dramatically from
one platform location to the next and from one frequency to
another. Instead of discretizing $\rho$ over all five variables at
once, we discretize it only with respect to the location in the
image window $\mathcal{Y}$, for one probing direction and frequency at a time. 
To do so, we separate the data over subsets defined by carefully
calibrated sub-apertures and sub-bands, and freeze the direction
and frequency dependence of the reflectivity for each subset. The grand optimization 
is divided this way into smaller optimizations for $Q$
unknowns, which are then coupled by requiring that the unknown vectors
share the same spatial support in the imaging window $\mathcal{Y}$.

\section{Reduction to the Multiple Measurement Vector framework}
\label{sect:MMVRed}
\begin{figure}[t]
\centering  
\begin{picture}(0,0)%
\includegraphics{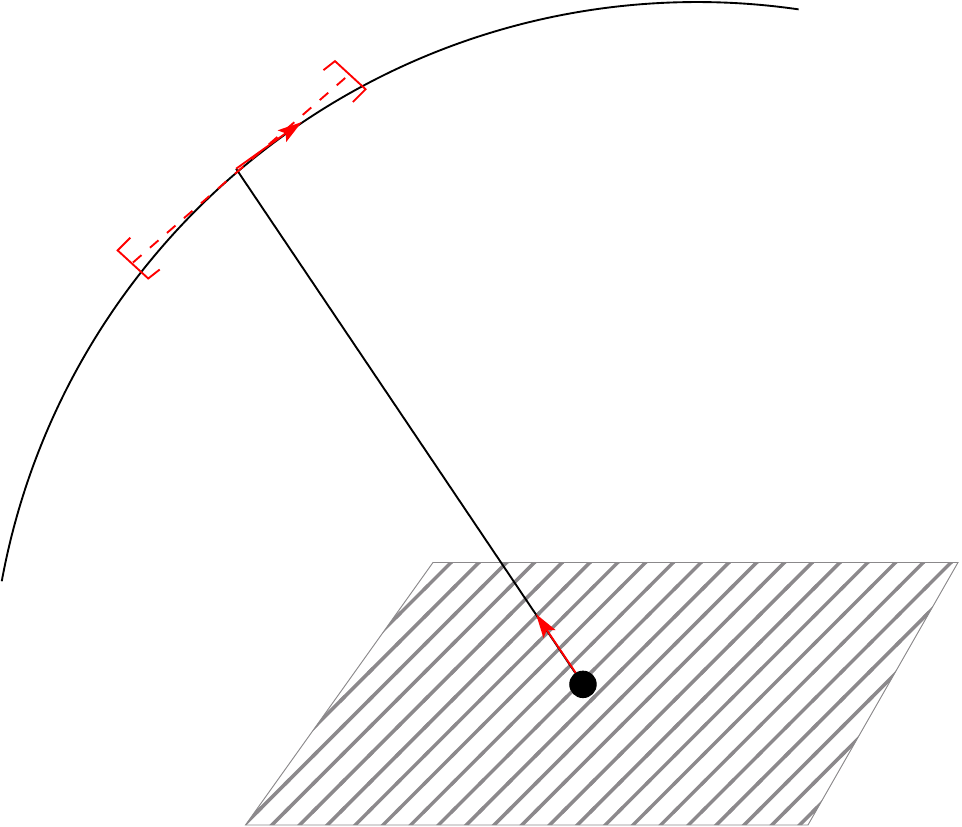}%
\end{picture}%
\setlength{\unitlength}{2368sp}%
\begingroup\makeatletter\ifx\SetFigFont\undefined%
\gdef\SetFigFont#1#2#3#4#5{%
  \reset@font\fontsize{#1}{#2pt}%
  \fontfamily{#3}\fontseries{#4}\fontshape{#5}%
  \selectfont}%
\fi\endgroup%
\begin{picture}(7677,6611)(1036,-3073)
\put(5917,-1954){\makebox(0,0)[lb]{\smash{{\SetFigFont{7}{8.4}{\familydefault}{\mddefault}{\updefault}{\color[rgb]{0,0,0}$\vy_o$}%
}}}}
\put(5458,-1250){\makebox(0,0)[lb]{\smash{{\SetFigFont{7}{8.4}{\familydefault}{\mddefault}{\updefault}{\color[rgb]{0,0,0}$\vm_\alpha$}%
}}}}
\put(3142,2107){\makebox(0,0)[lb]{\smash{{\SetFigFont{7}{8.4}{\familydefault}{\mddefault}{\updefault}{\color[rgb]{0,0,0}$\vt_\alpha$}%
}}}}
\put(2560,2291){\makebox(0,0)[lb]{\smash{{\SetFigFont{7}{8.4}{\familydefault}{\mddefault}{\updefault}{\color[rgb]{0,0,0}$\vr(s_\alpha^\star)$}%
}}}}
\end{picture}%
\caption{Schematic of the geometry for one sub-aperture
centered at the location $\vr(s_q^\star)$ of the receive-transmit platform.}
  \label{fig:geometry}
\end{figure}
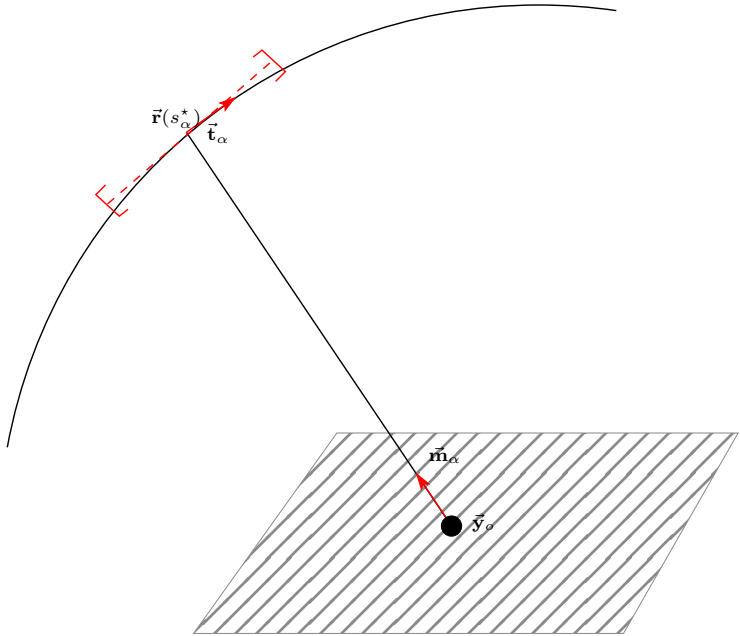
We present here an analysis of how we can write the linear relation between the
direction and frequency dependent reflectivity and the data as a
linear matrix system
\begin{equation}
  \boldsymbol{\mathbb{A}} {\bf X} = \bD,
  \label{eq:MMVEQ}
\end{equation}
where the unknown is the matrix ${\bf X}$ with $Q$
rows. The entries in the rows correspond to the discretization of this
reflectivity at the $Q$ grid points $\vy_q$ in $\mathcal{Y}$. Each
column of ${\bf X}$ depends on the reflectivity at
the backscattered direction defined by the center of a sub-aperture
and the center frequency of a sub-band. The data are segmented over
$\cN_\alpha$ sub-apertures and $\cN_\beta$ sub-bands and are grouped in
the matrix $\bD$.  The objective of this section is to describe the
data segmentation and derive the linear system (\ref{eq:MMVEQ}), which
can be inverted with the MMV approach as explained in section
\ref{sect:MMVAlg}.


We begin in section \ref{sect:single} with a single sub-aperture and
sub-band. We show in Lemma \ref{lem.1} that with proper calibration
of the sub-aperture and sub-band size, the reflectivity-to-data matrix
has a simple approximate form. Its entries have nearly constant
amplitudes while the phases depend linearly on the slow time and
frequency parametrizing the data subset. This simplification allows us
to transform the linear system via coordinate rotation to a reference one,
for all data subsets, as shown in section \ref{sect:MMVProb}.  The
matrix $\boldsymbol{\mathbb{A}}$ in (\ref{eq:MMVEQ}) corresponds to
the reference sub-aperture and sub-band, and the statement of the result is in
Proposition \ref{lem.2}.

\subsection{The sub-aperture and sub-band segmentation}
\label{sect:setup}
We enumerate the sub-apertures by $\alpha = 1, \ldots, \cN_\alpha$,
and denote by $s_\alpha^\star$ the slow time that corresponds to their center
location $\vr(s_\alpha^\star)$.  The choice of the sub-aperture size $a$ is
important, and we address it in the next section.  For now it suffices
to say that it is small enough so that we can approximate it by a line
segment, as illustrated in Figure \ref{fig:geometry}. The unit tangent
vector along the trajectory, at the center of the sub-aperture, is denoted
by $\vt_\alpha$, and the platform motion will be assumed uniform, at
speed $V \vt_\alpha$.  The unit vector from the reference location
$\vy_o$ in the image window to $\vr(s_\alpha^\star)$ is $\vm_\alpha$. We call
it the range vector for the $\alpha$ sub-aperture. The range (distance) to
the imaging window is
\begin{equation}
  L_\alpha = |\vr(s_\alpha^\star)-\vy_o|.
\end{equation}
Each sub-aperture is parametrized by the slow time offset from
$s_\alpha^\star$, denoted by
\begin{equation}
  \Delta s = s - s_\alpha^\star \in \Big[-\frac{a}{2V},\frac{a}{2V}
    \Big].
  \label{eq:M1}
\end{equation}
We do not index it by $\alpha$ because it belongs to the same interval
for each sub-aperture. The discretization of $\Delta s$ is at the slow
time sample spacing $h_s$, and there are
\[n_s = \frac{a}{V h_s} + 1\]
sample points, where $a/(V h_s)$ is rounded to an integer.  Similarly, we
divide the bandwidth in $\cN_\beta$ sub-bands of support $b \le B$,
centered at $\om_\beta^\star$, and let $\Delta \om$ be the frequency
offset
\begin{equation}
  \Delta \om = \om - \om_\beta^\star \in \Big[-\pi b,\pi b\Big].
  \label{eq:M3}
\end{equation}
We sample the sub-band with $n_\om$ points.

The reflectivity dependence on the direction and frequency is denoted
by the superscript pair $(\alpha,\beta)$, and by discretizing it with
the $Q$ points in $\mathcal{Y}$ we obtain the vector of unknowns
$\brho^{(\alpha,\beta)} \in \mathbb{C}^{Q}$. It is mapped to the data
vector $\bd^{(\alpha,\beta)}$ with entries given by the samples of
$d(s_\alpha^\star + \Delta s, \om_\beta^\star + \Delta \om)$. The
mapping is via the $n_s n_\om \times Q$ reflectivity-to-data matrix $\bA^{(\alpha,\beta)}$
described in Lemma \ref{lem.1}.
\subsection{Reflectivity-to-data model for a single sub-aperture and sub-band}
\label{sect:single}
Here we explain how we can choose the size of the sub-apertures and
frequency sub-bands so that we can simplify
the reflectivity-to-data matrix. The calibration depends on the size
of the imaging window $\mathcal{Y}$, which is quantified with two length scales
\begin{equation}
  \cY = \max_{q=1,\ldots Q} |(\vy_q-\vy_o) \cdot \vm_\alpha|,
  \label{eq:M4}
\end{equation}
and 
\begin{equation}  \cYp =
  \max_{q=1,\ldots Q} |\mathbb{P}_\alpha(\vy_q-\vy_o)|.
  \label{eq:M4p}
\end{equation}
Here $\mathbb{P}_\alpha = I - \vm_\alpha \vm_\alpha^T$ is the
projection on the cross-range plane orthogonal to $\vm_\alpha$, and
$I$ is the identity matrix. The length scale $\cY$ gives the size
of $\mathcal{Y}$ viewed from the range direction $\vm_\alpha$, and
$\cYp$ is the cross-range size.

The first constraints on the aperture $a$ and the cross-range size
$\cYp$ of the imaging window state that they are not too small, and
thus imaging with adequate resolution can be done with the data subset. Explicitly, we ask that
for all $\alpha = 1, \ldots, \cN_\alpha$,
\begin{align}
  \label{eq:M5}
  \frac{a^2}{\la_o L_\alpha} \gtrsim \frac{a \cYp}{\la_o L_\alpha}
  \gtrsim \frac{(\cYp)^2}{\la_o L_\alpha} & \gtrsim 1.
\end{align}
The inequalities on the left involve { two Fresnel numbers $a^2/(\la_o L_\alpha)$ and $(\cYp)^2/(\la_o L_\alpha)$, whose magnitudes define the imaging regime. If these numbers were small, we would be in a Fraunhofer diffraction regime, with approximately planar wavefronts  on 
the scale of the sub-aperture and of the size of the imaging window. We consider a Fresnel diffraction regime, 
where these numbers are larger and we can get  better resolution of images.}
The cross-range resolution is $\la_o
L_\alpha/a$, and naturally, the middle inequality in (\ref{eq:M5})
says that the image window is larger than the resolution limit. In
the range direction we suppose that
\begin{equation}
  \cY \gtrsim \frac{c}{b} \gg \la_o,
  \label{eq:M6}
\end{equation}
where $c/b$ is the range resolution for the sub-bands, and we used that $b \le B
\ll \om_o$.

{While we would like to have $a$ and $b$ large so as to get good spatial
resolution of the unknown reflectivity, we recall that $\rho$ is
frozen in our discretization in the small frequency sub-band and in  the narrow cone of opening angle 
of the order $a/L_\alpha$, with axis defined by the center $\vr(s^\star_\alpha)$ of the sub-aperture and 
the reference point $\vy_o$. 
The larger $a$ and $b$ are, the coarser the estimation of the
direction and frequency dependence of $\rho$. The more rapid  the variation of $\rho$ with direction and frequency,
the smaller $a$ and $b$ should be to represent it, at the expense of resolution. }

There is also a trade-off between resolution and the
complexity of the inversion algorithm. By constraining $a$ and $b$ so
that
\begin{align}
  \frac{b}{\om_o} \frac{\cYp}{\la_o L_\alpha/a}  \ll 1, \label{eq:M8}
\end{align}
and
\begin{align}
  \frac{a^2\cY}{\la_o L_\alpha^2} \ll 1, \qquad \frac{a^2\cYp}{\la_o
    L_\alpha^2} \ll 1, \label{eq:M10}
\end{align}
we can simplify the mapping between the reflectivity and the data
subset, as stated in Lemma \ref{lem.1}. This simplification allows us
to use the efficient MMV framework to solve the large optimization
problem for the entire data set, by considering jointly the smaller
problems for the segmented data in an automatic way. The key
observation here is that the unknown reflectivities for each data
subset share the same spatial support. This is what the MMV formalism
is designed to capture.

The next lemma gives the form of the reflectivity-to-data matrix in the linear system
\begin{equation}
  \bA^{(\alpha,\beta)} \brho^{\alpha,\beta} = \bd^{(\alpha,\beta)},
  \label{eq:lem.1p}
\end{equation}
for the $(\alpha,\beta)$ data subset. It is an approximation of the
system (\ref{eq:L1.1}) restricted to the rows indexed by the $n_s$
slow times in the $\alpha-$aperture and the $n_\om$ frequencies in the
$\beta-$band. The expression of $\bA^{(\alpha,\beta)}$ is derived in
appendix \ref{ap:proofs}.
\begin{lemma}
  \label{lem.1}
  Under the assumptions (\ref{eq:M5})-(\ref{eq:M10}), and with the
  pulse model (\ref{eq:L1.2}), the matrix
  $\bA^{(\alpha,\beta)}$ consists of $n_\om$ blocks
  $\bA^{(\alpha,\beta)}(\Delta \om_l)$ indexed by the frequency offset
  $\Delta \om_l$, for $l = 1, \ldots, n_\om$. Each block is an $n_s
  \times Q$ matrix with entries defined by
  \begin{align}
    A^{(\alpha,\beta)}_{j,q}(\Delta \om_l) = \frac{k_o^2 |\hat
      f(\om_o)|^2}{(4 \pi L_\alpha)^2} \exp \Big\{ -2 i (k_\beta +
    \Delta \om_l/c) \vm_\alpha \cdot \vy_q  \nonumber \\ -2 i k_\beta
    V \Delta s_j \frac{ \vt_\alpha \cdot \mathbb{P}_\alpha \Delta
      \vy_q}{L_\alpha} + i k_\beta \frac{\Delta \vy_q \cdot
      \mathbb{P}_\alpha \Delta \vy_q}{L_\alpha} \Big\},
    \label{eq:lem1}
  \end{align}
  where $k_\beta = \om_\beta^\star/c$, and $\Delta \vy_q = \vy_q -
  \vy_o$.
\end{lemma}

\subsection{Multiple sub-aperture and sub-band model as an MMV system}
\label{sect:MMVProb}
It remains to show how to write equations (\ref{eq:lem.1p}) in the
matrix form (\ref{eq:MMVEQ}) with a reflectivity-to-data matrix independent
of the sub-apertures and sub-bands.  This is accomplished via a
rotation, that brings all the sub-apertures to a single reference
sub-aperture.  But to do this, we need to know that each data subset has a
similar view of the image window. Mathematically, this is expressed
by the following two additional constraints on $a$ and $b$
\begin{align}
  \max_{1 \le \alpha \le \cN_\alpha, 1 \le q \le Q} \frac{b}{c} \big|
  (\vm_\alpha - \vm_1) \cdot \Delta \vy_q \big| \ll 1,
  \label{eq:lem2.1}
\end{align}
and
\begin{align}
  \max_{1 \le \alpha \le \cN_\alpha, 1\le \beta \le \cN_\beta, 1 \le q
    \le Q} \Big| \Big(\frac{a k_\beta}{L_\alpha} \vt_\alpha \cdot
  \mathbb{P}_\alpha -\frac{a k_1}{L_1} \vt_1 \cdot \mathbb{P}_1 \Big)
  \Delta \vy_q\Big| \ll 1,
  \label{eq:lem2.2}
\end{align}
The constraint (\ref{eq:lem2.1}) states that the imaging points remain
within the range resolution limit $b/c$ for all the apertures. The
constraint (\ref{eq:lem2.2}) states that the imaging points remain within
the cross-range resolution limits, as well. 

The derivation of the linear system (\ref{eq:MMVEQ}) is in appendix
\ref{ap:proofs} and the result is stated in the next proposition.
\begin{proposition}
  \label{lem.2}
  Under the same assumption as in Lemma \ref{lem.1} and in addition,
  supposing that conditions (\ref{eq:lem2.1}) and (\ref{eq:lem2.2})
  hold, we can combine the linear systems (\ref{eq:lem.1p}) in the
  matrix equation (\ref{eq:MMVEQ}). The reference sub-aperture and
  sub-band are indexed by $\alpha = 1$ and $\beta = 1$. The unknown
  matrix ${\bf X}$ has $Q$ rows and $\cN_\alpha
  \cN_\beta$ columns indexed by $(\alpha,\beta)$. Its entries are
  \begin{equation}
   X_q^{(\alpha,\beta)} = \rho_q^{(\alpha,\beta)} \exp \Big[
     - 2 i k_\beta \vm_\alpha \cdot \Delta \vy_q + i k_\beta
     \frac{\Delta \vy_q \cdot \mathbb{P}_\alpha \Delta
       \vy_q}{L_\alpha} \Big],
    \label{eq:lem2.4}
  \end{equation}
  where
  \begin{equation}
    \rho_q^{(\alpha,\beta)} =
    \rho(\vy_q,\vm_\alpha,
    \om_\beta^\star), \qquad
    \vm_\alpha=
    \frac{\vr(s_\alpha^\star)-\vy_o}{|\vr(s_\alpha^\star)-\vy_o|}.
   \label{eq:lem2.5}
  \end{equation}
  The data matrix $\bD$ has $n_s n_\om$ rows and $\cN_\alpha \cN_\beta$
  columns indexed by $(\alpha,\beta)$. We organize the equations in
  blocks indexed by the frequency $\Delta \om_l$, for $l = 1, \ldots,
  n_\om$. The entries of $\bD$ are defined in terms of the down-ramped
  data vectors $\bd^{(\alpha,\beta)}$ as
  \begin{equation}
    D_j^{(\alpha,\beta)}(\Delta \om_l) = \frac{(4 \pi
      L_\alpha)^2}{k_o^2 |\hat f(\om_o)|^2} d^{(\alpha,\beta)}(\Delta
    s_j,\Delta \om_l),
    \label{eq:lem2.6}
  \end{equation}
  where we recall that
  \begin{equation}
    d^{(\alpha,\beta)}(\Delta s_j,\Delta \om_l) = d\big(s_\alpha^\star
    + \Delta s_j,\om_\beta^\star + \Delta \om_l\big),
    \label{eq:lem2.7}
  \end{equation}
  and $d(s,\om)$ is  defined in (\ref{eq:L1.3}).
  The reflectivity to data matrix $\boldsymbol{\mathbb{A}}$ has
  $n_\om$ blocks indexed by $\Delta \om_l$, denoted by
  $\boldsymbol{\mathbb{A}}(\Delta \om_l)$. Each block is an $n_s
  \times Q$ matrix with entries
  \begin{equation}
    \mathbb{A}_{j,q}(\Delta \omega_l) = \exp \Big[ -2 i \frac{\Delta
        \om_l}{c} \vm_1 \cdot \Delta \vy_q - 2 i k_1 \frac{V \Delta
        s_j}{L_1} \vt_1 \cdot \mathbb{P}_1 \Delta \vy_q \Big].
    \label{eq:lem2.8}
  \end{equation}
\end{proposition}

Note that the product of the reflectivity-to-data matrix $\boldsymbol{\mathbb{A}}$ with each column of
${\bf X}$ can be interpreted, up to a constant
multiplicative factor, as a Fourier transform with respect to the
range offset $\vm_1 \cdot \Delta \vy$ and cross-range offset $\vt_1
\cdot \mathbb{P}_1 \Delta \vy$ in $\mathcal{Y}$. Equation
(\ref{eq:lem2.4}) shows that the columns of ${\bf X}$
differ from each other by a linear phase factor in $\Delta \vy$, which
amounts to a rotation of the coordinate system of the
$\alpha$ sub-aperture, and a quadratic factor which corrects for Fresnel
diffraction effects. Thus, the linear system (\ref{eq:MMVEQ}) gives
roughly the Fourier transform of the reflectivity $\rho$ for different
range direction views, and the imaging problem is to invert it to
estimate $\rho$.

\section{Inversion algorithm}
\label{sect:MMVAlg}
Here we describe the algorithm that estimates the reflectivity by
inverting the linear system (\ref{eq:MMVEQ}). By construction, the
columns of the $Q \times \cN_\alpha \cN_\beta$ unknown matrix
${\bf X}$ have the same spatial support, because they
represent the same spatial reflectivity function. Thus, we formulate the
inversion as a common support recovery problem for unknown matrices with
relatively few nonzero rows
\cite{rao1998sparse,chen06,cotter2005sparse,eldar10}. This 
Multiple Measurement Vector (MMV) formulation has been studied in
\cite{eldar10,chen06,rao1998sparse} and has been used successfully for
source localization with passive arrays of sensors in
\cite{malioutov05} and for imaging strong scattering scenes, where
multiple scattering effects cannot be neglected, in
\cite{chai2014imaging}.

In the MMV framework the support of the unknown matrix
${\bf X}$ is quantified by the number of nonzero
rows, that is the row-wise $\ell_0$ norm of
${\bf X}$.  If we define the set
\begin{equation}
  \operatorname{rowsupp}({\bf X})= \{q = 1,\ldots, Q
  ~ ~\mbox{s.t.} ~ ~ \|{\bf e}_q^T
  {\bf X}\|_{\ell_2}\neq0\},
  \label{eq:rowsup}
\end{equation}
where ${\bf e}_q^T {\bf X}$ is the $q-$th row of
${\bf X}$ and ${\bf e}_q$ is the vector with entry
$1$ in the $q-$th row and zeros elsewhere, then the row-wise $\ell_0$
norm of ${\bf X}$ is the cardinality of
$\operatorname{rowsupp}({\bf X})$,
\[\Xi_0({\bf X})=|
\operatorname{rowsupp}({\bf X})|.
\]
To estimate ${\bf X}$ we must to solve the
optimization problem
\begin{equation}\label{MMV.NP}
  \min\Xi_0({\bf X})\quad\text{s.t.}\quad
  \boldsymbol{\mathbb{A}} {\bf X} = \bD,
\end{equation}
but this is an NP hard problem. We solve instead the convex problem
\begin{equation}\label{MMV21}
\min J_{2,1}({\bf X})\quad\text{s.t.}\quad
\boldsymbol{\mathbb{A}} {\bf X} = \bD,
\end{equation}
which gives, under certain conditions on the model matrix $\boldsymbol{\mathbb{A}}$
\cite{eldar10,chai2014imaging}, the same solution as (\ref{MMV.NP}). In (\ref{MMV21})
$J_{2,1}$ denotes the $(2,1)$-norm 
\begin{equation}
J_{2,1}({\bf X})=\sum_{q=1}^m\|{\bf
  e}_q^T{\bf X}\|_{\ell_2},
\label{eq:Jpq}
\end{equation}
which is the $\ell_1$ norm of the vector formed by the $\ell_2$ norms
of the rows of ${\bf X}$.  Furthermore, because data
are noisy in practice, we replace the equality constraint in
(\ref{MMV21}) by $\|\boldsymbol{\mathbb{A}} {\bf X} -
\bD\|_F < \epsilon$, where $\| \cdot \|_F$ is the Frobenius norm and
$\epsilon$ is a tolerance commensurate with the noise
level of the data.

There are different algorithms for solving \eqref{MMV21} or its
reformulation for noisy data.  We use an extension of an iterative
shrinkage-thresholding algorithm, called GeLMA, proposed in
\cite{moscoso12} for matrix-vector equations. This algorithm is very
efficient for solving $\ell_1$-minimization problems, and has the
advantage that the solution does not depend on the regularization
parameter used to promote minimal support solutions, see \cite{moscoso12} for
details.
\begin{algorithm}
\begin{algorithmic}
\REQUIRE Set ${\bf X}=\vect0$, $\vect\Zc=\vect 0$,
and pick the step size $\mu$ and the regularization parameter
$\gamma$.  \REPEAT \STATE Compute the residual
$\boldsymbol{\mathcal{E}} = \bD - \boldsymbol{\mathbb{A}}
{\bf X}$ \STATE
${\bf X}\Leftarrow{\bf X} + \mu
\bA^\ast(\vect\Zc + \boldsymbol{\mathcal{E}})$ \STATE ${\bf
  e}_q^T{\bf X} \Leftarrow\operatorname{sign}(\| {\bf
  e}_q^T{\bf X}\|_{\ell_2}-\mu\gamma)\frac{\|{\bf
    e}_q^T{\bf X} \|_{\ell_2}-\mu\gamma}{\|{\bf
    e}_q^T{\bf X}\|_{\ell_2}}{\bf
  e}_q^T{\bf X}$, $~ ~q=1,\ldots,Q$ \STATE
$\vect\Zc\Leftarrow\vect\Zc + \gamma\boldsymbol{\mathcal{E}}$
\UNTIL{Convergence}
\end{algorithmic}
\caption{GeLMA-MMV}
\label{algo}
\end{algorithm}

After estimating ${\bf X}$ with Algorithm
\ref{algo}, we recover the discretized direction and frequency
dependent reflectivity using equation (\ref{eq:lem2.4}),
\begin{equation}
  \rho_q^{(\alpha,\beta)} = X_q^{(\alpha,\beta)} \exp \Big[
    2 i k_\beta \vm_\alpha \cdot \Delta \vy_q - i k_\beta \frac{\Delta
      \vy_q \cdot \mathbb{P}_\alpha \Delta \vy_q}{L_\alpha} \Big],
  \label{eq:result}
\end{equation}
for the imaging points $\vy_q = \vy_o + \Delta \vy_q$ indexed by $q =
1, \ldots, Q$, the sub-apertures indexed by $\alpha = 1, \ldots,
\cN_\alpha$ and frequency sub-bands indexed by $\beta = 1, \ldots,
\cN_\beta$.

\section{Numerical simulations}
\label{sect:NUM}
We begin in section \ref{sect:NumSet} with the numerical setup, which
is in the regime of the GOTCHA Volumetric data set
\cite{gotcha} for X-band persistent surveillance SAR. Then we present
in sections \ref{sect:results1D} and \ref{sect:results2D} the simulation results.

\subsection{Imaging in the X-band (GOTCHA) SAR regime}
\label{sect:NumSet}
The numerical simulations generate the data with the model
(\ref{eq:F2}), for various scattering scenes. The regime of parameters
is that of the GOTCHA data set, where the platform trajectory is circular at
height $H=7.3$km, with radius $R=7.1$km and speed $V=70$m/s.  The
signal $f(t)$ is sent every $1.05$m along the trajectory, which gives
a slow time spacing
$h_s=0.015$s. The carrier frequency is $\om_o/(2 \pi) = 9.6$GHz and
the bandwidth is $B=622$MHz. The waves propagate at electromagnetic
speed $c = 3 \cdot 10^8$m/s, so the wavelength is $\la_o =
3.12$cm. The image window $\mathcal{Y}$ is at the ground level,
below the center of the flight trajectory, and the distance from the
platform to its center $\vy_o$ is $L = 10.18$km. It is 
a square, with side length $Y = Y^\perp$ of the order
of $40$m.  The size of the sub-apertures is $a = 42$m and the width of each
sub-band is $b = B/15$.

Given these parameters, the nominal resolution limits are 
\[
  \la_o L/a =7.56\mbox{m}, \qquad c/b = 7.23\mbox{m}.
\]
The image window $\mathcal{Y}$ is discretized in uniform steps
$h = 2$m in range and $h^\perp = 1$m in cross-range, and the reflectivity
is modeled as piecewise constant on the imaging grid. {The 
image discretization affects the quality of the reconstruction with $\ell_1$ optimization.
It must be coarse enough so that uniqueness of the $\ell_1$ minimizer holds, and yet 
fine enough so that modeling errors due to off-grid placement of the unknown
are controlled. We refer to \cite{borcea2015resolution} for a study of this trade-off.}

{To illustrate the performance of the algorithm, we present in the next 
two sections results for various imaging scenes consisting of  small scatterers 
supported on one pixel of the imaging grid, or over multiple adjacent pixels. The latter 
is for representing larger scatterers for which the direction dependent reflectivity can be motivated 
by Snell's law of reflection at their surface. }

The results presented in the next sections compare the images obtained with 
reverse time migration and the algorithm proposed in this paper,
hereby referred to as the MMV algorithm.  The migration image is
computed with the formula
\begin{equation}
  \cI(\vy) = \frac{(4 \pi)^2}{k_o^2 |\hat f(\om_o)|^2 n_s n_\om h
    h^\perp} \sum_{j = 1}^{n_s} \sum_{l=1}^{n_\om} d(s_j,\om_l)
  |\vr(s_j)-\vy|^2 e^{-2 i \om_l
      \big[\tau(s_j,\vy)-\tau(s_j,\vy_o)\big]},
  \label{eq:MIGN}
\end{equation}
which is a weighted version of (\ref{eq:migration}), where the weights
are chosen so as to provide a quantitative estimate of the unknown $\rho$. That
is to say, when we substitute the data model in (\ref{eq:MIGN}), under
the assumption of an isotropic and frequency independent reflectivity
we get that $\cI(\vy)$ peaks at the true location of the scatterers and its
value at the peaks equals the true reflectivity there.

Let us verify the assumptions (\ref{eq:M5})-(\ref{eq:M10}) with the
GOTCHA parameters.  The Fresnel numbers are larger than one, as stated in
(\ref{eq:M5}),
\[
  \frac{a^2}{\la_o L} = 5.55 \quad \mbox{and} \quad
  \frac{(Y^\perp)^2}{\la_o L} = 5.04.
\]
The size of the imaging region and the range resolution satisfy
(\ref{eq:M6}). Moreover,
\[
\frac{b}{\om_o} \frac{Y^\perp}{\la_o L/a} = 0.0036,
\]
which is consistent with (\ref{eq:M8}), and (\ref{eq:M10}) is
satisfied as well,
\[
\frac{a^2 Y^\perp}{\la_o L^2} = \frac{a^2 Y}{\la_o L^2} = 0.022.
\]

\begin{figure}[t]
  \begin{center}
    \includegraphics[width=8.cm,height=9.cm]{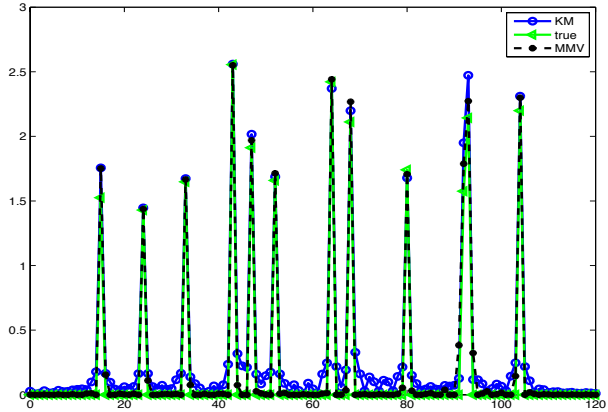}
    \end{center}
 \vspace{-0.5in}
\caption{Estimation of an isotropic, frequency independent
  reflectivity as a function of cross-range, using $N_\alpha = 8$
  consecutive, non-overlapping apertures. The exact reflectivity is
  shown with the full green line, the migration result with the blue
  line and the MMV inversion result with the broken line. The abscissa
  is cross-range in meters.}
\label{fig:isotropic}
\end{figure}

\subsection{Single frequency results}
\label{sect:results1D}
We begin with imaging results at the carrier frequency, where we
assume we know the range of the scatterers and seek to reconstruct
their reflectivity as a function of cross-range and direction.  The
image window extends over $120$m in cross-range, and it is sampled
in steps $h^\perp = 1$m, where we recall that $\la_o L/a = 7.56$m.

\begin{figure}[t]
\begin{minipage}{1.1\textwidth}
  \begin{center}
   \hspace{-1cm} \includegraphics[width=5.3cm,height=5.3cm]{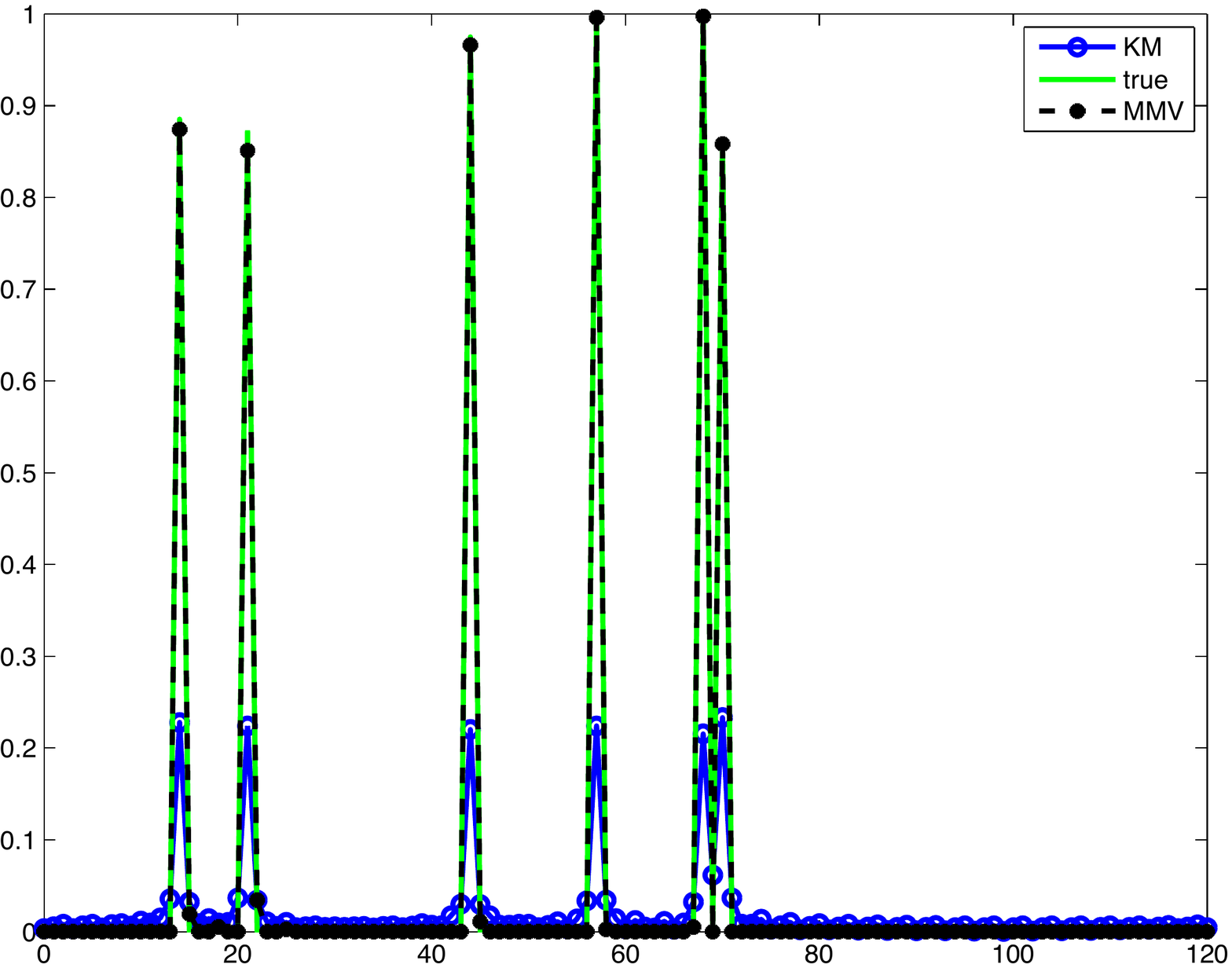}
\includegraphics[width=5.3cm,height=5.3cm]{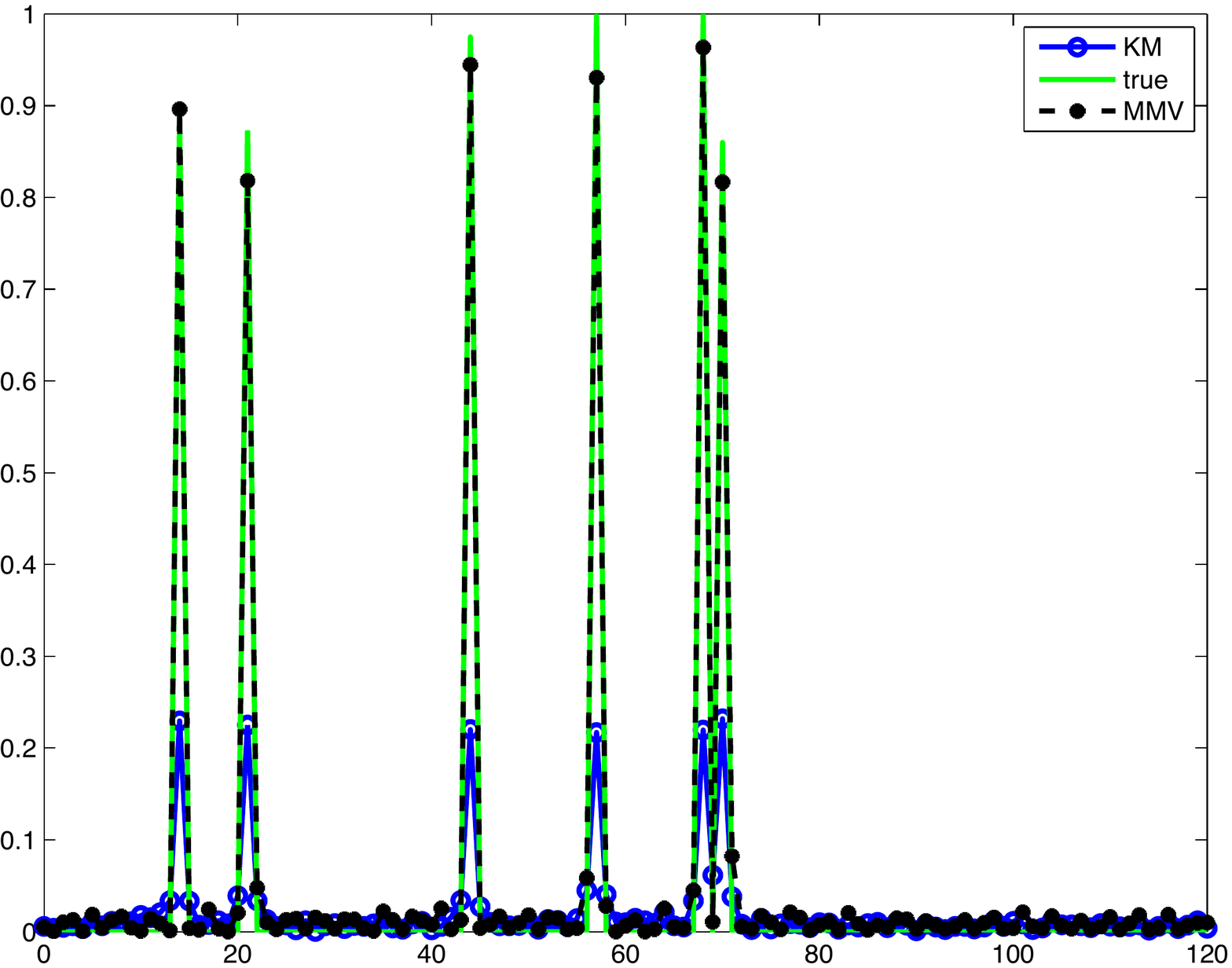}
\\ 
\vspace{-0.6in}
\hspace{-1cm}  \includegraphics[width=4.2cm]{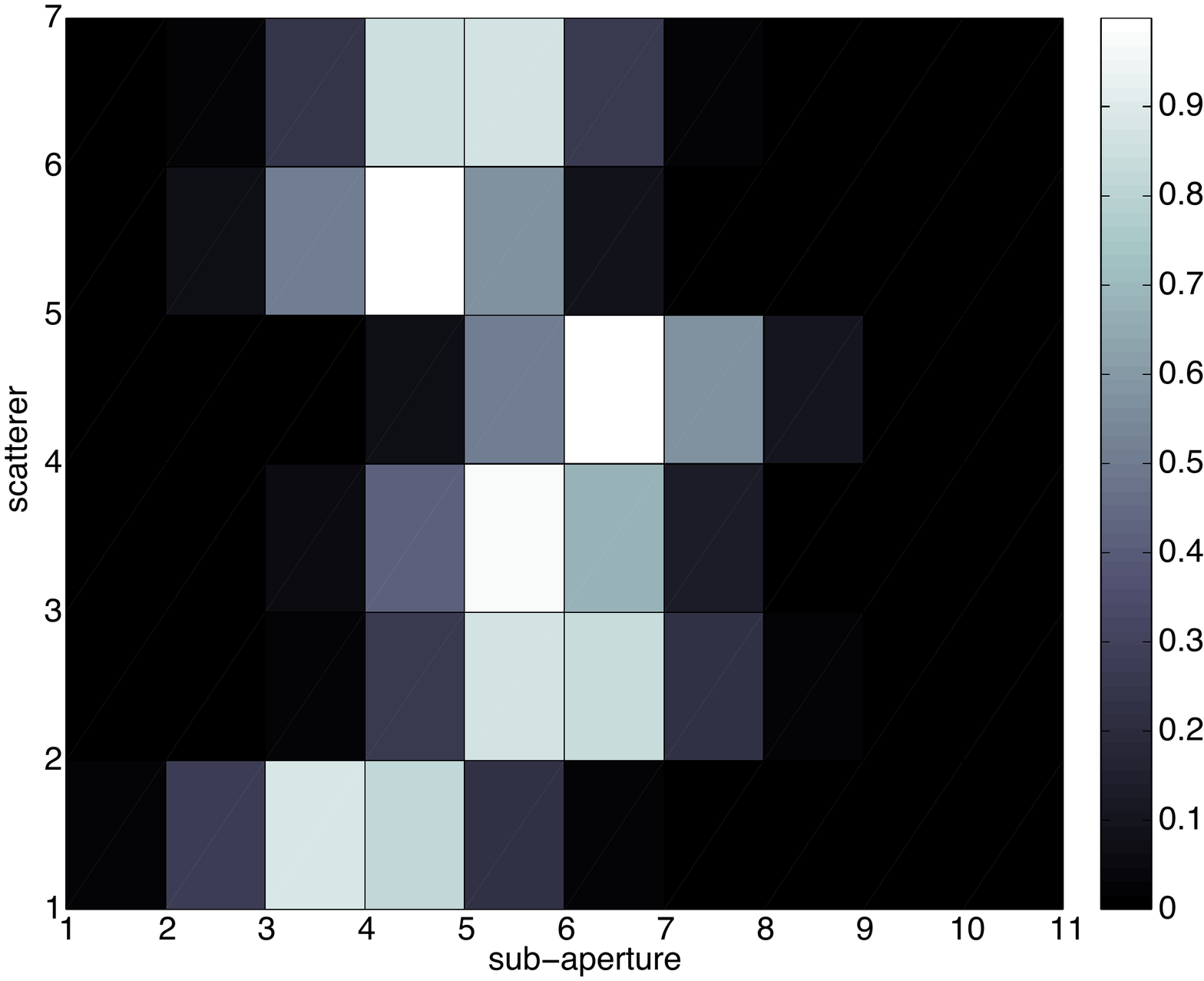}
\includegraphics[width=4.2cm]{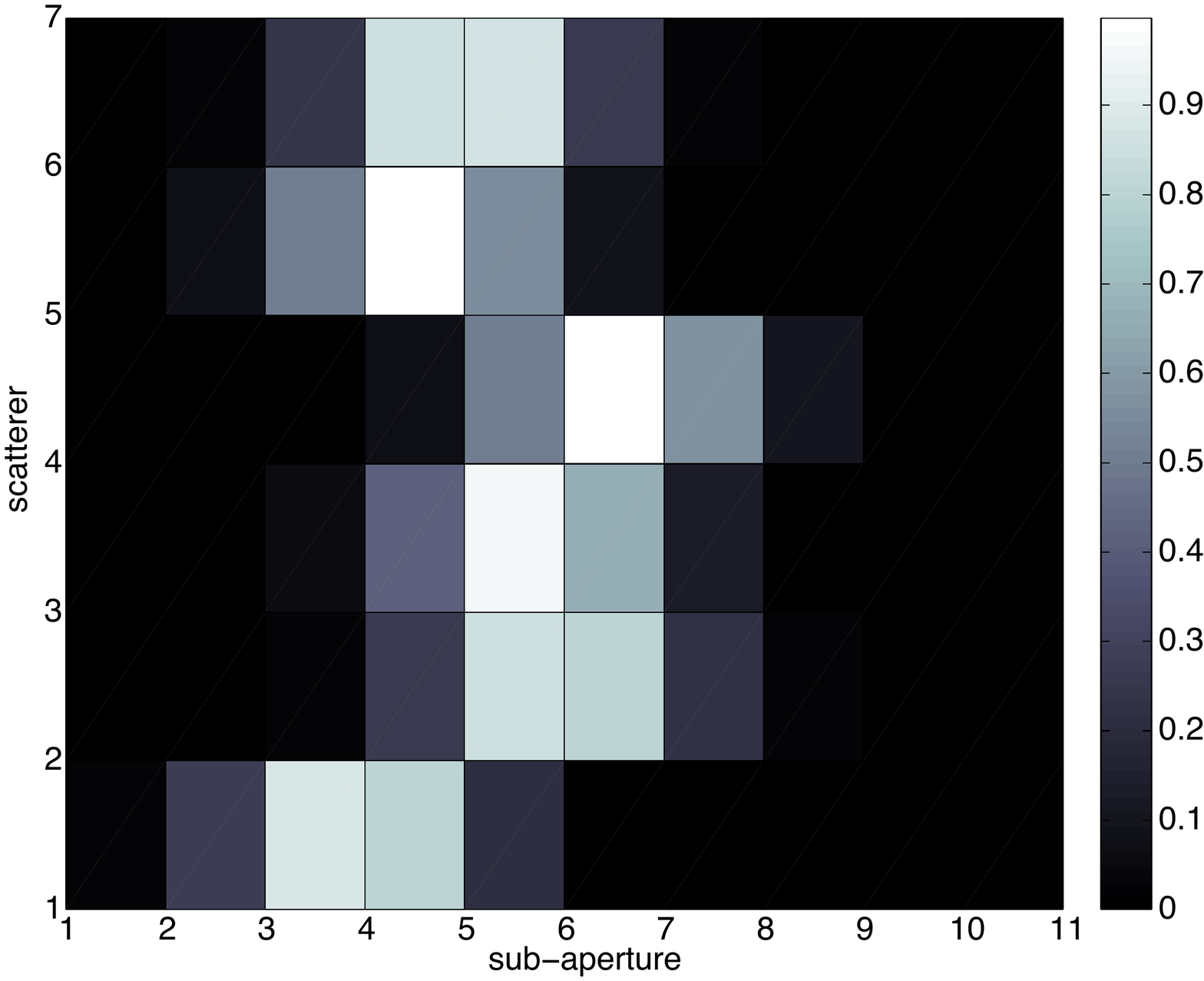}
\includegraphics[width=4.2cm]{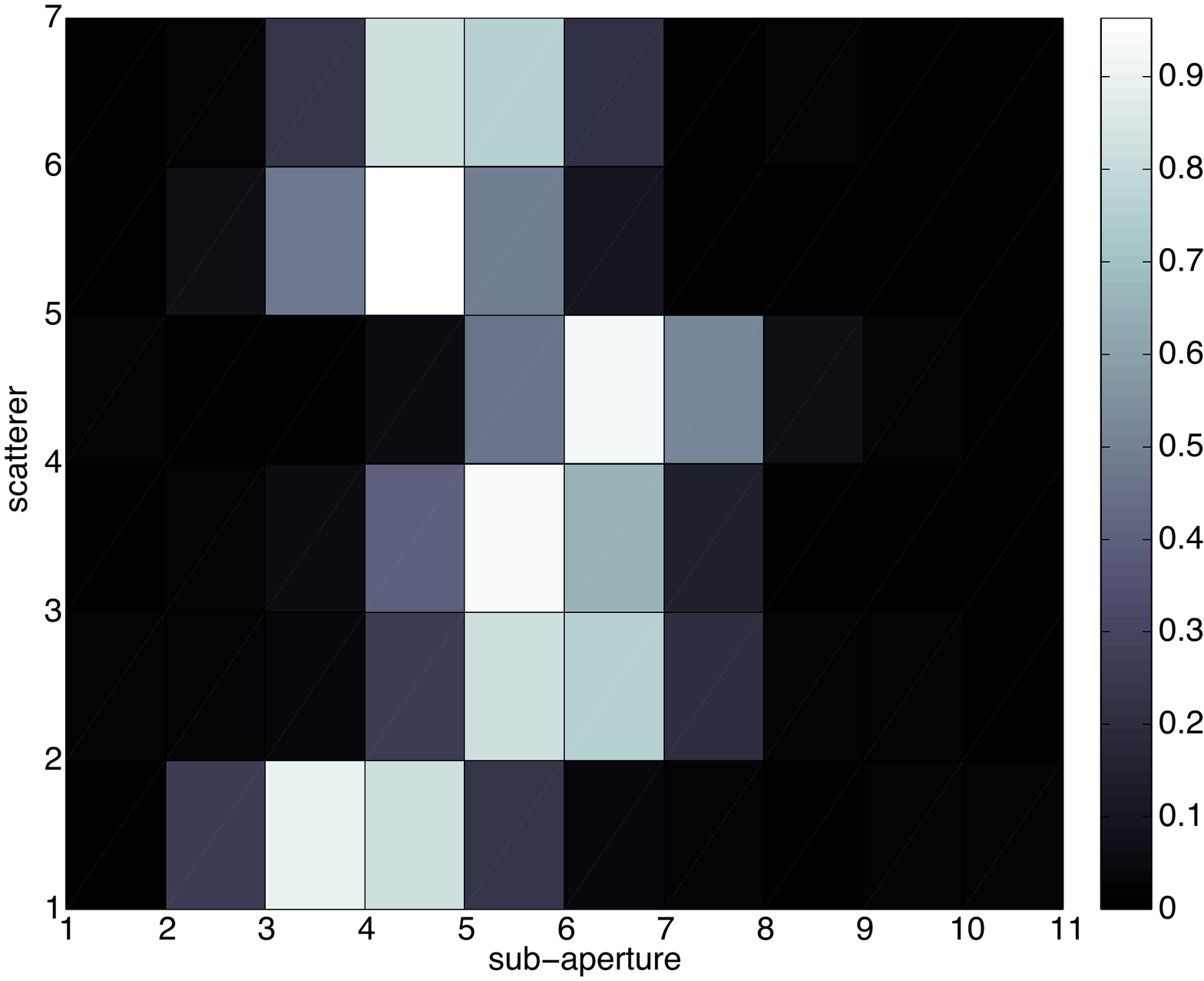}
  \end{center}
\end{minipage}
\vspace{-0.5in}
\caption{Estimation of the reflectivity as a function of direction and
  cross-range location for a scene with $6$ scatterers. The top plots
  show the reflectivity as a function of cross-range (the abscissa in
  meters), for the peak directions. The left plot is for noiseless
  data and the right plot is for data contaminated with $10\%$
  additive noise.  The green line is the exact peak value and the broken
  line the one obtained with MMV. The blue line is obtained with migration.
  The bottom plots display the reflectivity of each
  scatterer as a function of sub-aperture i.e., the slow time index
  $\alpha = 1, \ldots, 10$, where $10$ is the number of sub-apertures.
  The left plot is for the true reflectivity, the middle plot is for
  the noiseless reconstruction and the right plot is for the noisy
  reconstruction.}
\label{fig:1d_real1}
\end{figure}  

The first result displayed in Figure \ref{fig:isotropic} is for an
isotropic, frequency independent reflectivity of $11$ scatterers,
$\cN_\alpha = 8$ consecutive, non-overlapping apertures and noiseless
data. We display in green the true reflectivity, in blue the
reflectivity estimated with formula (\ref{eq:MIGN}), and with broken
line the result of the MMV inversion algorithm. In the legend we
abbreviate the migration formula result with the letters KM, standing
for Kirchhoff Migration. The figure shows that the MMV algorithm
reconstructs exactly the reflectivity, and that the weighted migration
formula (\ref{eq:MIGN}) does indeed give quantitative estimates of the
reflectivity. However, the migration estimates deteriorate when the
reflectivity is anisotropic and frequency dependent, as illustrated
next.

\begin{figure}[t]
\begin{minipage}{1.1\textwidth}
\includegraphics[width=5.3cm,height=5.3cm]{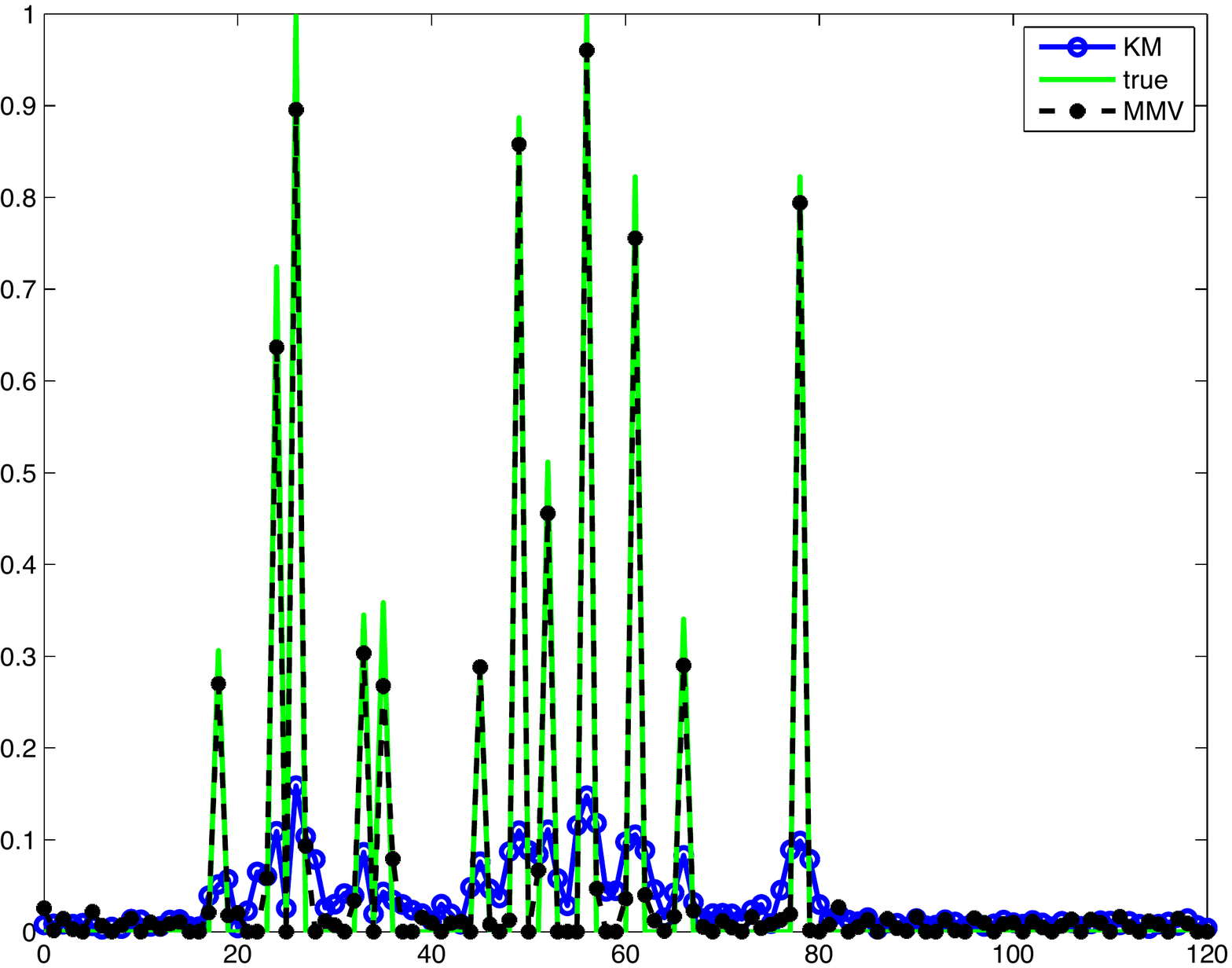} 
\includegraphics[width=4cm]{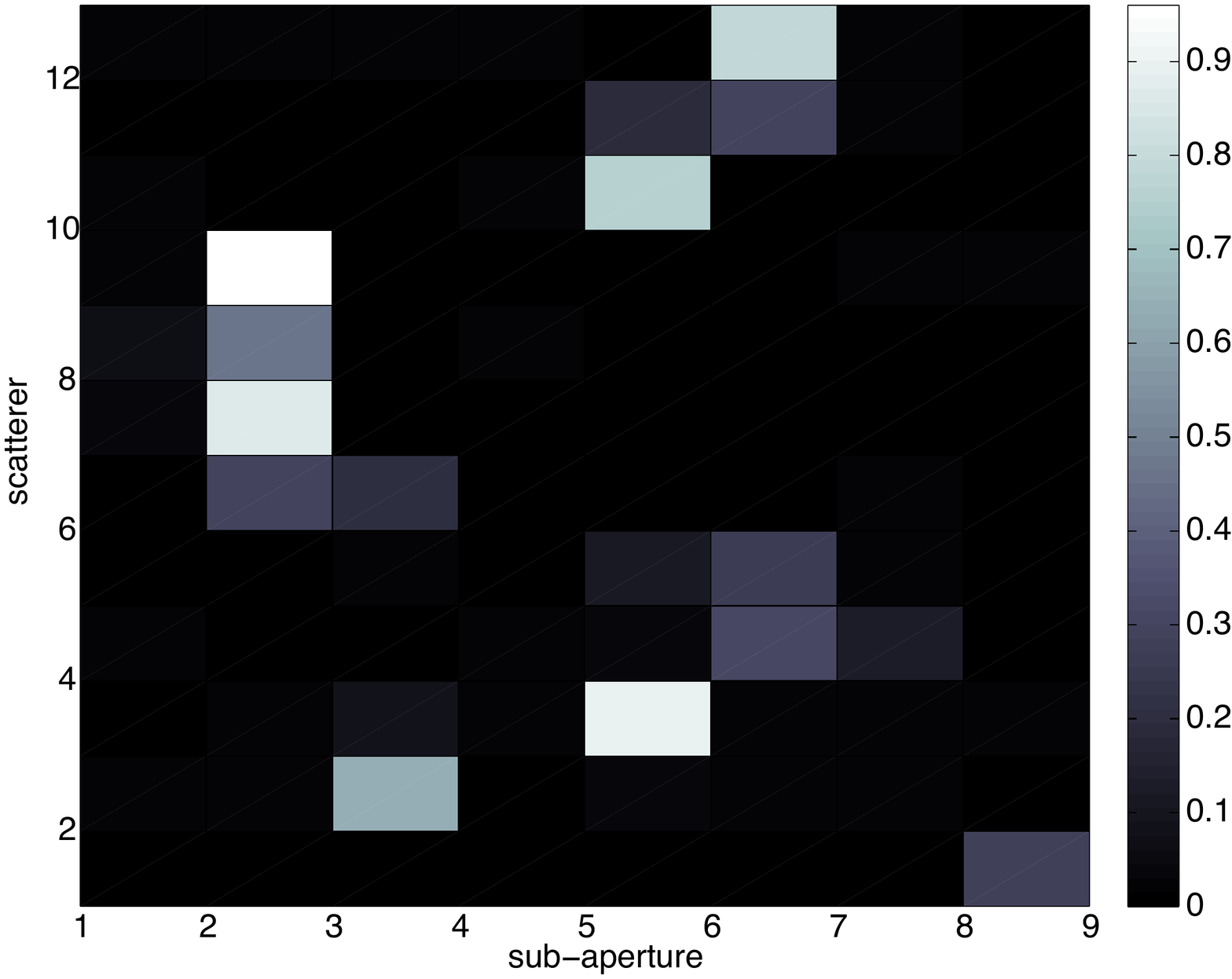}
\includegraphics[width=4cm]{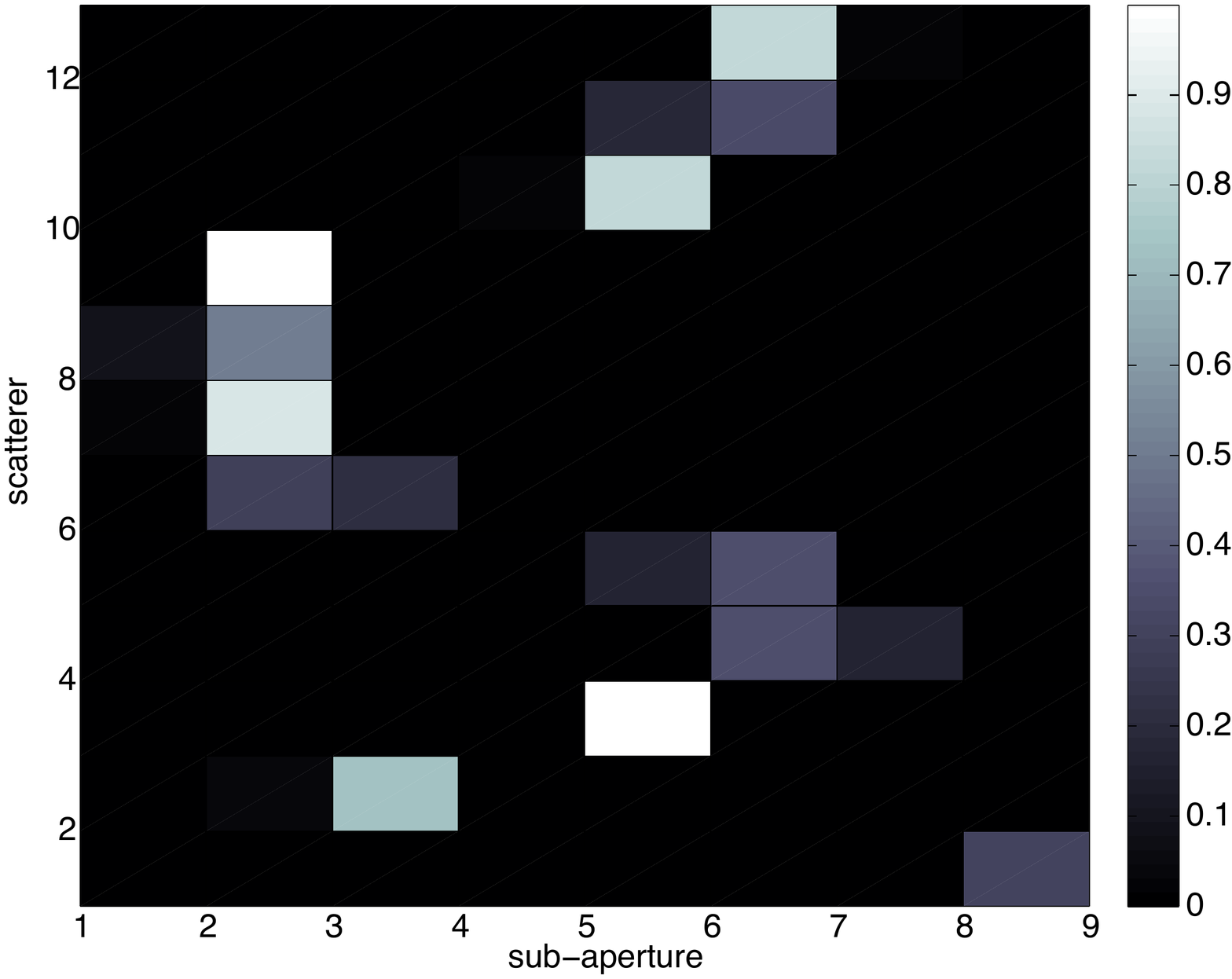}
\vspace{-0.6in}
\\ \vspace{-0.6in}
\includegraphics[width=5.3cm,height=5.3cm]{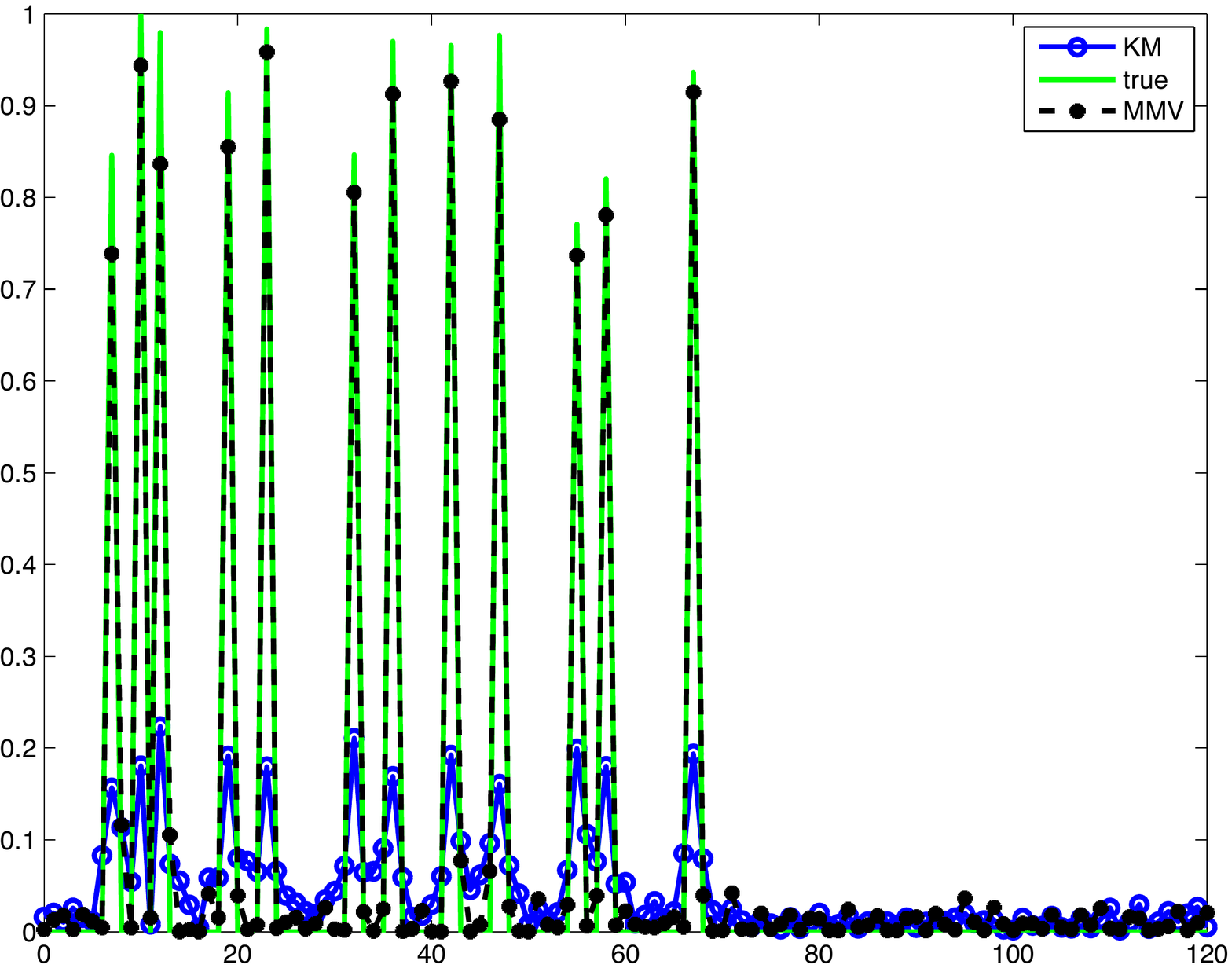} 
\includegraphics[width=4cm]{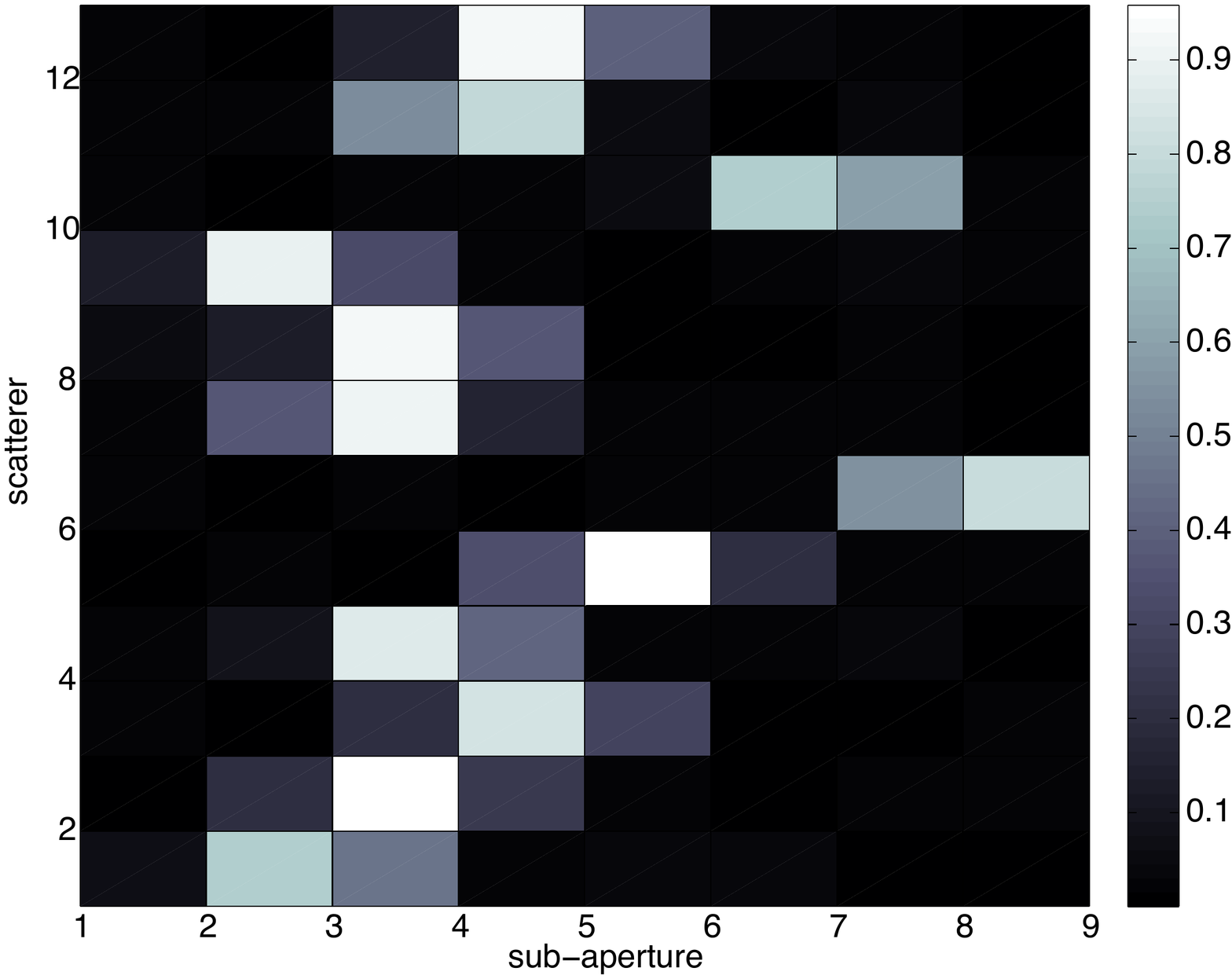}
\includegraphics[width=4cm]{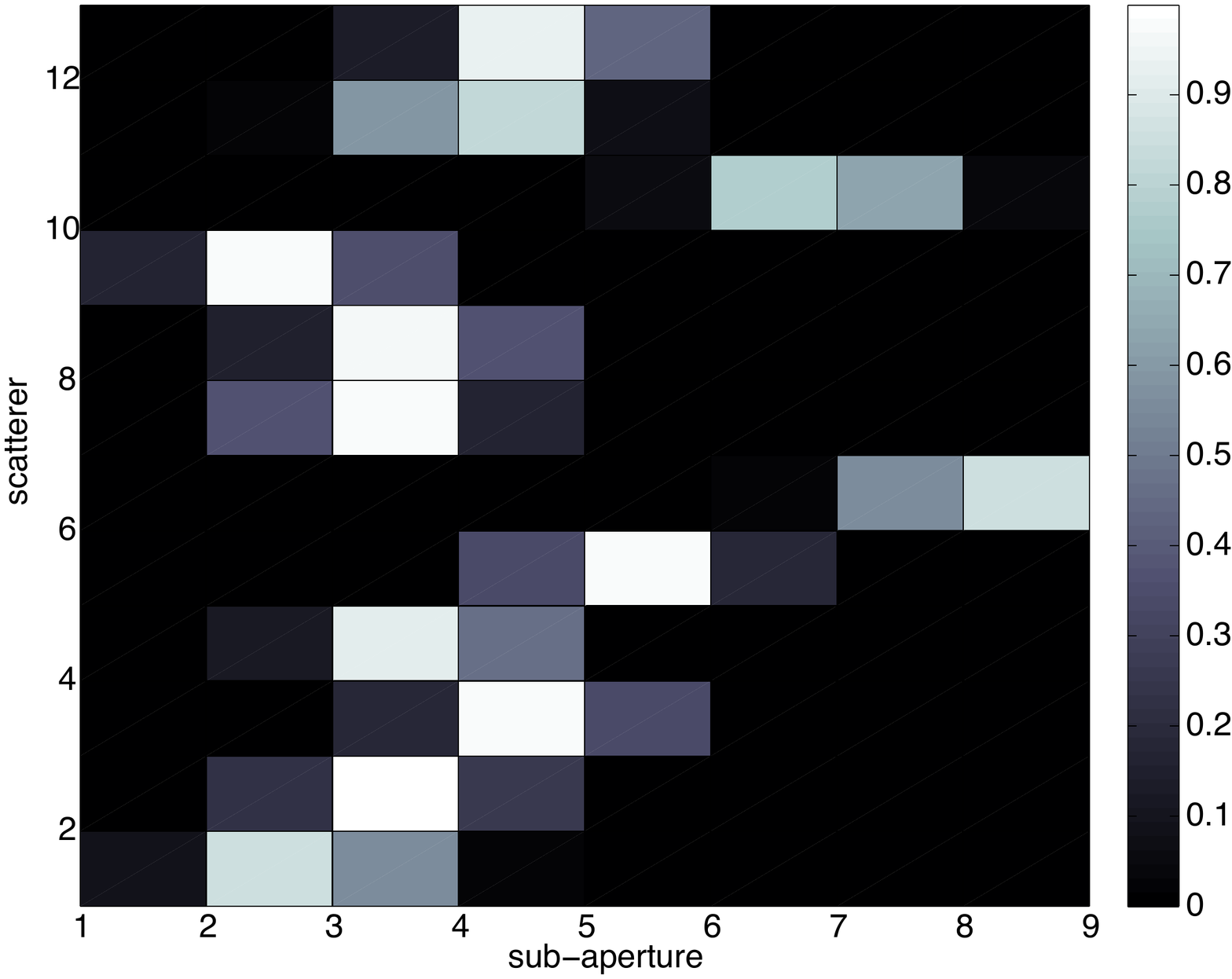}
\\ \vspace{-0.4in} \includegraphics[width=5.3cm,height=5.3cm]{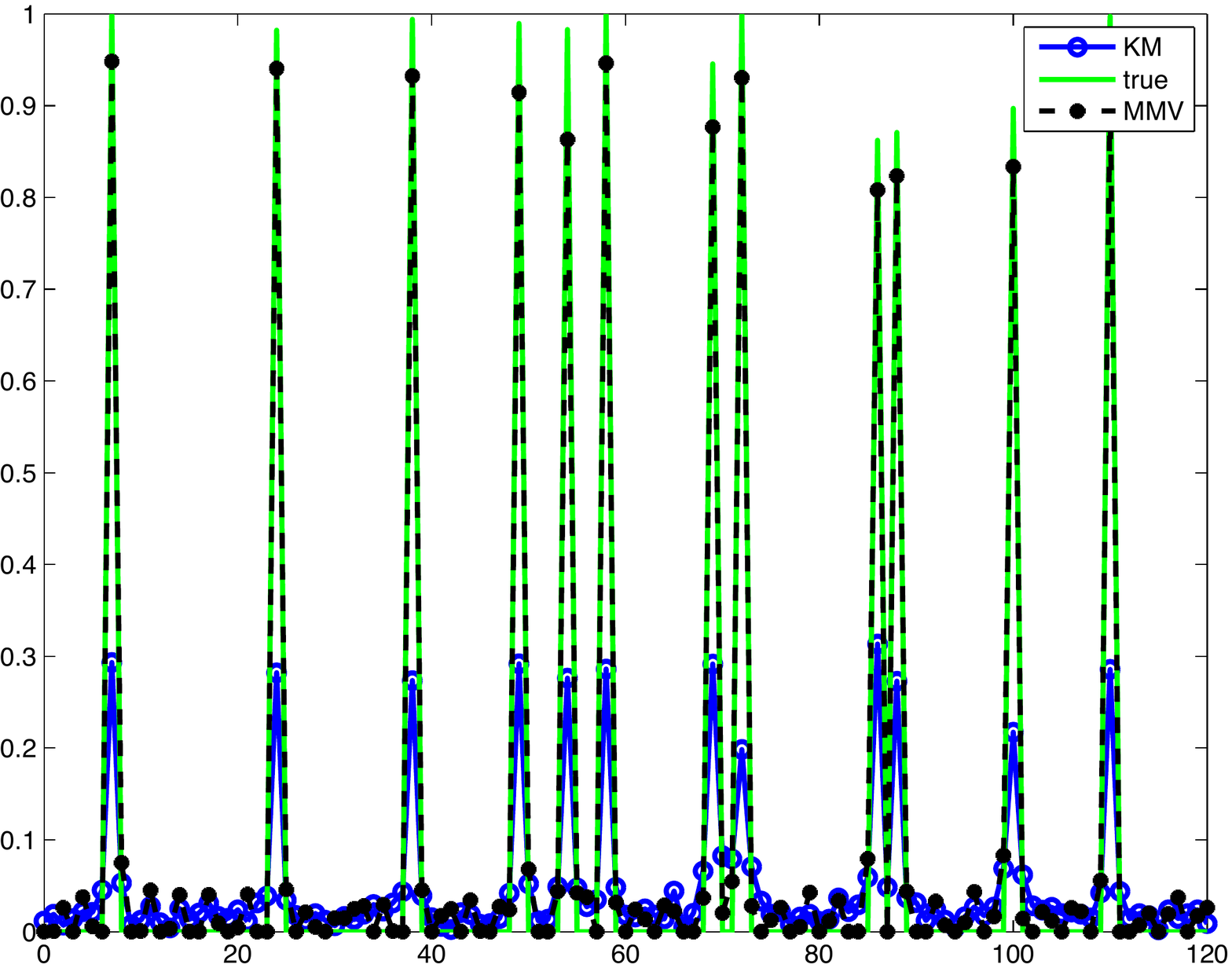} 
\includegraphics[width=4cm]{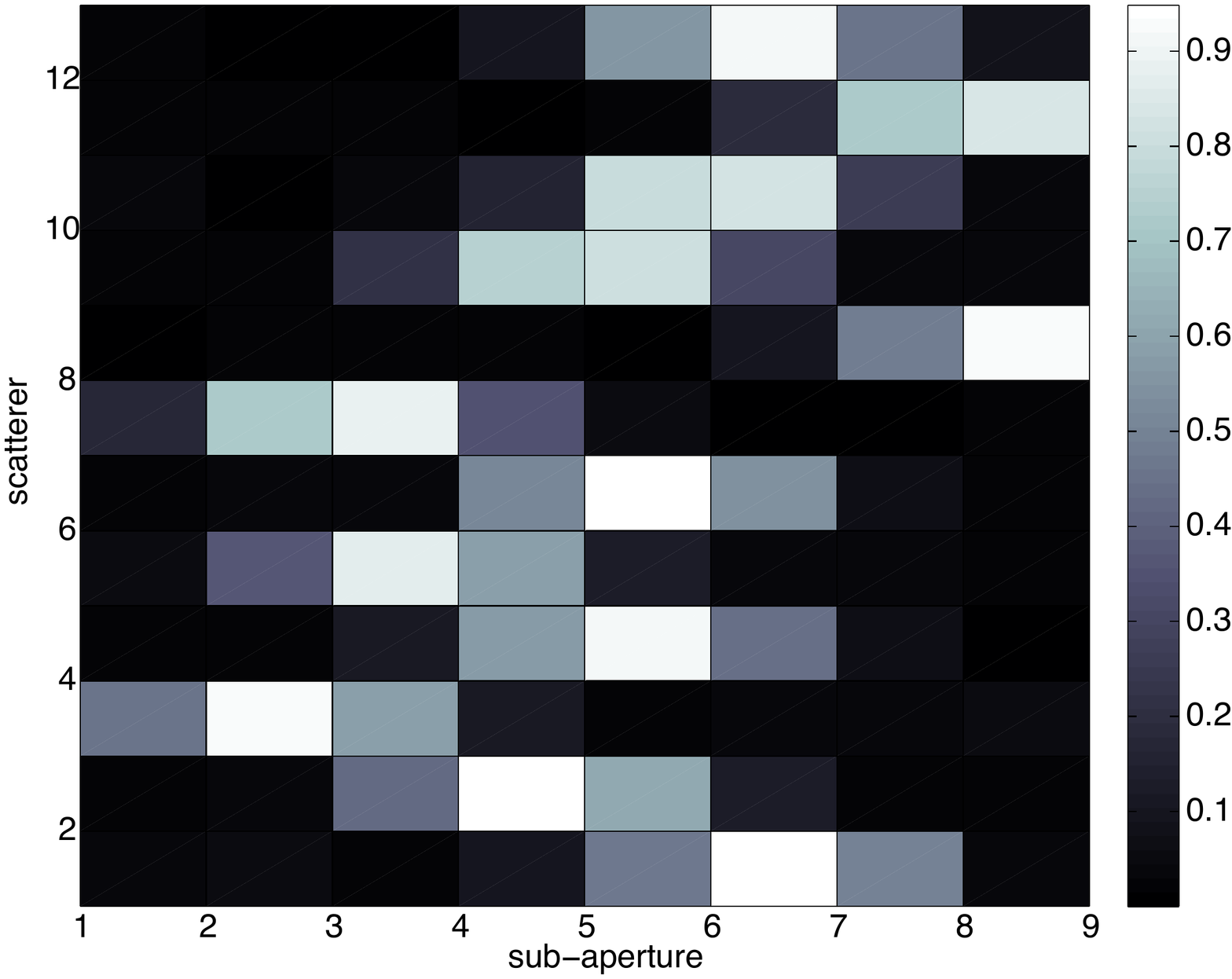}
\includegraphics[width=4cm]{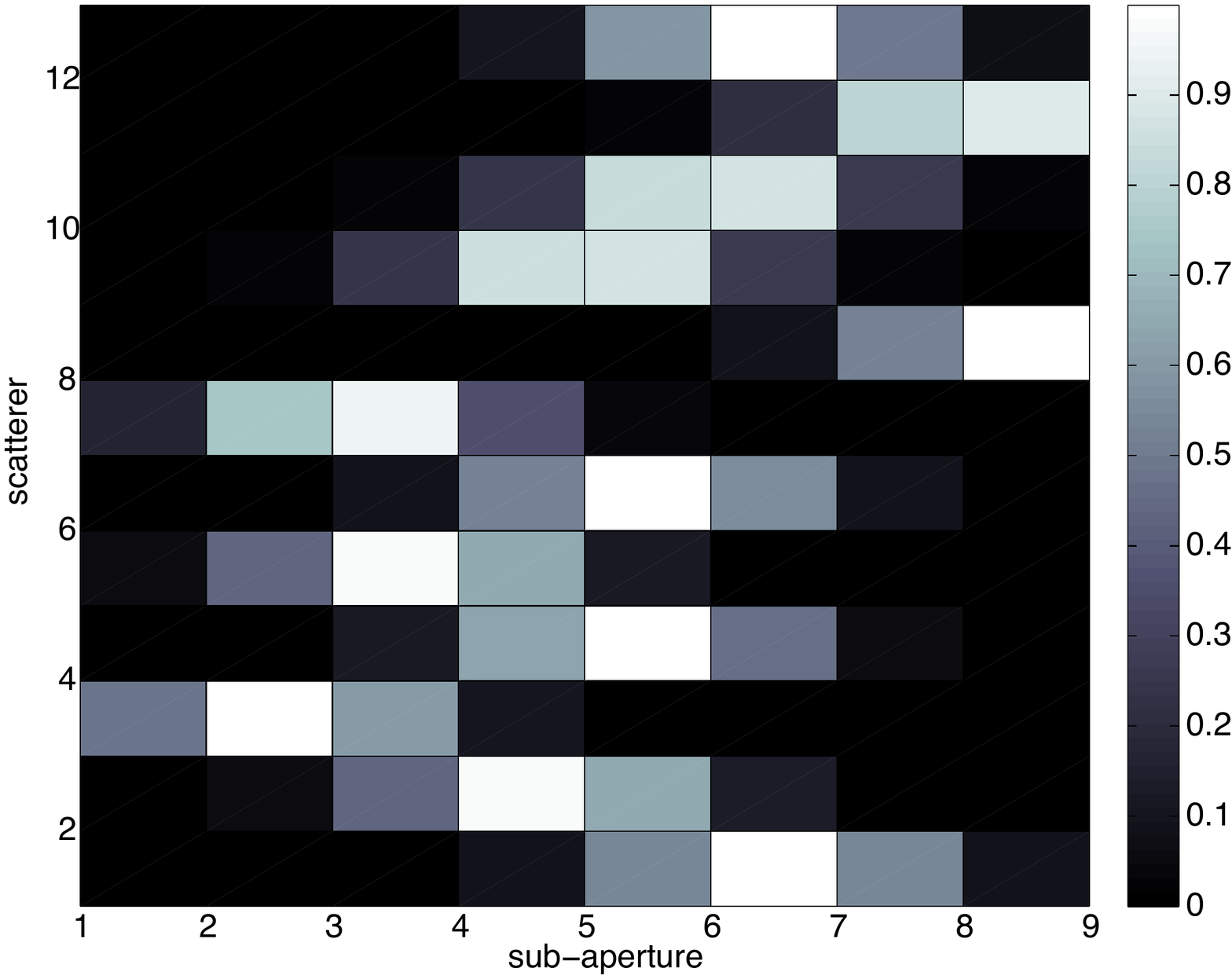}
\end{minipage}
\caption{Imaging results for anisotropic reflectivities, $\cN_\alpha =
  10$ sub-apertures and data contaminated with $10\%$ additive
  noise. From top to bottom we decrease the anisotropy. This can be
  seen from the right column plots which show the reflectivity of each
  scatterer for each sub-aperture. The middle column shows the
  reconstructed reflectivity as a function of
  direction with the MMV algorithm.}
\label{fig:1dapertures_real1}
\end{figure}

The results displayed in Figure \ref{fig:1d_real1} are 
obtained with $\cN_\alpha = 10$ consecutive, non-overlapping
apertures. The reflectivity depends on two variables: the cross-range
location and the scattering direction, parameterized by the slow time
$s_\alpha^\star$, for $\alpha = 1, \ldots, 10$. In discretized form it
gives a matrix $\boldsymbol{\mathcal{R}}_{\mbox{true}}$ with row
index corresponding to the pixel location in the image window, and
column index corresponding to the sub-aperture.
The reconstruction of this matrix is
denoted by ${\boldsymbol{\mathcal{R}}}$. The green and broken lines in
the top plots in the figure display the true and reconstructed
reflectivity at the peak direction, vs. cross-range. Explicitly, for
each pixel in the image i.e., each row $q$ in
$\boldsymbol{\mathcal{R}}_{\mbox{true}}$ or
${\boldsymbol{\mathcal{R}}}$, we display the maximal entry.  The
migration image of the reflectors is independent of the direction and
is plotted with the blue line.  The results show that we have $6$
small scatterers, which are well estimated by the MMV algorithm even
for noisy data. The migration method identifies correctly the
locations of the $6$ scatterers, but the reflectivity value is no longer
accurate because only a few sub-apertures see each
reflector, as we infer from the bottom plots described next. This also
implies a deterioration in the cross-range resolution which is more
visible in the next set of results in Figure
\ref{fig:1dapertures_real1}. Naturally, the migration image gives no
information about the direction dependence of the reflectivity.

In the bottom plots in Figure \ref{fig:1d_real1} we show the value of
the reflectivity of each scatterer as a function of direction,
parameterized by the slow time $s_\alpha^\star$. That is to say, we
identify first the row indexes $q$ in
$\boldsymbol{\mathcal{R}}_{\mbox{true}}$ or
${\boldsymbol{\mathcal{R}}}$ at which we have a strong
scatterer (see top plots) and then display those rows. The left plot
is for the true reflectivity, the middle is for the noiseless
reconstruction, and the right is for the noisy reconstruction. We
observe that the MMV method reconstructs the direction dependent
reflectivity exactly in the noiseless case, and very well in the noisy
case.

In Figure \ref{fig:1dapertures_real1} we illustrate the effect of the
anisotropy of the reflectivity on the imaging process. We display the
results the same way as in in the previous figure. The point is to
notice that while the MMV method estimates accurately the direction
dependent reflectivity in all cases, the migration method performs
poorly when the anisotropy is strong, meaning that each scatterer is
seen only by one sub-aperture at a time (top plots). The resolution is
not that corresponding to the actual aperture of $10a = 420$m, but
that for a single sub-aperture of $a = 42$m. The middle and bottom row
plots show how migration images improve when the anisotropy of the
reflectivity is weaker and more sub-apertures see each scatterer.

\subsection{Multiple frequency results}
\label{sect:results2D}
{Now we consider multiple frequency sub-bands and thus seek to estimate the
reflectivity as a function of range, cross-range, direction and
frequency.} We have $\cN_\om$ sub-bands of width $b$, and we sample
each of them at $n_\om = 15$ frequencies. The number of sub-apertures is
$\cN_\alpha = 8$. The imaging region is a square of side $40$m and it
is sampled in cross-range in steps $h^\perp = 1$m and in range in
steps $h = 2$m.  We denote, as before, by
$\boldsymbol{\mathcal{R}}_{\mbox{true}}$ the true matrix of
discretized reflectivities and by $\boldsymbol{\mathcal{R}}$ the
reconstructed ones. These are matrices of size $Q \times \cN_\alpha \cN_\om$ and we
display them in the image window $\mathcal{Y}$ as follows: For each
pixel in the image window i.e., a row $q$ in
$\boldsymbol{\mathcal{R}}_{\mbox{true}}$ or
$\boldsymbol{\mathcal{R}}$, we display the maximum entry, the peak
value of the reflectivity at point $\vy_q$ over directions and
frequencies. Once we identify the location of the scatterers from these
images, i.e., determine their associated rows, we display the entries
in these rows, to illustrate the direction and frequency dependence of
their reflectivity. These are the middle and right plots in the
figures.

\begin{figure}[t]
\begin{minipage}{1.1\textwidth}
\hspace*{4cm}
\includegraphics[width=4.2cm]{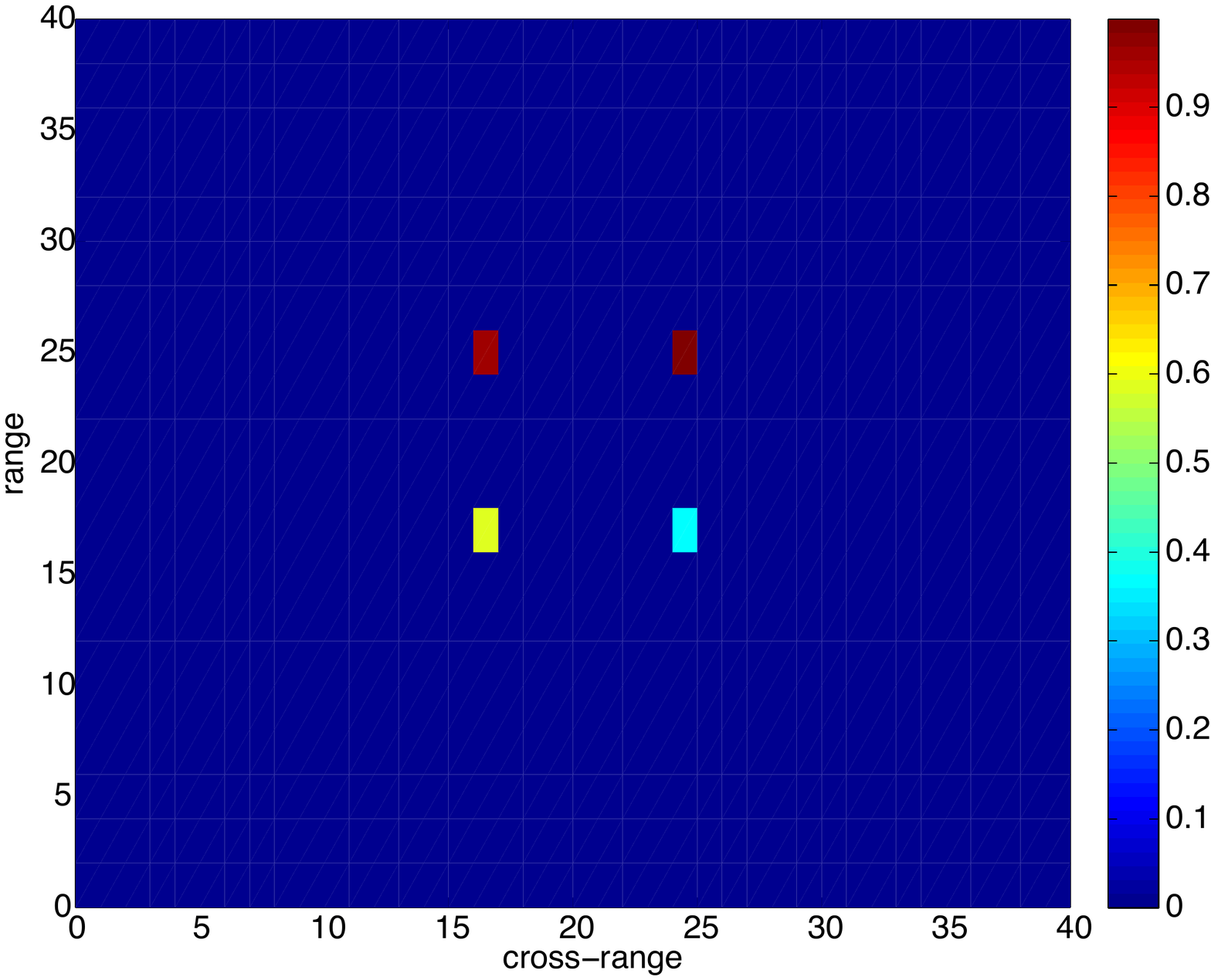}
\includegraphics[width=4.2cm]{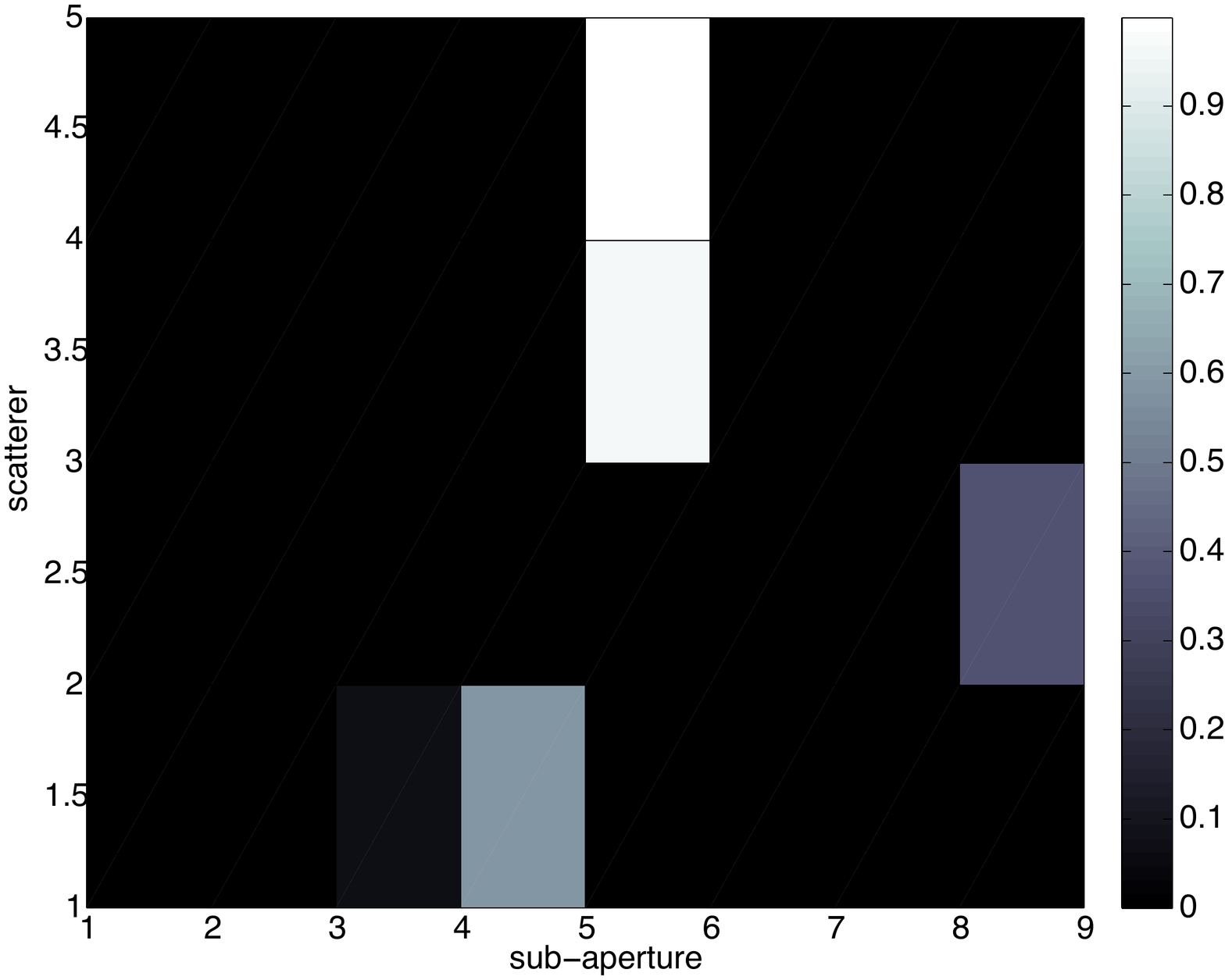}
\vspace{-0.6in}
\\
\includegraphics[width=4.2cm]{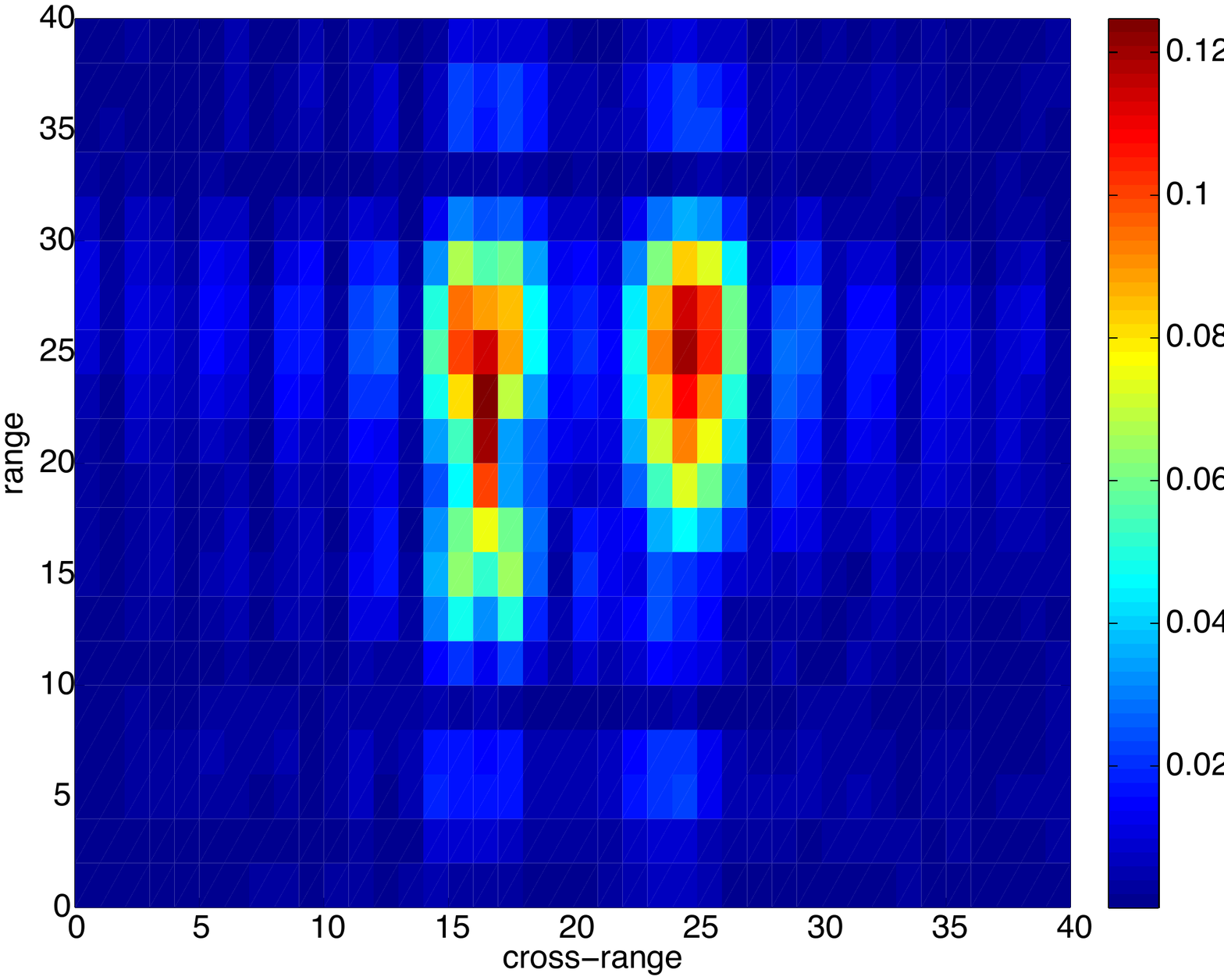} 
\includegraphics[width=4.2cm]{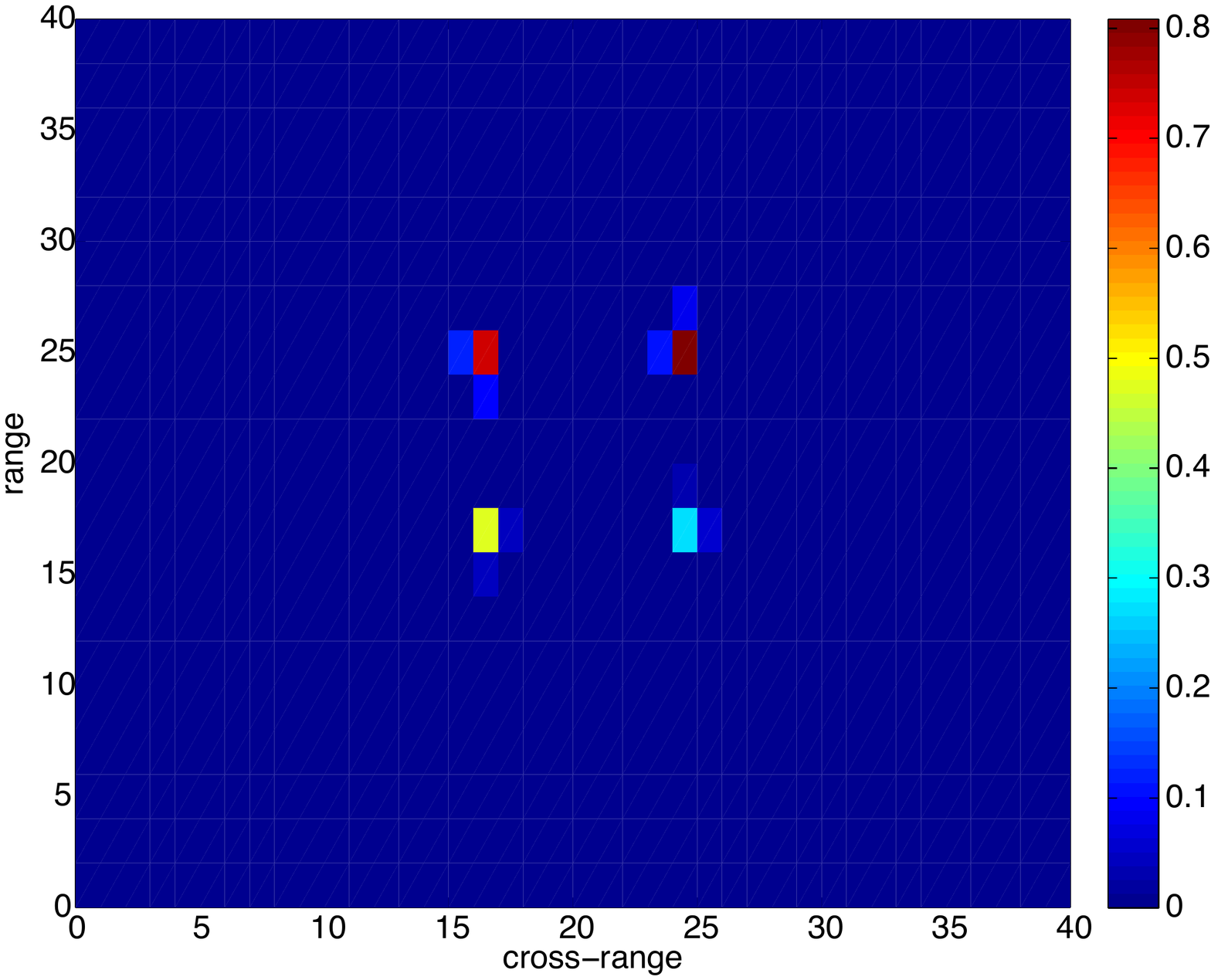} 
\includegraphics[width=4.2cm]{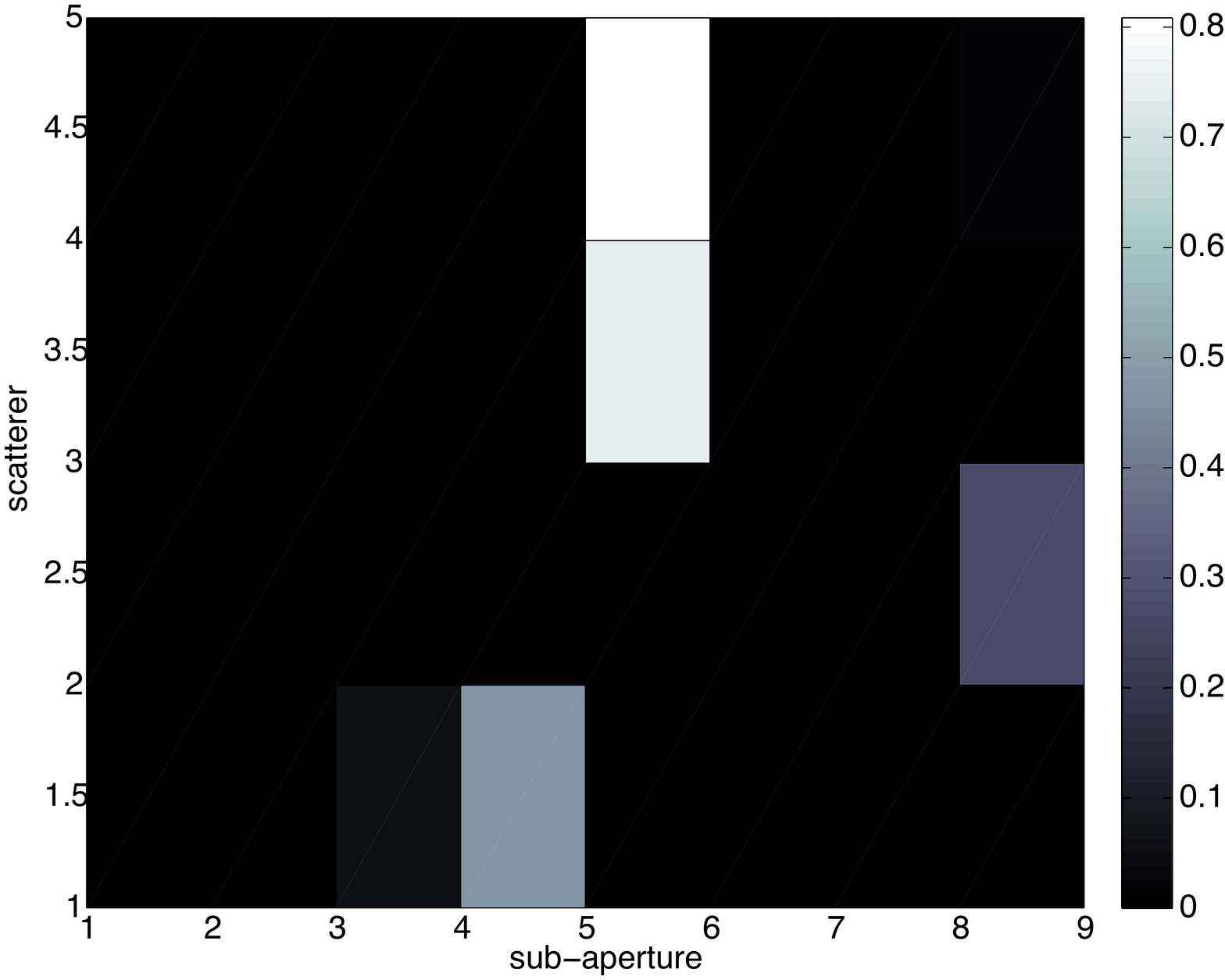}
\\
\end{minipage}
\vspace{-0.6in}
\caption{Results for a single frequency sub-band and $\cN_\alpha=8$
  consecutive, non-overlapping sub-apertures. On the top we show the
  true reflectivity as a function of location (middle) and direction (right).
  On the bottom we show the reconstructed reflectivity with
  20\% additive noise.  Left plot is the migration image, middle plot
  is the MMV image and the right plot is the directional dependence of
  the reflectivity reconstructed with the MMV algorithm. The axes in
  the left images are cross-range and range in meters. The abscissa in
  the right plots are sub-aperture index $\alpha = 1, \ldots,
  \cN_\alpha$ and the ordinate is the index of the scatterer (from $1$
  to $4$).}
\label{fig:2d_real2}
\end{figure}  
We begin in Figure \ref{fig:2d_real2} with a single frequency sub-band
($\cN_\om = 1$), $\cN_\alpha = 8$ consecutive, non-overlapping
sub-apertures and data contaminated with $20\%$ additive noise.  The
anisotropic reflectivity model has four scatterers, as illustrated in the
top plots. Each scatterer is seen by a single sub-aperture. The
reconstructed reflectivity is shown in the bottom plots. On the left
we show the migration image, which is blurry and is unable to locate the
weaker scatterers. The MMV algorithm gives an excellent reconstruction
as shown in the middle and right plots.

\begin{figure}[t]
\begin{minipage}{1.1\textwidth}
\vspace{-0.6in}
\includegraphics[width=4.2cm]{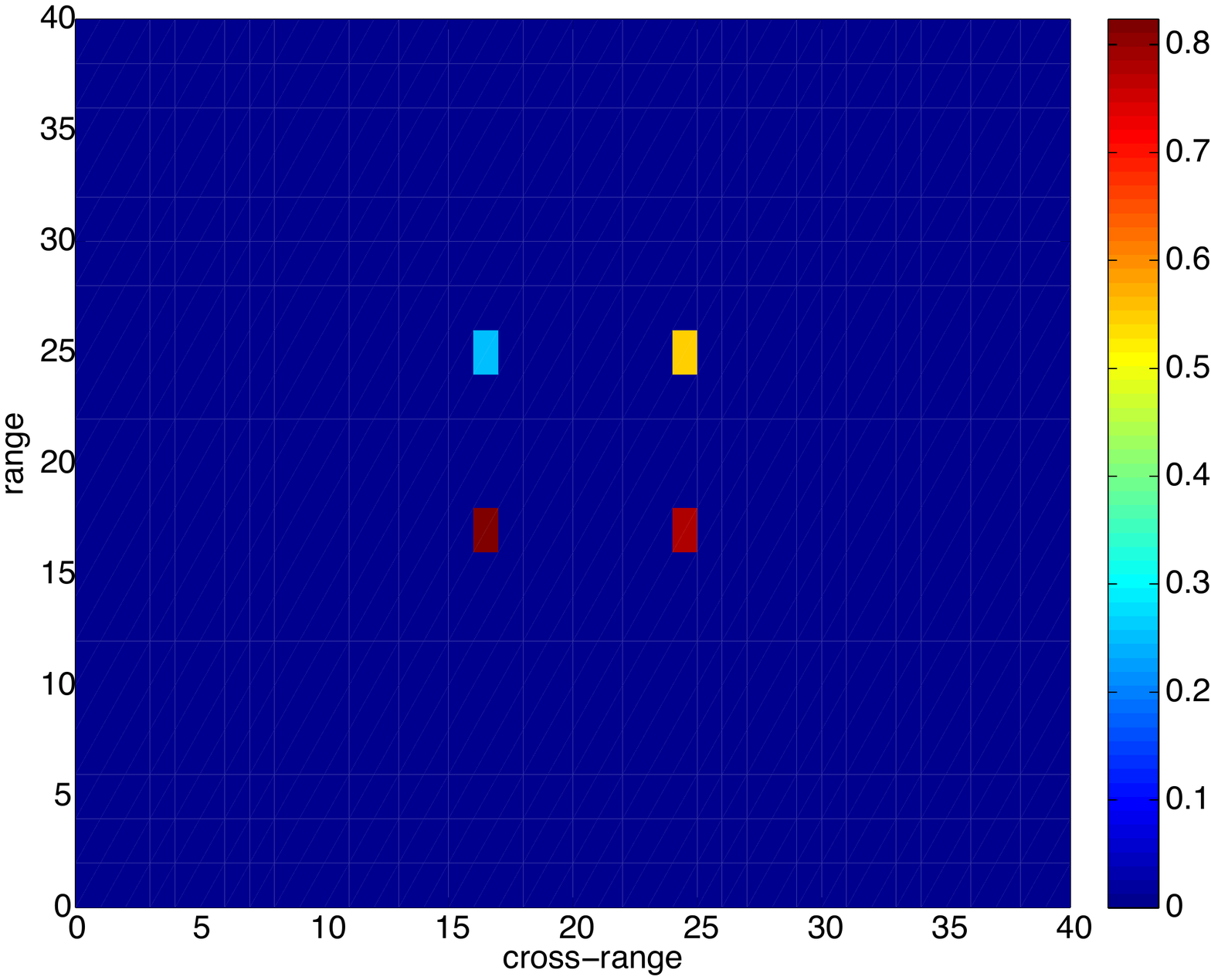}
\includegraphics[width=4.2cm]{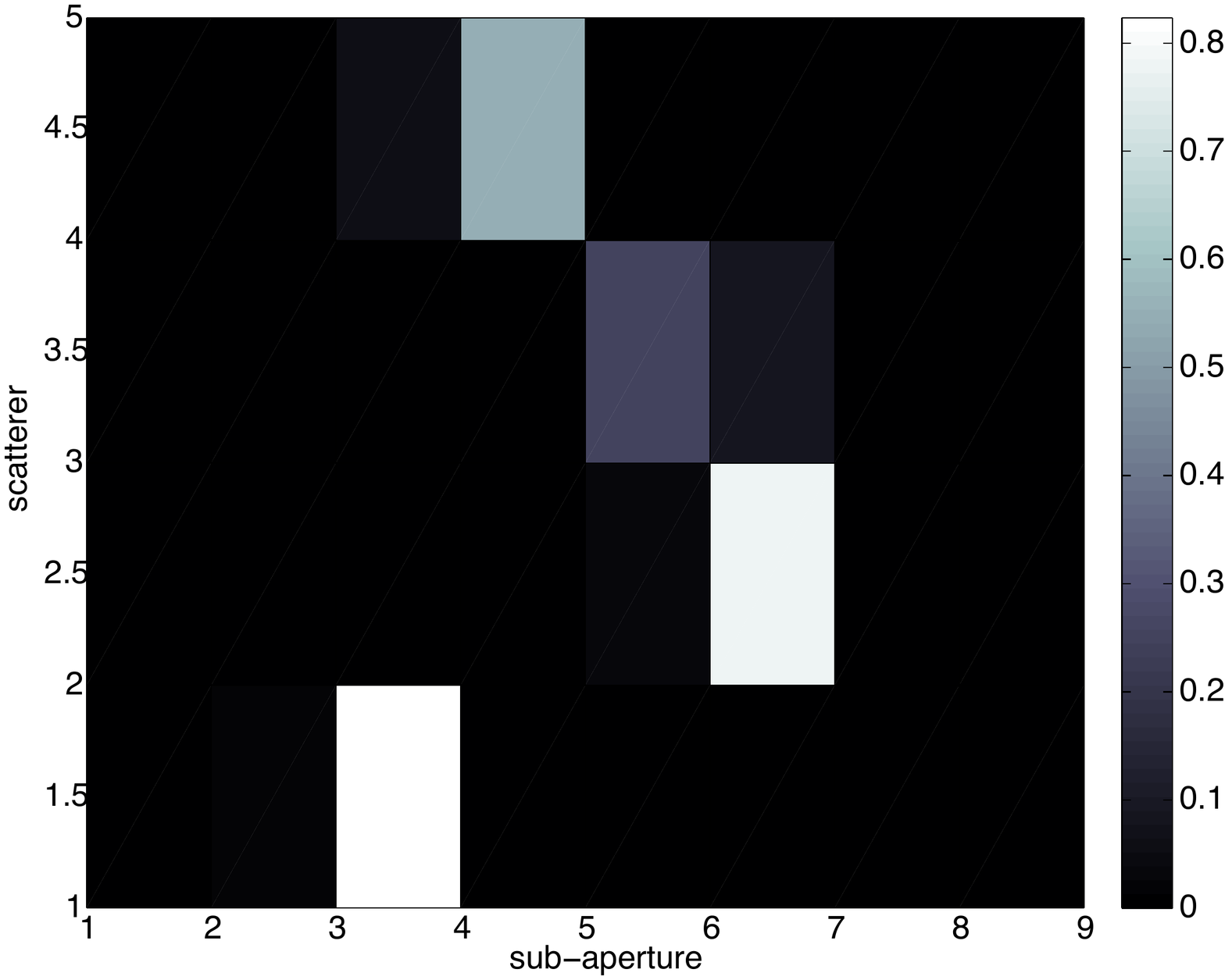}
\includegraphics[width=4.2cm]{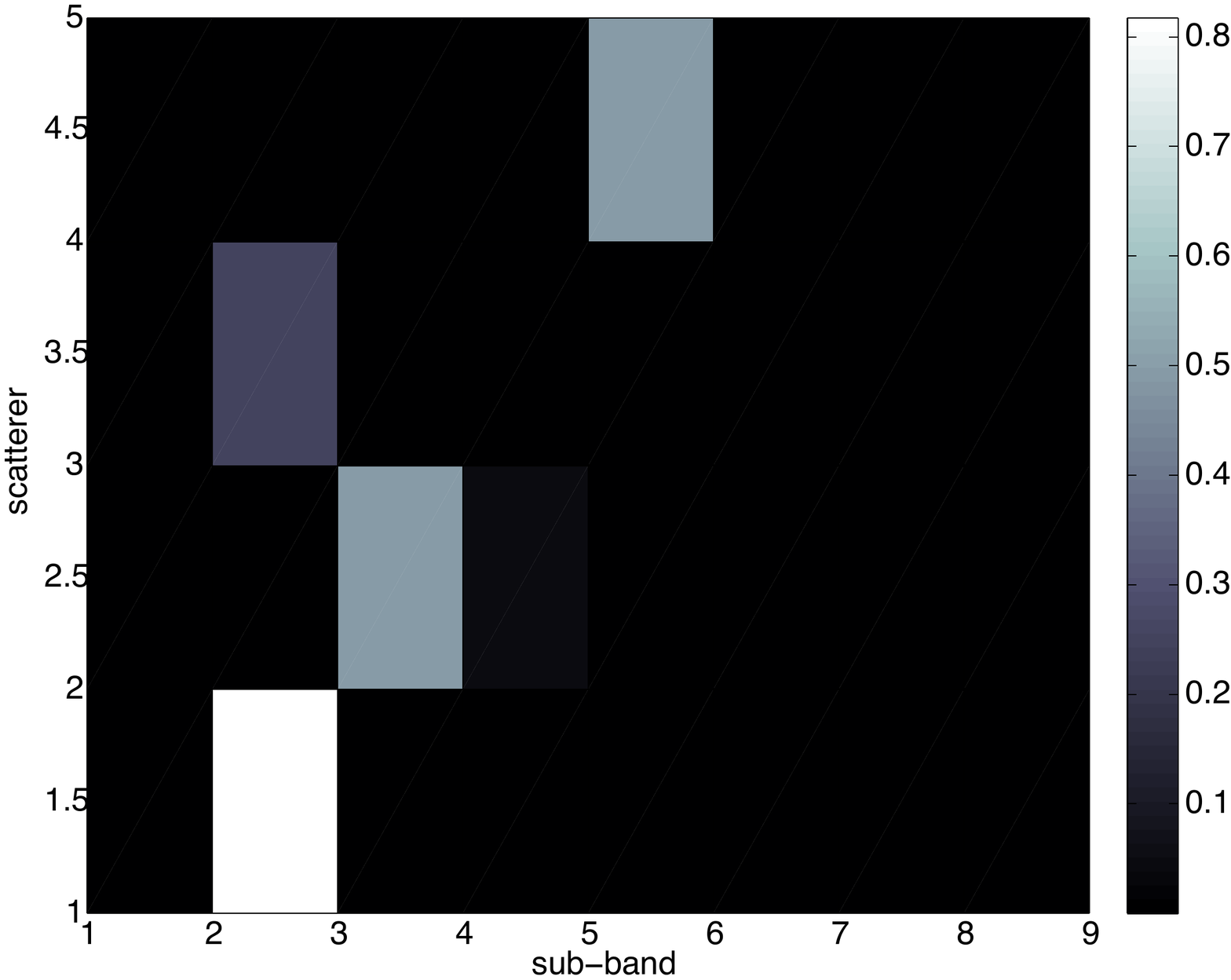}
\vspace{-0.6in}\\
\vspace{-0.6in}
\includegraphics[width=4.2cm]{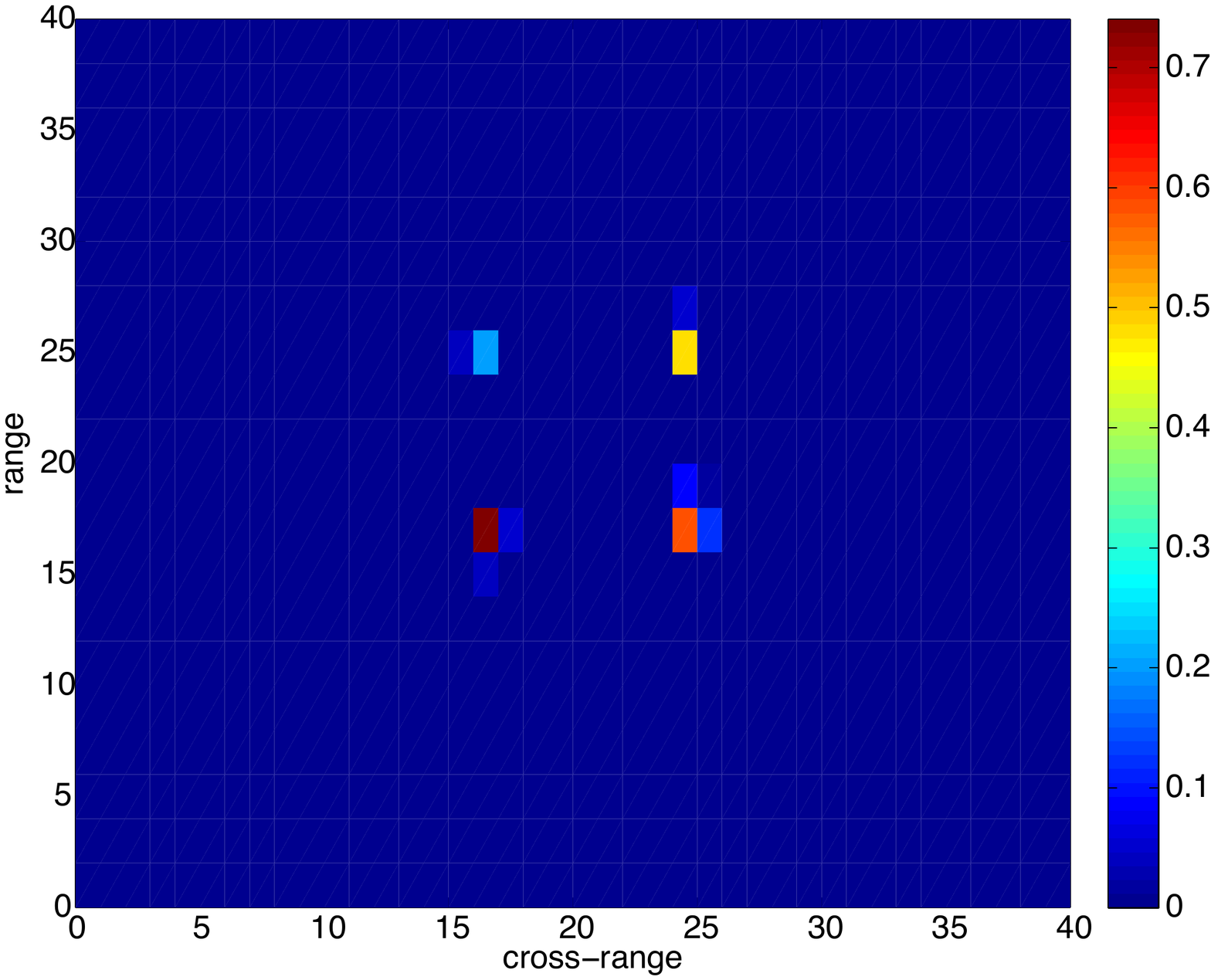}
\includegraphics[width=4.2cm]{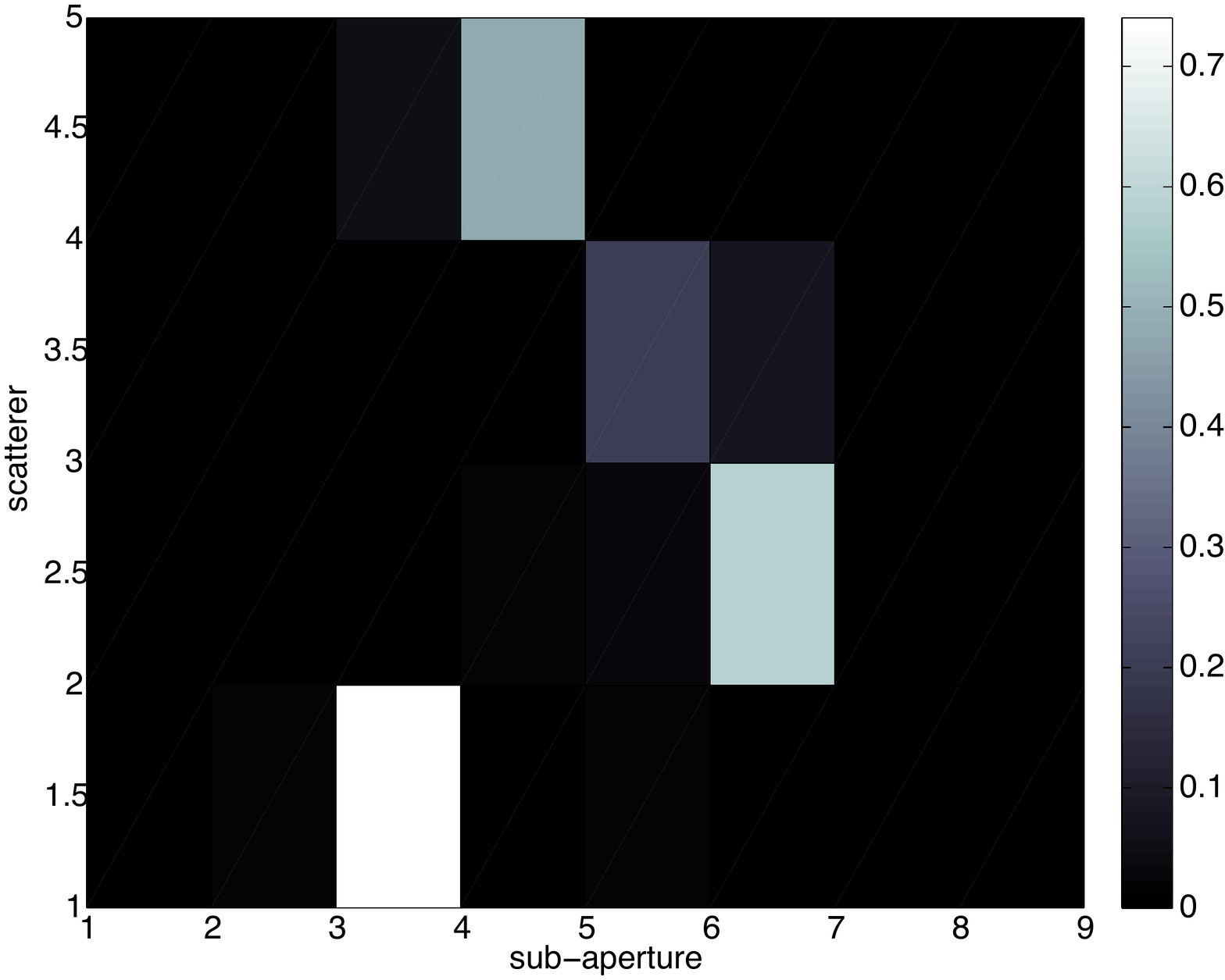}
\includegraphics[width=4.2cm]{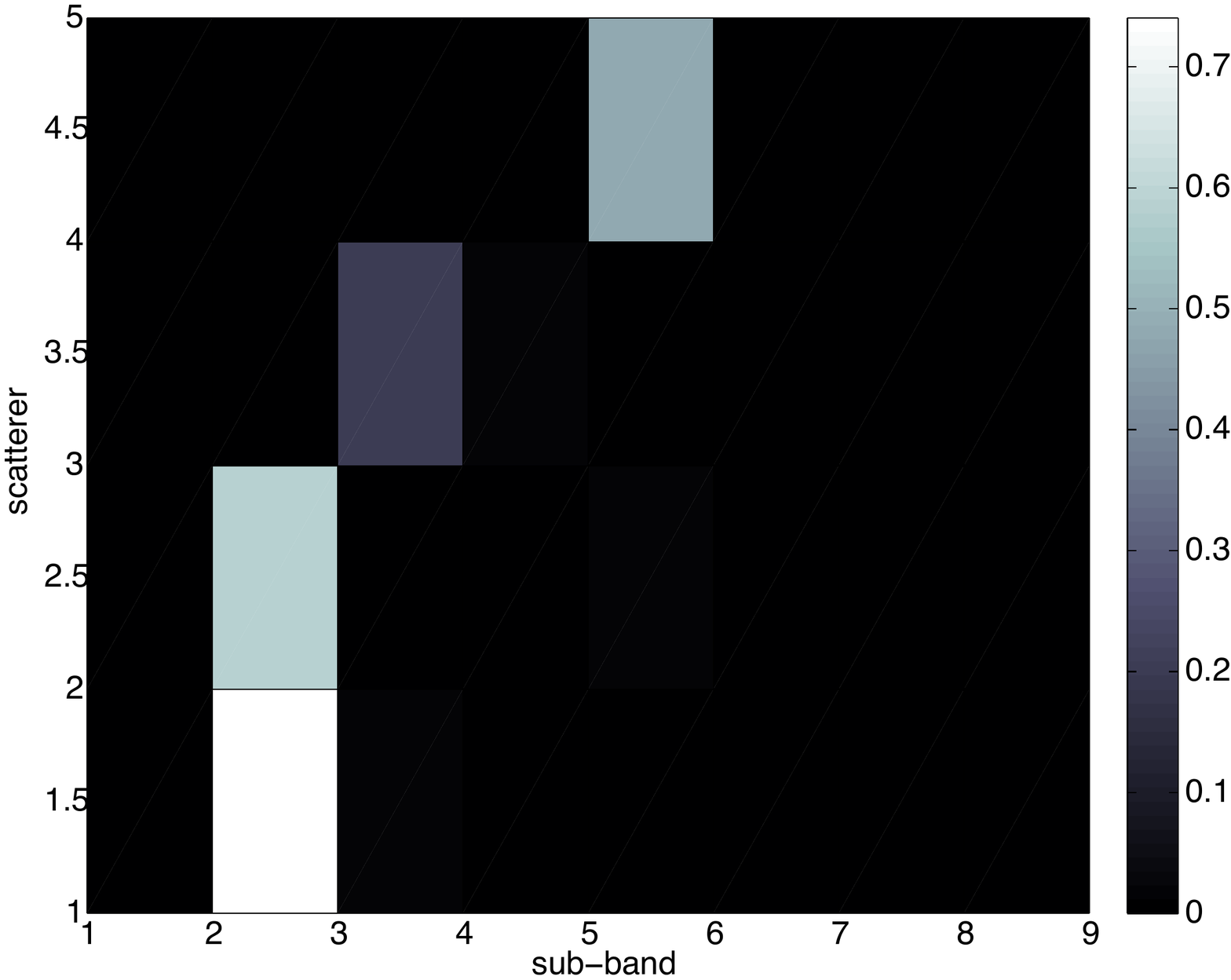}
\\ \vspace{-0.4in}\includegraphics[width=4.2cm]{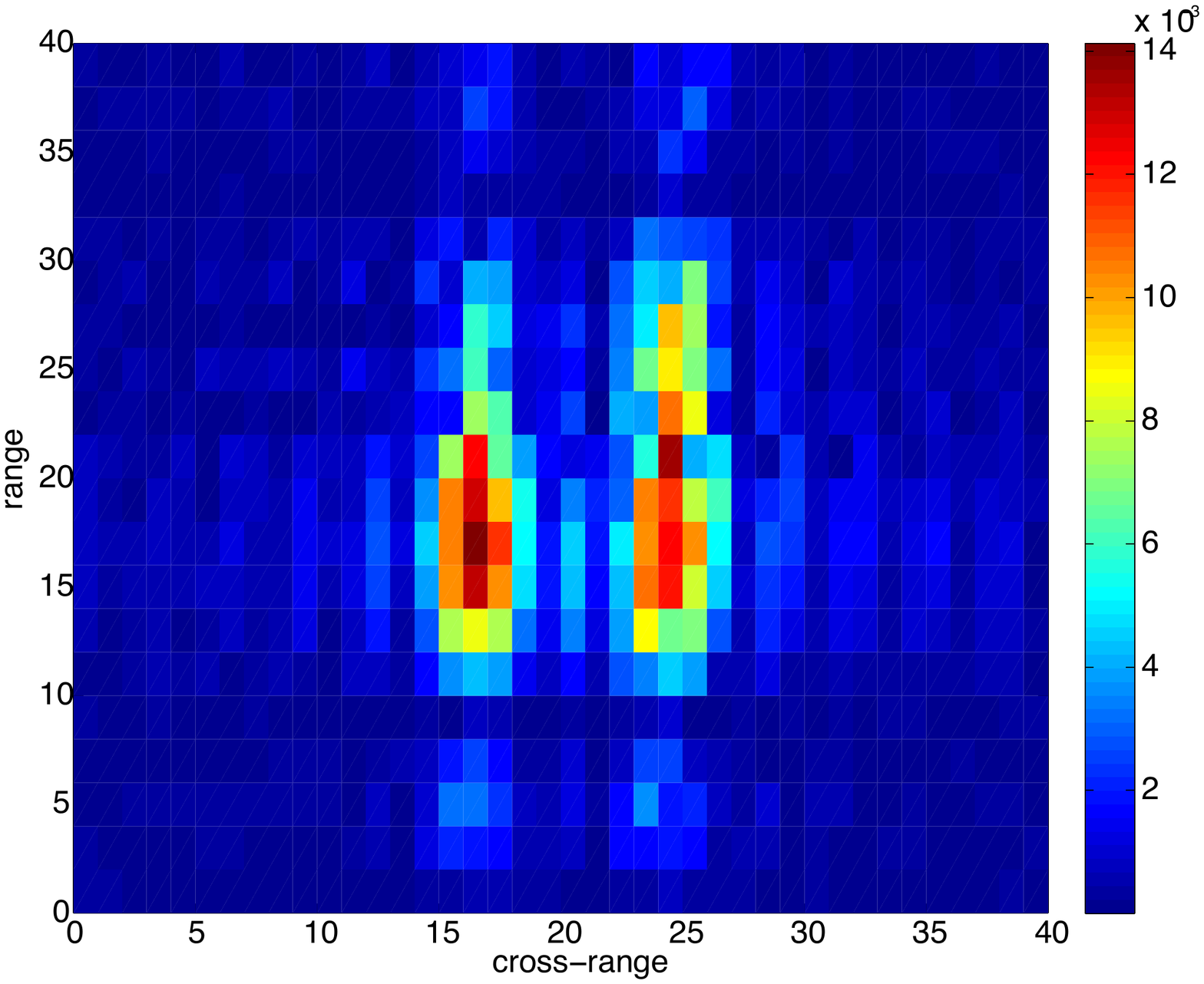}
\end{minipage}
\caption{Results with $\cN_\om = 8$ frequency intervals, $\cN_\alpha =
  8$ apertures and data contaminated with $20\%$ noise.  On the top we
  show the true reflectivity as a function of location (left) and
  direction (middle), and frequency (right).  In the middle row we show the reconstructed
  reflectivity with the MMV algorithm. The plot in the bottom row is
  the migration image.  The axes in the left images are cross-range
  and range in meters. The abscissa in the middle and right plots are
  the sub-aperture and sub-band index and the ordinate
  is the index of the scatterer (from $1$ to $4$). }
\label{fig:2df_3}
\end{figure}

\begin{figure}[t]
\includegraphics[width=4.2cm]{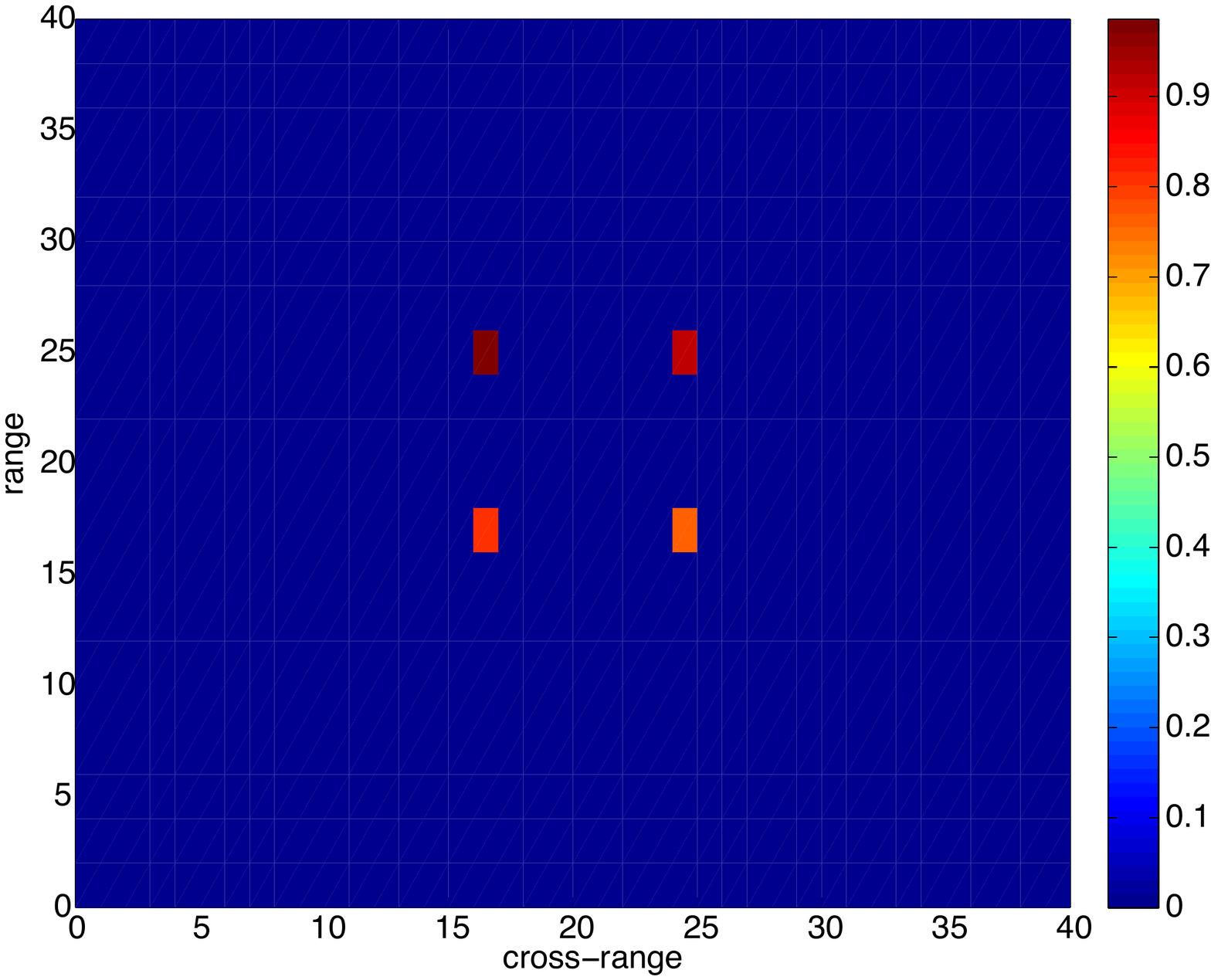}
\includegraphics[width=4.2cm]{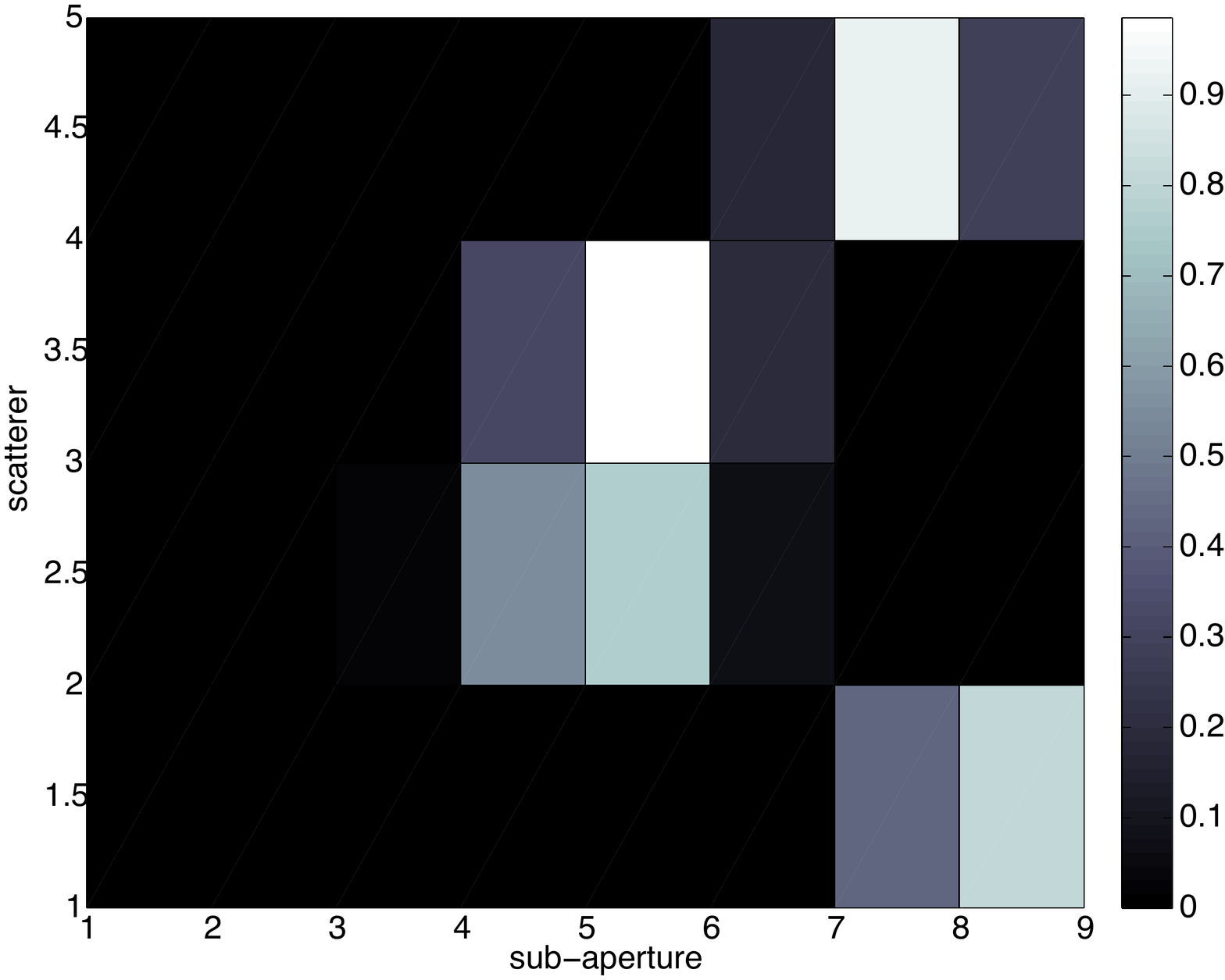}
\includegraphics[width=4.2cm]{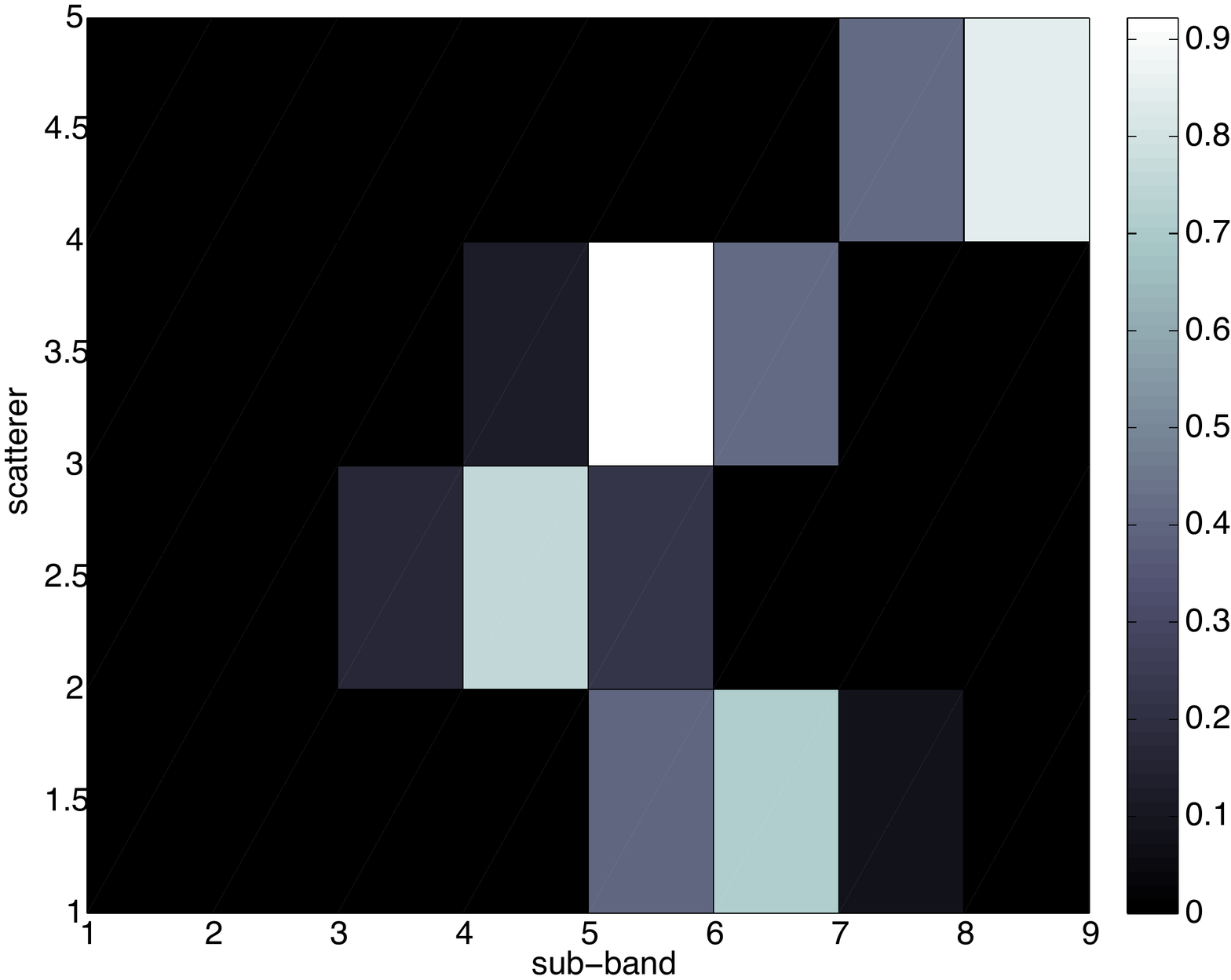}
\vspace{-0.6in}\\
\includegraphics[width=4.2cm]{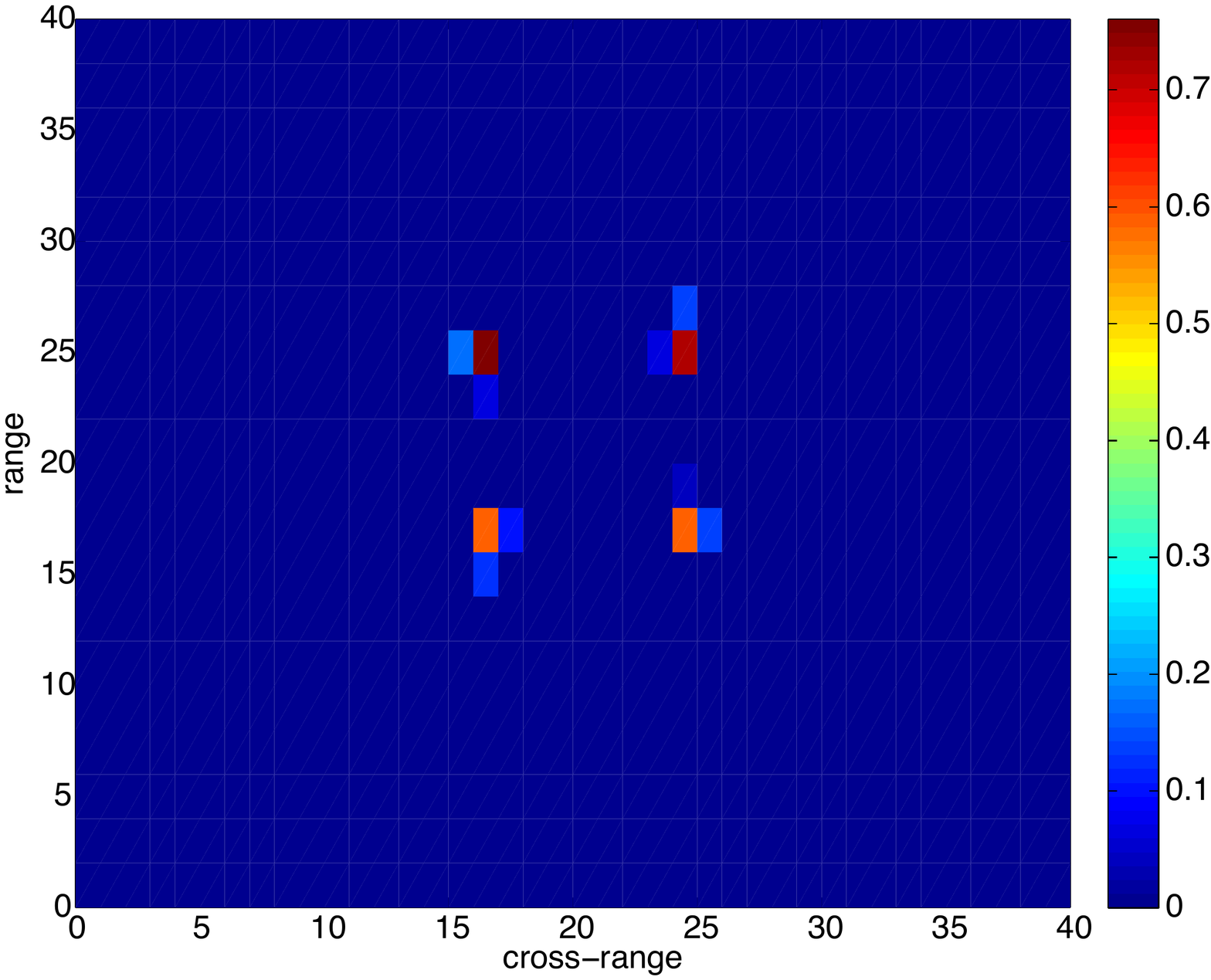} 
\includegraphics[width=4.2cm]{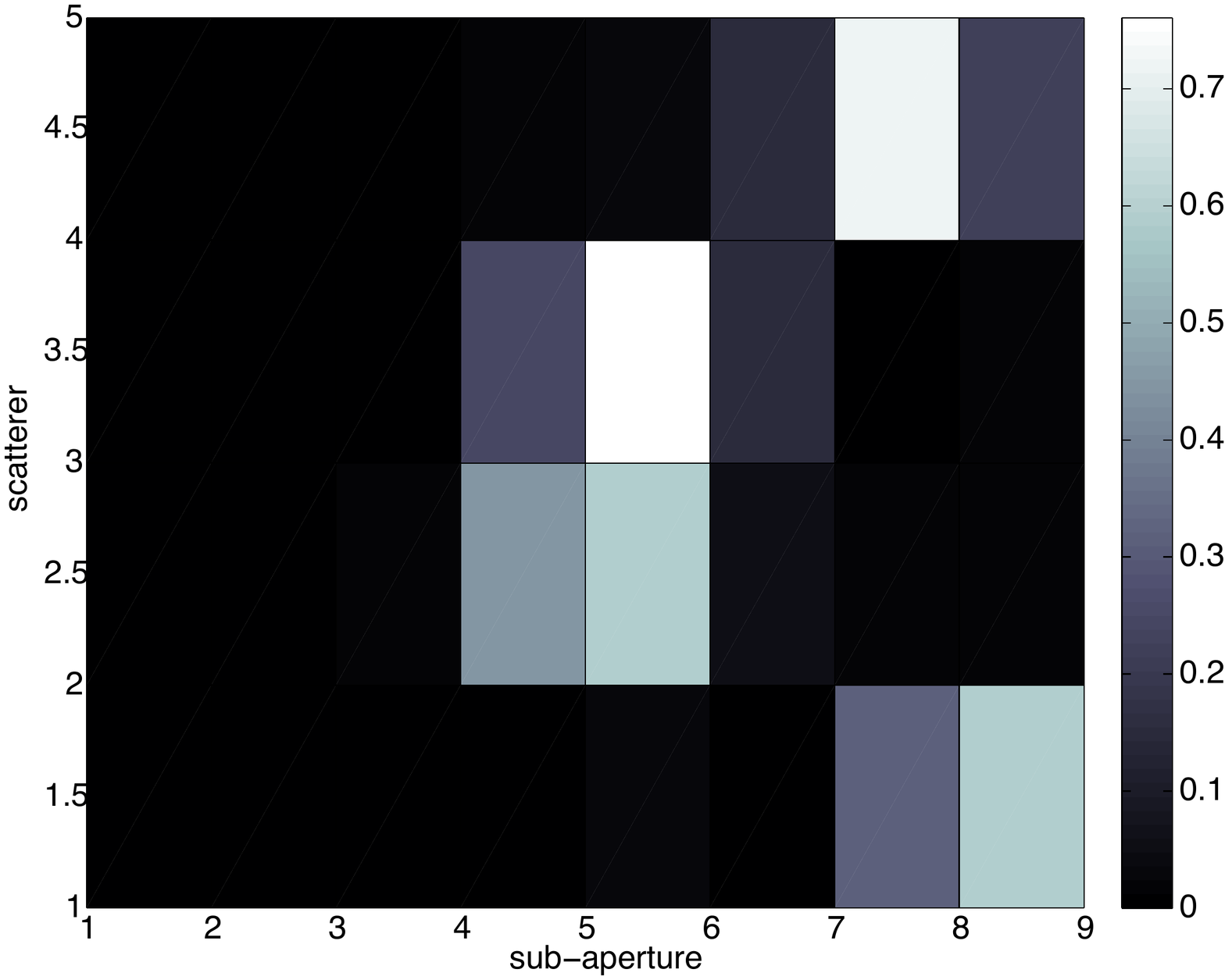}
\includegraphics[width=4.2cm]{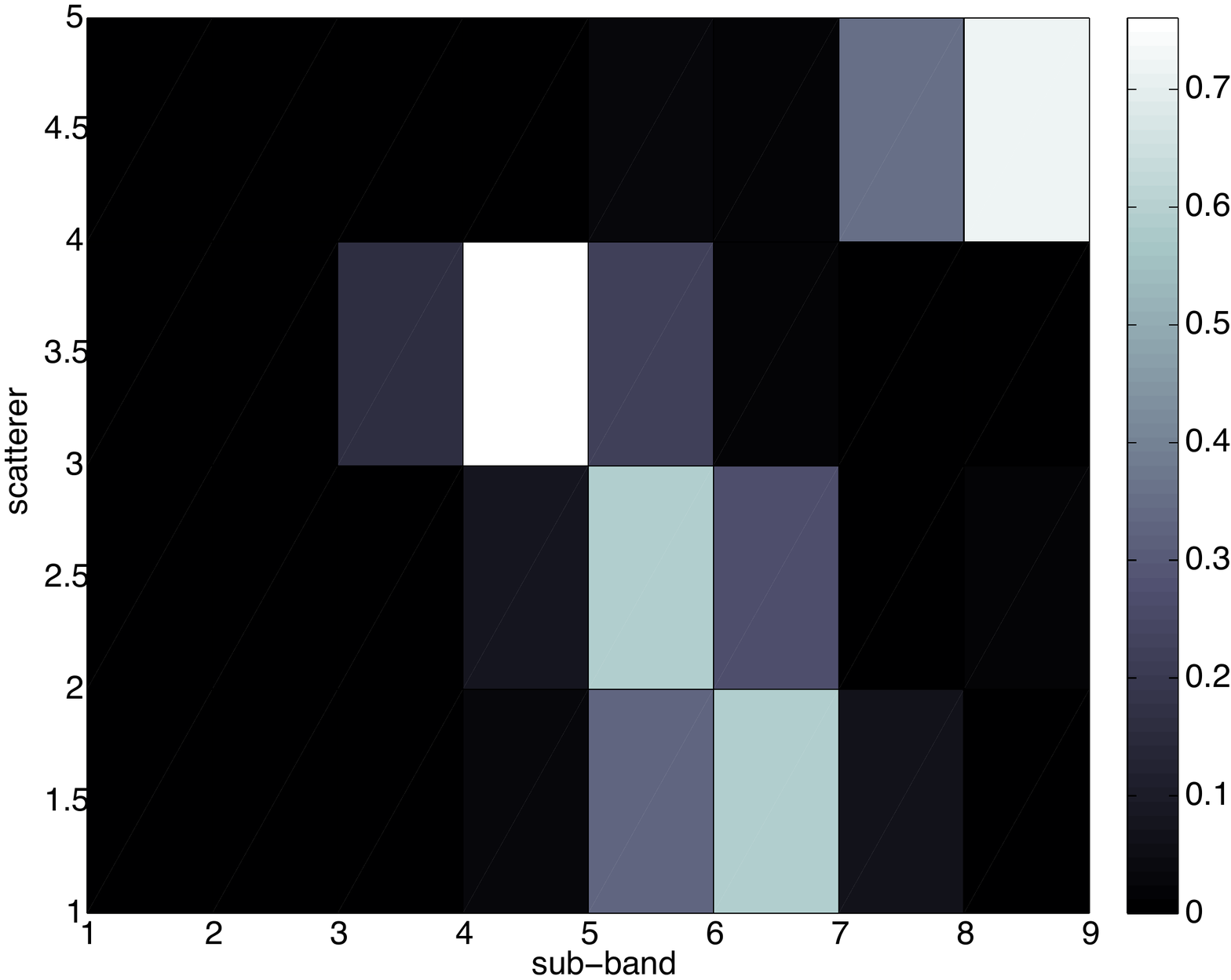}
\vspace{-0.6in}\\ \vspace{-0.4in}
\includegraphics[width=4.2cm]{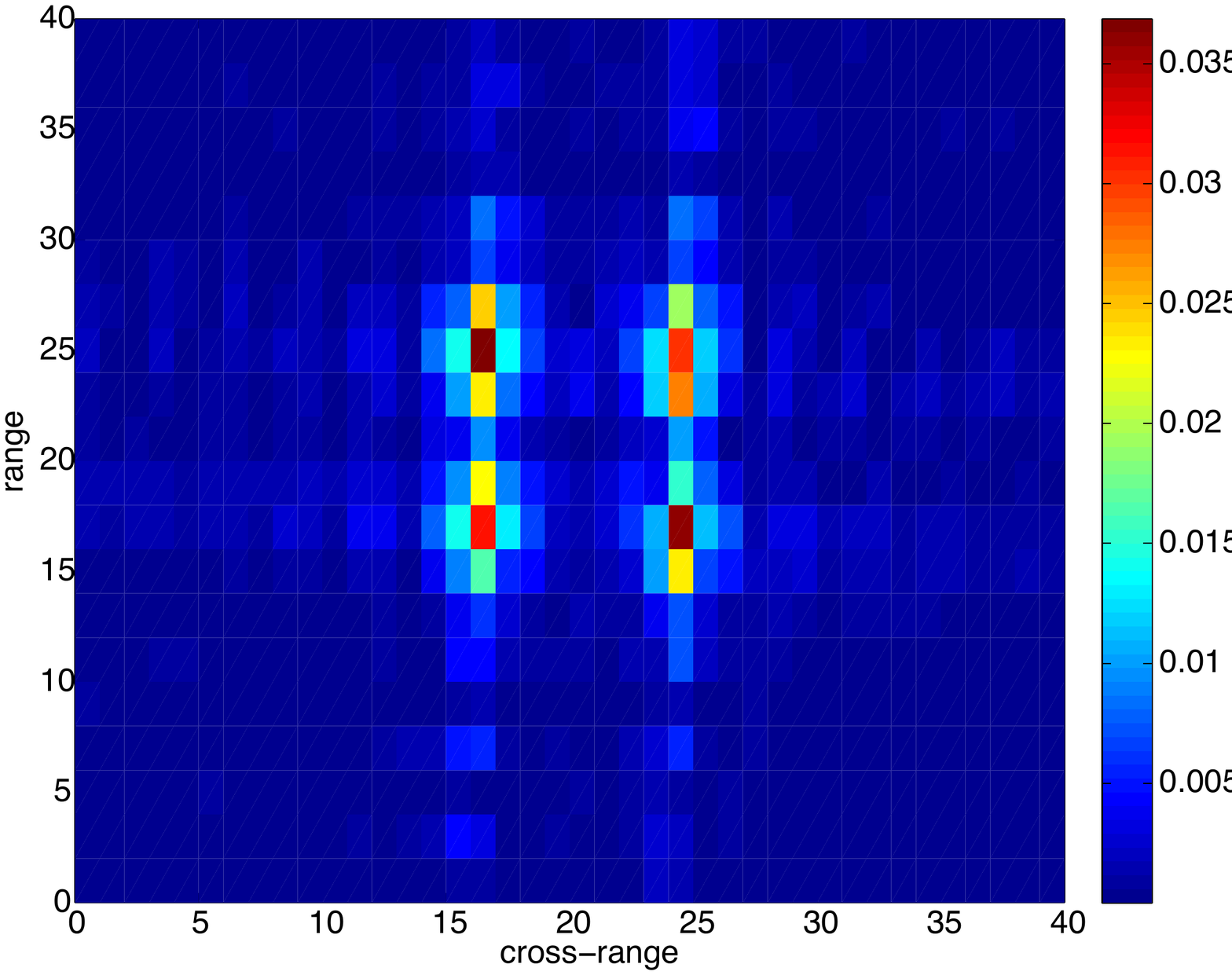} 
\caption{Results with $\cN_\om = 8$ frequency intervals, $\cN_\alpha =
  8$ apertures and data contaminated with $20\%$ noise.  On the top we
  show the true reflectivity as a function of location (left), direction (middle), and
  frequency (right).  In the middle row we show the reconstructed
  reflectivity with the MMV algorithm. The plot in the bottom row is
  the migration image.  The axes in the left images are cross-range
  and range in meters. The abscissa in the middle and right plots are
  the sub-aperture and sub-band index and the ordinate
  is the index of the scatterer.}
\label{fig:2df_4}
\end{figure}

The results in Figures \ref{fig:2df_3} and \ref{fig:2df_4} are for
$\cN_\om = 8$ consecutive, non-overlapping frequency bands and
$\cN_\alpha = 8$ consecutive, non-overlapping sub-apertures. The
difference between the figures is the strength of the scatterers and
their anisotropy. The results in Figure \ref{fig:2df_3} show that the
MMV algorithm reconstructs well the location of the scatterers and the
direction dependence of their reflectivity. The frequency dependence
of the weaker scatterers is not that accurate, likely because the
bandwidth is small and all frequencies are similar to the carrier. 
As in Figure \ref{fig:2d_real2}, the migration image is blurrier and
does not locate the weak scatterers. Figure \ref{fig:2df_4} shows that
the migration image improves when all scatterers are of approximately
the same strength and they have weaker anisotropy.

\begin{figure}[t]
\includegraphics[width=4.2cm]{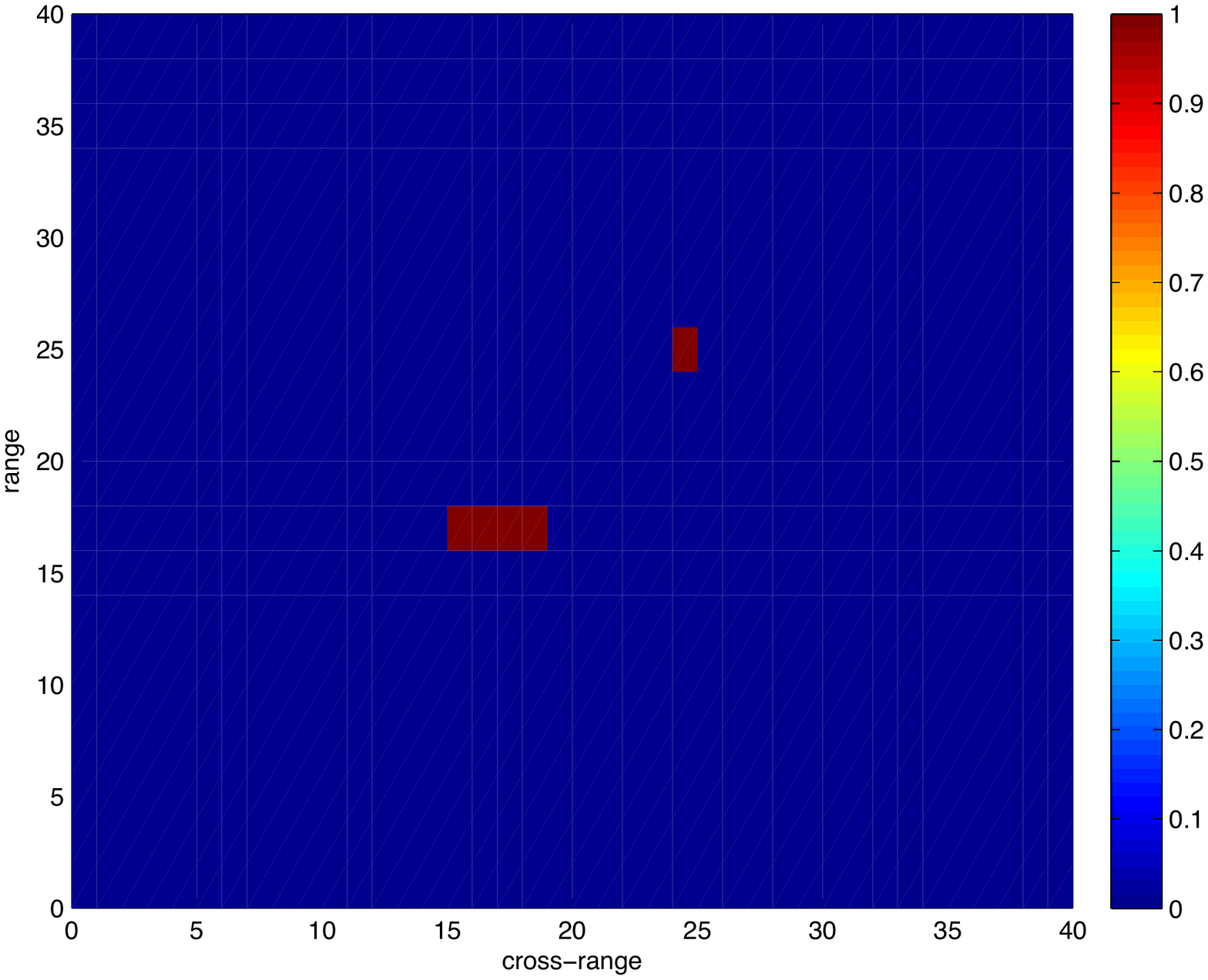}
\includegraphics[width=4.2cm]{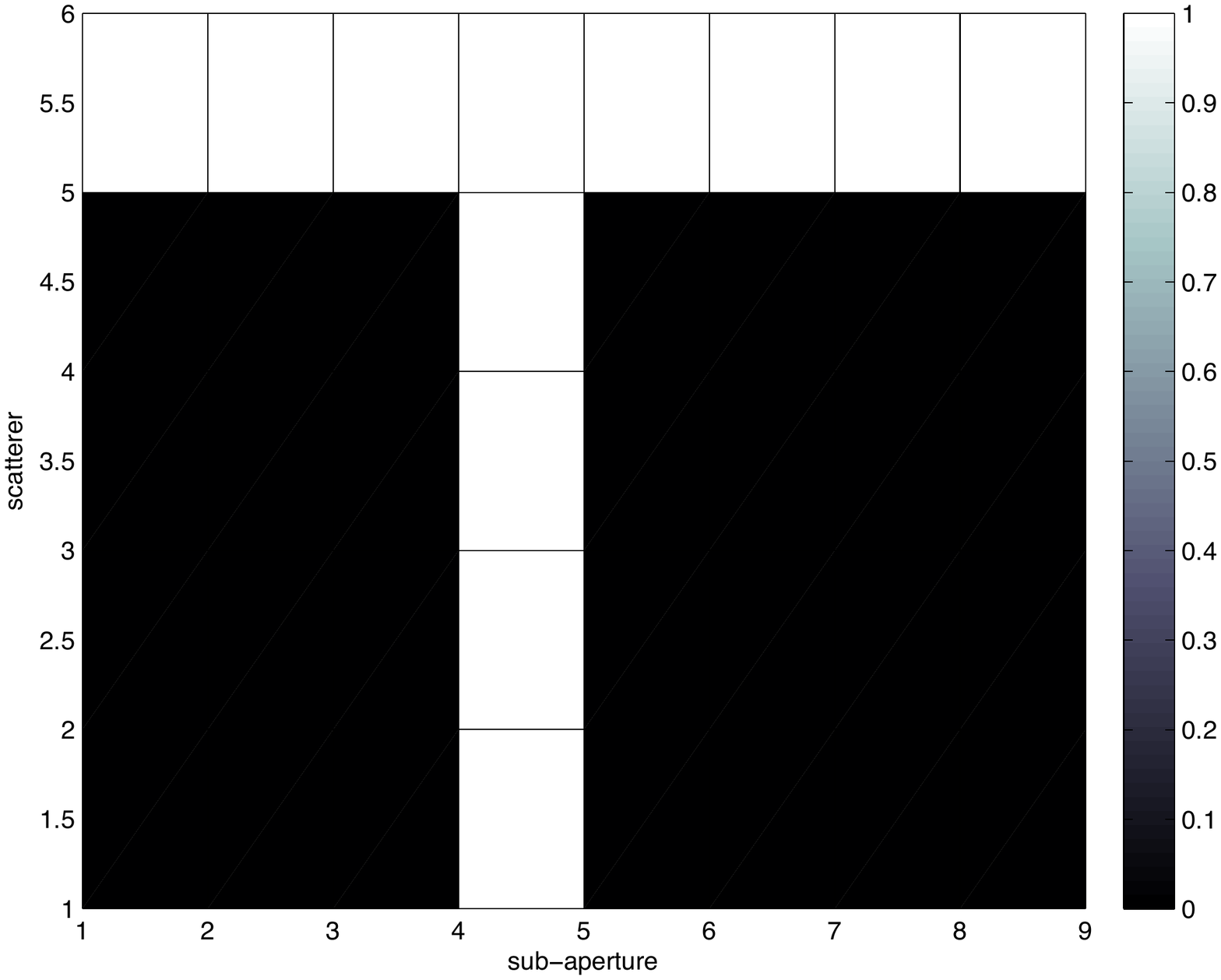}
\includegraphics[width=4.2cm]{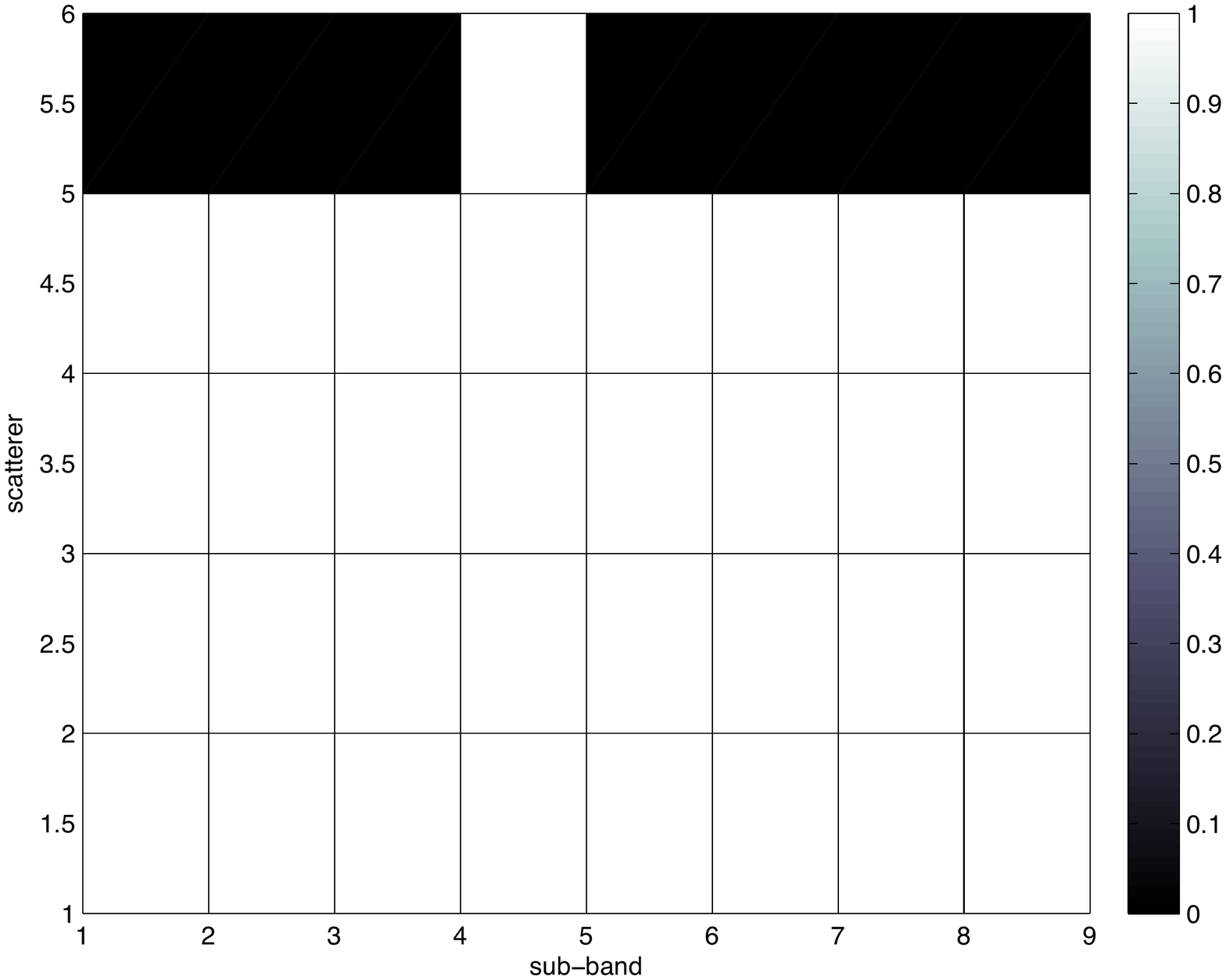}
\vspace{-0.6in}
\\
\includegraphics[width=4.2cm]{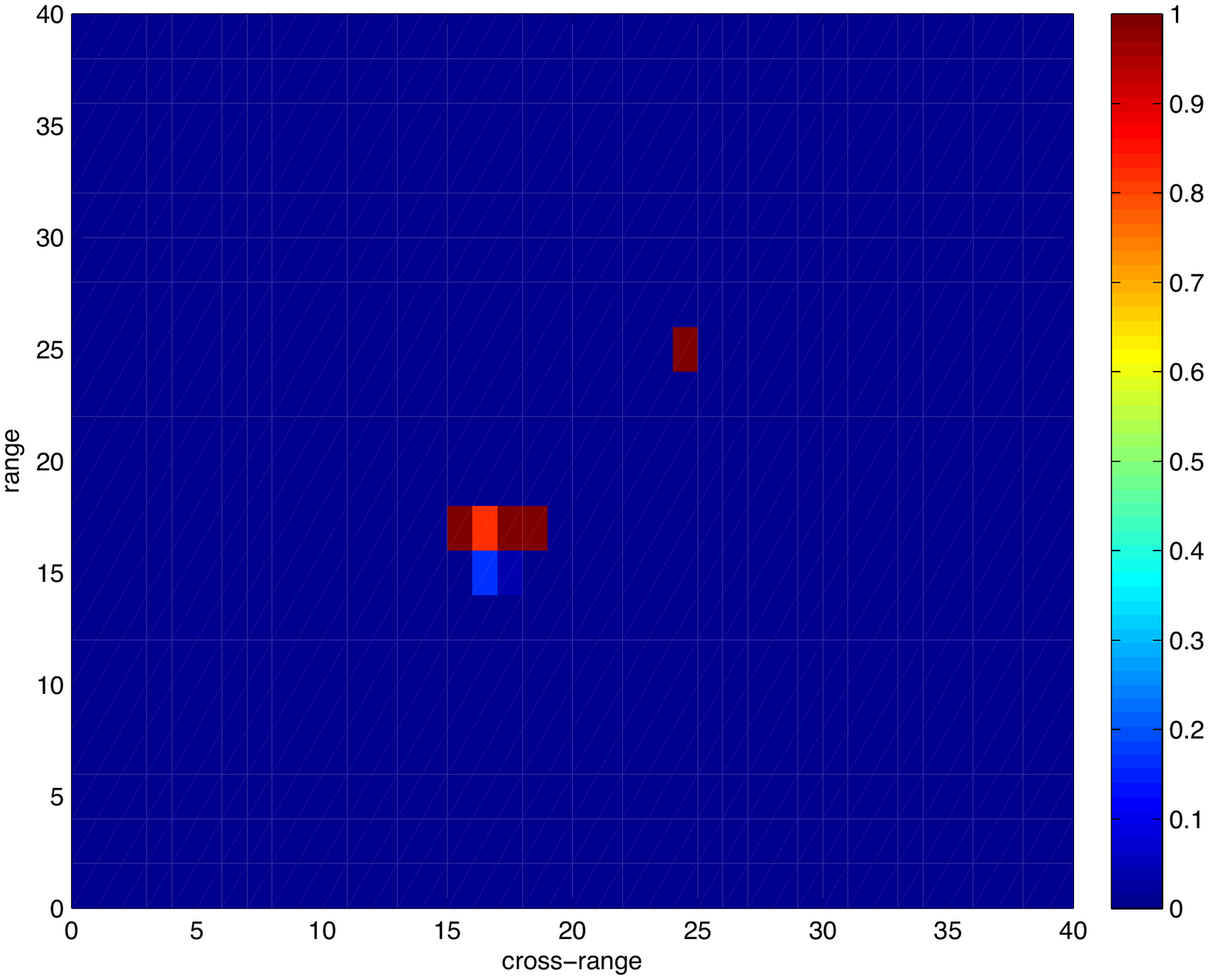}
\includegraphics[width=4.2cm]{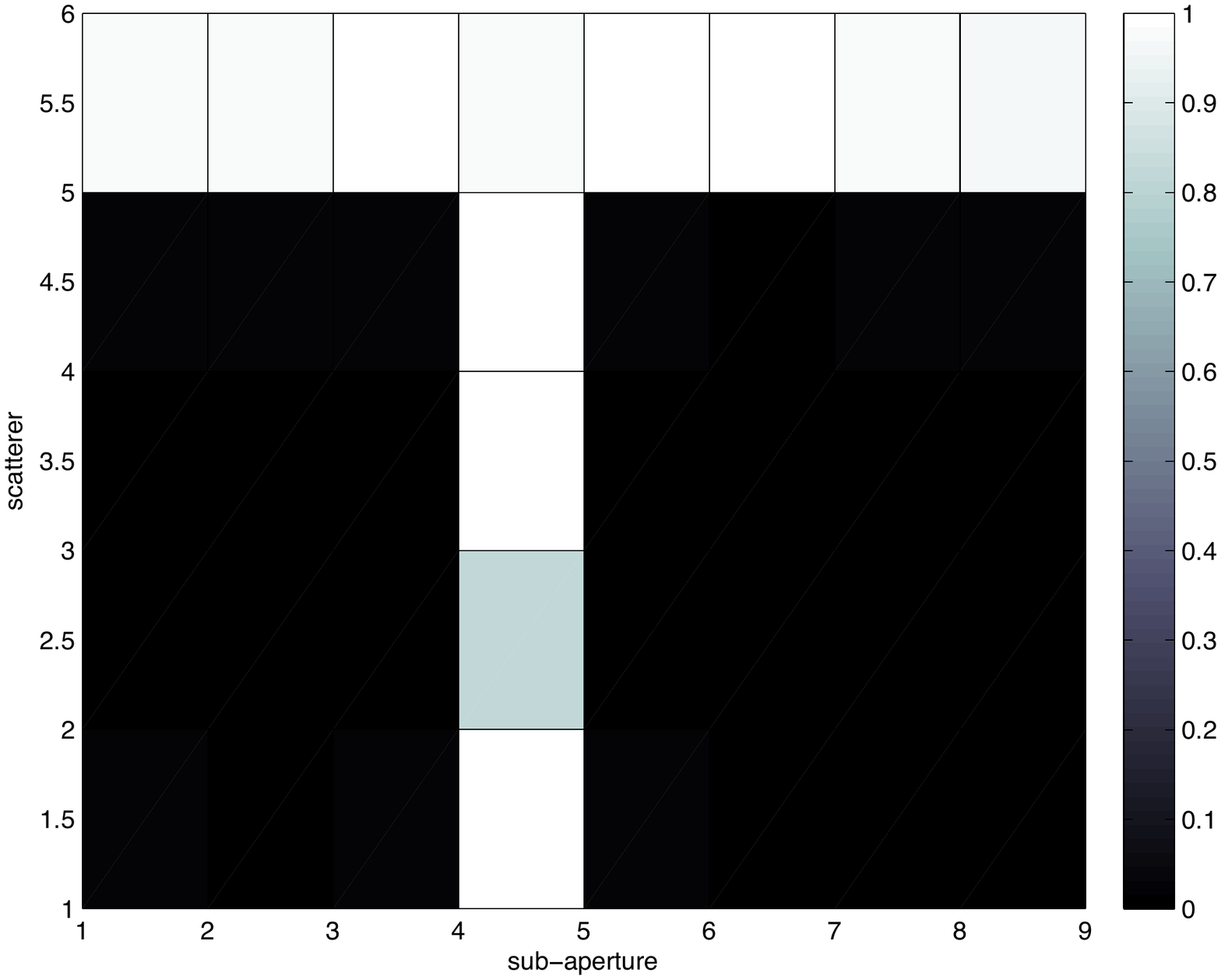}
\includegraphics[width=4.2cm]{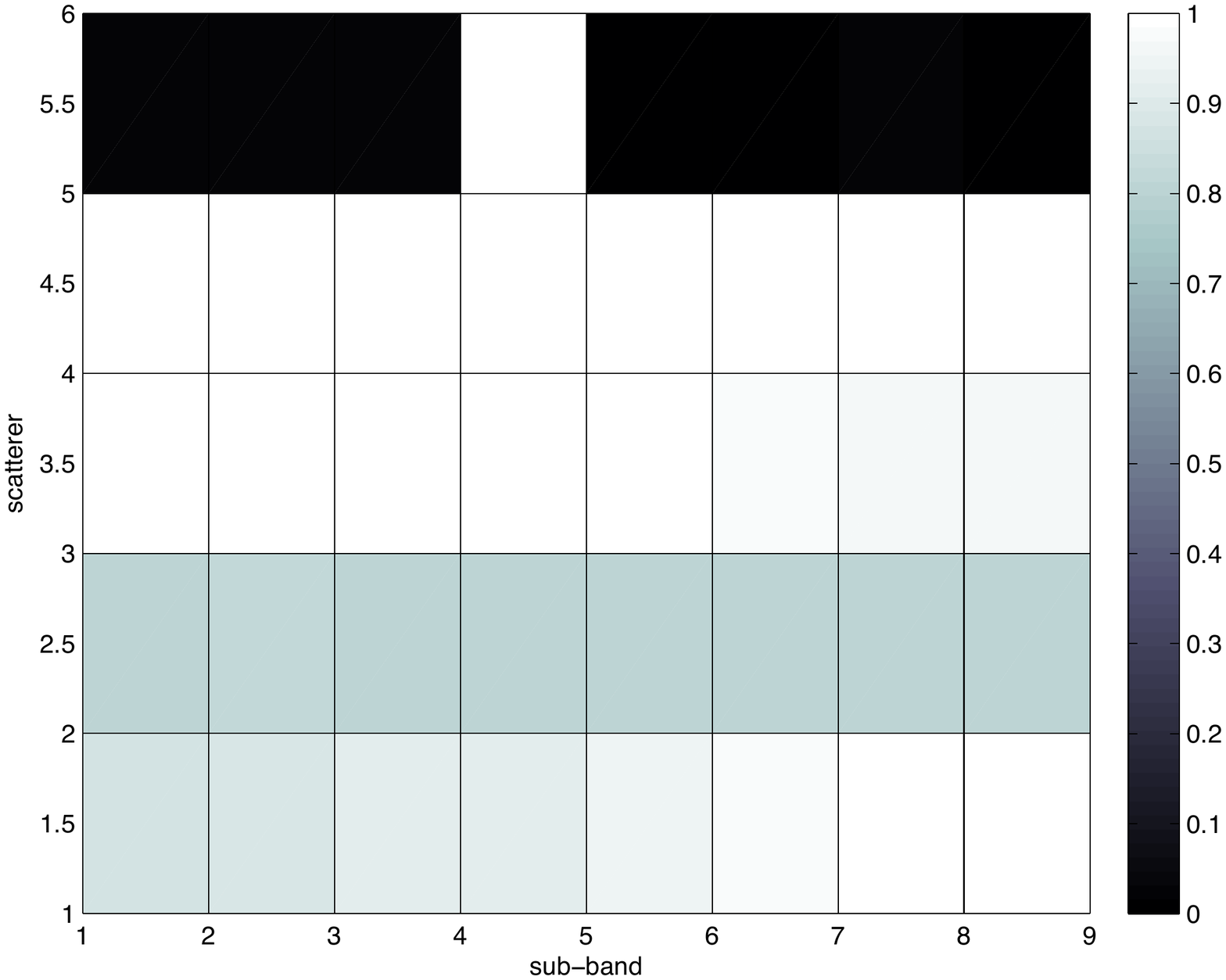}
\vspace{-0.6in}\\
\includegraphics[width=4.2cm]{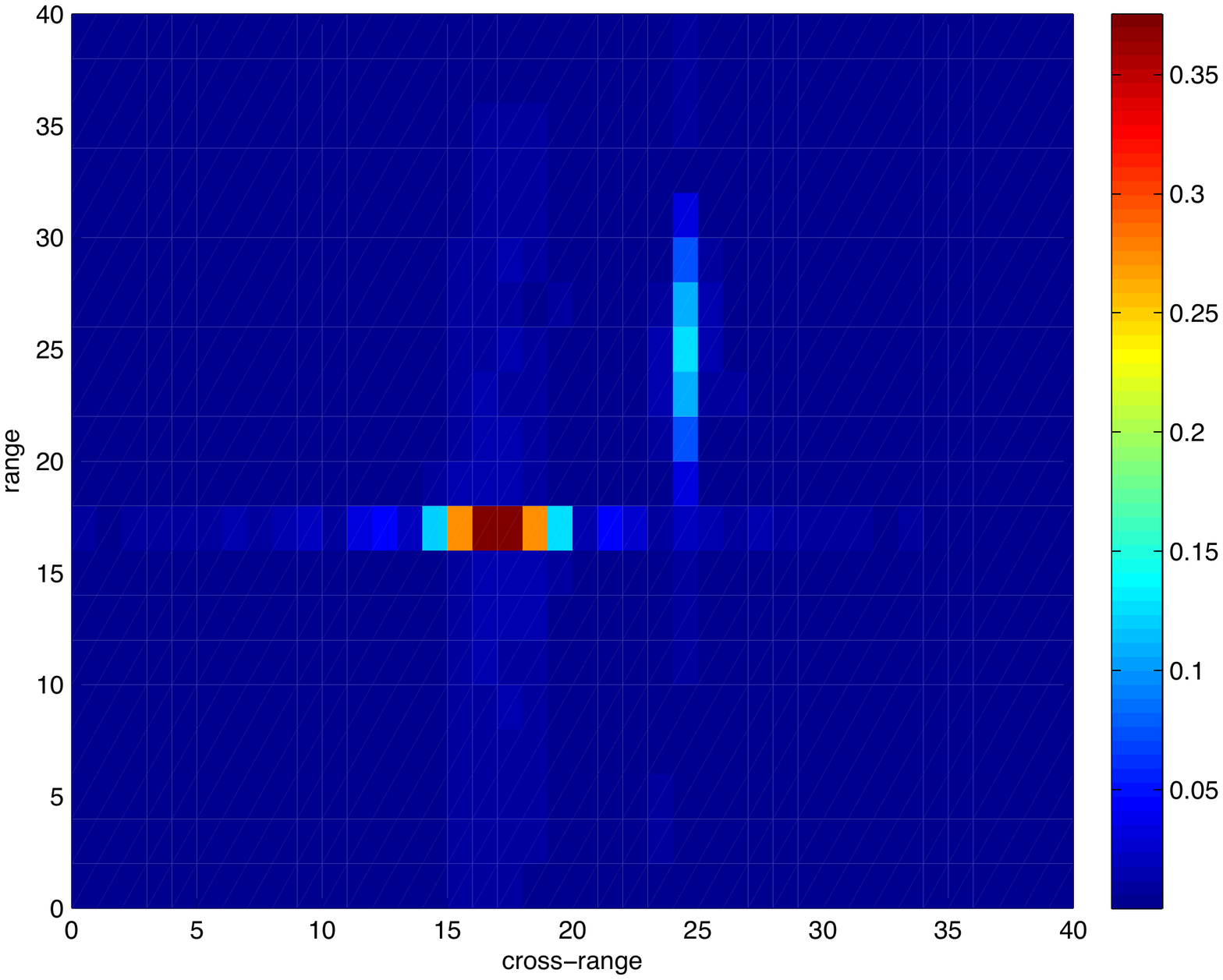}
\vspace{-0.4in}
\caption{Results with $\cN_\om = 8$ frequency intervals, $\cN_\alpha =
  8$ apertures and data contaminated with $20\%$ noise.  On the top we
  show the true reflectivity as a function of location (left), direction (middle), and
  frequency (right).  The reflectivity is supported in $5$ pixels, four of them are adjacent and represent
  an extended scatterer, which is frequency independent. The other pixel supports a small
  scatterer with frequency dependent reflectivity. 
  In the middle row we show the reconstructed
  reflectivity with the MMV algorithm. The plot in the bottom row is
  the migration image.  The axes in the left images are cross-range
  and range in meters. The abscissa in the middle and right plots are
  the sub-aperture and sub-band index and the ordinate
  is the index of  the pixel in the support of the scatterer. }
\label{fig:new2}
\end{figure}

{The last illustration considers a larger scatterer with direction dependent reflectivity, supported 
over four adjacent pixels,  and a small isotropic scatterer with frequency dependent 
reflectivity. The data is contaminated with $20\%$ additive noise.  We note that the migration method gives the 
correct location of the large scatterer, but  not the value of its reflectivity. Moreover, it gives a blurry image
of the small scatterer. The MMV algorithm determines well the support of both scatterers, as well as  
accurate estimates of the reflectivity as a function of direction and frequency.}

\section{Doppler effects}
\label{ap:Doppler}
All the results up to now use the start-stop approximation of the data
model, which neglects the motion of the platform over the fast time
recording window. Here we extend them to regimes where Doppler effects
are important. We begin in section \ref{sect:Dop1} with the derivation
of the generalized data model that includes Doppler effects, and an
assessment of the validity of the start-stop approximation. Then we
explain in section \ref{sect:Dop2} how to incorporate these effects in
our imaging algorithm.

\subsection{Data model with Doppler effects}
\label{sect:Dop1}
For simplicity we first derive the data model for an isotropic
reflectivity $\rho = \rho(\vy)$.  Then we extend it in the obvious way to direction and
frequency dependent reflectivities in a sub-aperture indexed by $\alpha$
and sub-band indexed by $\beta$, with reflectivity
$\rho^{(\alpha,\beta)}(\vy)$.

The scattered wave $u(s,t)$ recorded at the transmit-receive platform is given by
\begin{align}
  u(s,t) &= - \int_{\Omega} d \vy \frac{\rho(\vy)}{c^2} \int_0^t dt_1 \int_0^{t_1}
  d t_2 \, f''(t_2) G(t_1-t_2,\vr(s+t_2),\vy)G(t-t_1,\vy,\vr(s+t)), \nonumber
  \\ &= - \frac{1}{c^2}\int_{\Omega} d \vy \rho(\vy) \, \frac{f''\big(t_2(t)\big)}{(4 \pi)^2 |\vr\big(s+t_2(t)\big)-\vy|
    |\vr(s+t)-\vy|} \label{eq:DM1}
\end{align}
where $t_2(t)$ is the solution of the equation 
\begin{equation}
t_2 + \frac{|\vr(s+t_2) - \vy|}{c} = t - \frac{|\vr(s+t)-\vy|}{c},
\label{eq:DM1p}
\end{equation}
and 
we used the expression of the Green's function of the wave equation
\[
  G(t,\vr,\vy) = \frac{\delta \big[t-|\vr-\vy|/c\big]}{4 \pi |\vr-\vy|},
\]
and the single scattering approximation. The expression (\ref{eq:DM1}) is
simply the spherical wave emitted from $\vr(s+t_2)$, over the duration $t_2$ of the pulse, scattered isotropically at $\vy$,
and then recorded at $\vr(t+s)$. Up to the single scattering approximation,
this is an exact formula. Expanding with respect to $t$ 
the arguments in (\ref{eq:DM1}) and (\ref{eq:DM1p}) we obtain
\begin{align}
  u(s,t) &= - \frac{1}{c^2} \int_\Omega d \vy \rho(\vy)
  \frac{1}{(4 \pi |\vr(s)-\vy|)^2\big(1 + O(Vt/L)\big)} \nonumber \times 
  \\ 
  & \hspace{-0.2in}f''\Big[\Big( t \Big(1 - \gamma(s,\vy) +
      O\Big(\frac{V}{c}\frac{Vt}{R}\Big)\Big) - 2 \tau(s,\vy)\Big)/\Big(1+\gamma(s,\vy)+ O\Big(\frac{V}{c}\frac{Vt}{R}\Big)\Big)\Big],
      \label{eq:DM2p}
\end{align}
where we introduced the Doppler factor $\gamma$ defined by
\begin{align}
  \gamma(s,\vy) = \frac{\vr'(s)}{c} \cdot \vm(s,\vy), \qquad
  \vm(s,\vy) = \frac{\vr(s)-\vy}{|\vr(s)-\vr(y)|}.
  \label{eq:DM3}
\end{align}
We assume that the platform is moving at constant speed $V$
along a trajectory with unit tangent denoted by $\vt(s)$, and
with radius of curvature $R$ assumed comparable to the range $L$. Thus,
\[
\gamma(s,\vy) = O \Big(\frac{V}{c}\Big) \ll 1,
\]
because the platform speed is typically much smaller than $c$, the wave speed, 
and we can neglect the residual in
(\ref{eq:DM2p}) which is even smaller than $\gamma$, because over the
duration of the fast time window the platform travels a small distance
compared with the radius of curvature $V t \ll R \sim L$. We have thus the data model
\begin{equation}
u(s,t) 
\approx - \frac{1}{c^2} \int_\Omega d \vy
  \rho(\vy) \frac{f''\Big[ t \big(1 - 2\gamma(s,\vy)\big) - 2 \tau(s,\vy)\big(1-\gamma(s,\vy)\big)
    \Big]}{(4 \pi |\vr(s)-\vy|)^2},
  \label{eq:DM2}
 \end{equation}
 which includes first order Doppler effects.

The start-stop approximation is valid when the Doppler factor in the
argument of $f''$ in (\ref{eq:DM2}) is negligible. Although $\gamma$ is small, 
$f''$ oscillates at the carrier frequency $\om_o$ which is large
and, depending on the scale of the fast time $t$, the Doppler factor
may play a role. Recall that $t$ is limited by the slow time spacing  
$h_s$. In practice the duration of the fast time window may be much smaller than $h_s$, although 
it must be large enough so that the platform can receive the echoes delayed by
the travel time, $2 \tau(s,\vy)$. Explicitly,
\[
t = O(L/c) + O(1/B),
\]
where $L/c$ is the scale of the travel time and $1/B$ is the scale of the
duration of the signal.

We conclude that  the start stop approximation holds when
\[
\om_o t \gamma(s,\vy) = O \Big(\frac{\om_o L}{c} \frac{V}{c} \Big) +
O \Big(\frac{\om_o}{B} \frac{V}{c} \Big) \ll 1.
\]
In the GOTCHA regime, considered in the numerical simulations in
section \ref{sect:NUM}, we have
\[
\frac{\om_o L}{c} \frac{V}{c} = 0.469, \qquad \frac{\om_o}{B}
\frac{V}{c} = 2.3 \cdot 10^{-5},
\]
so $\om_o t \gamma(s,\vy)$ is slightly less than one. We may include
it in the data model, but it amounts to a constant additive phase that
has no effect in imaging. 
To see this, let us take the Fourier transform with respect to $t$ in
(\ref{eq:DM2})
\begin{equation}
  \hat u(s,\om) \approx k^2
  \int_\Omega d \vy \rho(\vy) \hat f \Big[\om \big(1 + 2 \gamma(s,\vy)
    \big) \Big] \frac{\exp \big[ 2 i \om \big(1 + \gamma(s,\vy) \big)
  \tau(s,\vy)\big]}{(4 \pi |\vr(s)-\vy|)^2},
    \label{eq:DMn2}
\end{equation}
and expand the arguments over the slow time $s$ and imaging point
$\vy$. We use the approximation
\begin{equation}
\label{eq:DD3}
\vr'(s) \approx V \Big[\vt(s^\star) - \vec{\bf n}(s^\star) \frac{V \Delta
    s}{R}\Big],
\end{equation}
where $\Delta s$ is the slow time offset from the center $s^\star$ of
the aperture, and $\vt(s^\star)$ is the unit tangent to the trajectory of the
platform at the center point. The second term in (\ref{eq:DD3})
accounts for the curved platform trajectory, with unit vector $\vec{\bf
  n}(s^\star)$ orthogonal to $\vt$, in the plane defined by $\vt$ and the
center of curvature, and $R$ the radius of curvature.  We also have
\[
 | \vm(s,\vy) -
  \vm(s^\star,\vy_o)| =  O\Big(\frac{V |\Delta s|}{L} \Big) +
  O\Big(\frac{|Y^\perp|}{L} \Big),
\]
and
\[
  \om_o \tau(s,\vy) = \om_o \tau(s^\star,\vy_o) + O (k_o V\Delta s) 
  + O(k_o \Delta y).
\]
{Substituting in (\ref{eq:DMn2}) and using the parameters of the GOTCHA
regime, we see that,
\[
\omega \gamma(s,\vy) \tau(s,\vy) \approx \om_o \gamma(s^\star,\vy_o^) \tau(s^\star,\vy_o),
\]
so indeed, the Doppler effect amounts to 
a constant phase term.}

\subsection{Imaging algorithm with Doppler effects}
\label{sect:Dop2}
The model of the down-ramped data with the Doppler correction follows
from (\ref{eq:DMn2}),
\begin{align}
d(s,\om) &= \overline{\hat f \big[ \om \big( 1+2\gamma(s,\vy_o)\big)
    \big]} \hat u(s,\om) \exp \big[-2 i \om \big(1+\gamma(s,\vy_o)\big)
  \tau (s,\vy_o) \big] \nonumber \\ &\approx k^2 \overline{\hat f
  \big[ \om \big( 1+2\gamma(s,\vy_o) \big) \big]} \int_\Omega d \vy \,
\hat f \big[ \om \big( 1+2\gamma(s,\vy) \big) \big]\rho(\vy)
\times \nonumber
\\ &~ ~~~ \frac{\exp \big[ 2 i \om \big(1+\gamma(s,\vy)\big)
    \tau(s,\vy) - 2 i \om\big(1+\gamma(s,\vy_o)\big)
    \tau(s,\vy_o)\big]}{\big(4 \pi |\vr(s)-\vy|\big)^2}.
\label{eq:D1}
\end{align}
We are interested in direction and frequency dependent reflectivities,
so to use formula (\ref{eq:D1}), we consider next the $\alpha-$th
sub-aperture and the $\beta-$th sub-band, where we can replace $\rho$
by $\rho^{(\alpha,\beta)}(\vy)$. The data is denoted by
$d^{(\alpha,\beta)}(\Delta s,\Delta \om),$ where $\Delta s = s -
s_\alpha^\star$ and $\Delta \om = \om - \om_\beta^\star.$ The goal of the
section is to include Doppler effects in the statements of Lemma \ref{lem.1} and
Proposition \ref{lem.2}, which are the basis of our imaging
algorithm.

We begin with the observation that
\begin{equation}
\om \gamma(s,\vy) \tau(s,\vy)= \frac{\om}{c} \frac{\vr'(s)}{c} \cdot
(\vr(s)-\vy) = \om \gamma(s,\vy_o) \tau(s,\vy_o) - \frac{\om}{c}
\frac{\vr'(s)}{c} \cdot \Delta \vy,
\end{equation}
where $\Delta \vy = \vy - \vy_o$, and $\vr'(s)$ is given by
(\ref{eq:DD3}), and  assume henceforth that
\begin{equation}
\frac{V}{c} \frac{\cYp}{L_\alpha} \ll \frac{b}{\om_o} \ll 1.
\label{eq:D4}
\end{equation}
This is consistent with our previous assumptions because $\cYp \ll
L_\alpha$ and $V \ll c$, and allows us to approximate the Doppler factor in
the argument of the Fourier transform of the signal in (\ref{eq:D1})
by its value at the reference point. Then, using equation
(\ref{eq:L1.2}) and noting also that
\[
|\vr(s)-\vy| = L_\alpha \Big[ 1 + O \Big(\frac{a}{L_\alpha}\Big) +
  O\Big(\frac{\cYp}{L_\alpha}\Big)\Big], \qquad k = k_o \Big[ 1 +
  O\Big( \frac{b}{\om_o} \Big)\Big],
\]
we can simplify the amplitude factor in (\ref{eq:D1}) as
\begin{align}
\frac{k^2 \overline{ \hat f \big[ \om \big( 1+2\gamma(s,\vy_o) \big)
      \big] }\hat f \big[ \om \big( 1+2\gamma(s,\vy) \big) \big]}{(4
  \pi |\vr(s)-\vy|)^2} \approx \frac{k_o^2 |\hat f(\om_o)|^2}{(4 \pi
  L_\alpha)^2} ,
\label{eq:DD4}
\end{align}
and obtain 
\begin{align}
d^{(\alpha,\beta)}(\Delta s, \Delta \om) \approx \frac{k_o^2 |\hat
  f(\om_o)|^2}{(4 \pi L_\alpha)^2} \sum_{q=1}^Q
\rho^{(\alpha,\beta)}_q \exp \Big[ - 2 i ( k_\beta + \Delta k) \frac{\vr'(s_\alpha^\star+\Delta s)}{c} \cdot
  \Delta \vy_q +  \nonumber \\
  2 i (\om_\beta^\star + \Delta \om)\big[\tau(s_\alpha^\star+\Delta s,\vy_o + \Delta \vy_q) - \tau(s_\alpha^\star+\Delta s
    ,\vy_o)\big] \Big].
    \label{eq:datasp}
\end{align}
Here we have used that $k = k_\beta + \Delta k$, with center wavenumber $k_\beta =
\om_\beta^\star/c$ and offset $\Delta k = \Delta \om/c$. 

The difference between the travel times in the phase in (\ref{eq:datasp}) is approximated in the proof of
Lemma \ref{lem.1} in appendix \ref{ap:proofs}. It remains to expand
the first term in the phase, which is due to the Doppler factor.  We
use (\ref{eq:DD3}) and obtain
\begin{align*}
(k_\beta + \Delta k) \frac{\vr'(s_\alpha^\star + \Delta s)}{c} \cdot
  \Delta \vy_q = k_\beta \frac{V}{c} \Big[\vt_\alpha \cdot \Delta
    \vy_q - \frac{V \Delta s}{R} \vec{\bf n}_\alpha \cdot \Delta
    \vy\Big] + \Delta k \frac{V}{c} \vt_\alpha \cdot \Delta \vy_q +
  \\ O\Big(\frac{V}{c} \frac{a}{R} \frac{\vec{\bf n}_\alpha \cdot
    \Delta \vy_q}{c/b} \Big) ,
\end{align*}
with negligible residual under the assumption
\begin{equation}
\frac{V}{c} \frac{a}{R} \frac{\cY}{c/b} \ll 1.
\label{eq:DD5}
\end{equation}
Recall that $c/b$ is the range resolution, and although we want $\cY
\gg c/b$, the inequality (\ref{eq:DD5}) is easily satisfied because $
a \ll R \sim L_\alpha$ and $V \ll c$.

The generalization of the result  in Lemma \ref{lem.1} is as follows.  We have the linear system of equations 
\begin{equation}
\bA^{(\alpha,\beta)} \brho^{(\alpha,\beta)} = \bd^{(\alpha,\beta)},
\end{equation}
where 
the reflectivity vector 
$\brho^{(\alpha,\beta)}$ with entries $\rho_q^{(\alpha,\beta)}$ is mapped to 
the data vector $\bd^{(\alpha,\beta)}$ with entries $d^{(\alpha,\beta)}(\Delta s_j,\Delta \om_l)$ by
the reflectivity-to-data matrix $\bA^{(\alpha,\beta)}$.
The entries of $\bA^{(\alpha,\beta)}$ are given by 
\begin{align}
A_{j,q}^{(\alpha,\beta)}(\Delta \om_l) =  \frac{k_o^2 |\hat
  f(\om_o)|^2}{(4 \pi L_\alpha)^2} \exp \Big\{ -2 i (k_\beta +\Delta \om_l/c) \Big[ \vm_\alpha
  \cdot \Delta \vy_q + \frac{V}{c} \vt_\alpha \cdot \Delta \vy_q \Big]
\nonumber \\ \, -2 i \frac{k_\beta V
  \Delta s}{L_\alpha} \Big[\vt_\alpha \cdot \mathbb{P}_\alpha \Delta \vy_q - \frac{L_\alpha}{R} \frac{V}{c} 
  \vn_\alpha \cdot \Delta \vy_q\Big]
+ i k_\beta \frac{\Delta \vy_q \cdot \mathbb{P}_\alpha \Delta
  \vy_q}{L_\alpha} \Big\}.
\end{align}
The difference between this reflectivity-to-data matrix and the 
one given by (\ref{eq:lem1})
in Lemma \ref{lem.1} comes from the $V$ dependent terms in the square brackets in 
the phase, due to the Doppler effect.

We extend next the statement of Proposition \ref{lem.2}. We
proceed as in appendix \ref{ap:proofs}, and show that the matrix-matrix equation
(\ref{eq:MMVEQ}), $ \boldsymbol{\mathbb{A}} {\bf X} = \bD$, still applies, with the same definition (\ref{eq:lem2.6}) of the data
matrix $\bD$, 
 \begin{equation}
    D_j^{(\alpha,\beta)}(\Delta \om_l) = \frac{(4 \pi
      L_\alpha)^2}{k_o^2 |\hat f(\om_o)|^2} d^{(\alpha,\beta)}(\Delta
    s_j,\Delta \om_l),
 \nonumber
  \end{equation}
and with the unknown matrix
\begin{align}
X_{q}^{(\alpha,\beta)} = \rho_q^{(\alpha,\beta)} \exp \Big\{-2 i
k_\beta \Big[ \frac{V}{c} \vt_\alpha \cdot \Delta \vy_q + \vm_\alpha
  \cdot \Delta \vy_q \Big] + i k_\beta \frac{\Delta \vy_q \cdot
  \mathbb{P}_\alpha \Delta \vy_q}{L_\alpha} \Big\}.
\label{eq:NewX}
\end{align}
This is under the assumptions that 
\begin{align}
\max_{1 \le \alpha \le \cN_\alpha, 1 \le q \le Q} \frac{V}{c}
\frac{b}{c} \Big|[\vt_\alpha - \vt_1] \cdot \Delta \vy_q\Big| &\ll 1,  \\ 
\max_{1 \le \alpha \le \cN_\alpha, 1 \le q \le Q} \frac{V}{c} \frac{a}{\la_o R} \big| (\vn_\alpha-\vn_1) \cdot \Delta \vy_q \big| &\ll1, 
\end{align}
which are similar to (\ref{eq:lem2.1})-(\ref{eq:lem2.2}), and 
easier to satisfy for smaller $V$.  The expression of the entries of
the reflectivity-to-data matrix is a simple modification of that in
equation (\ref{eq:lem2.8}),
\begin{align}
\label{eq:matrixdp}
\mathbb{A}_{j,q}(\Delta \omega_l) = \exp \Big[ &-2 i \frac{\Delta
    \om_l}{c} \Big(\vm_1 \cdot \Delta \vy_q + \frac{V}{c} \vt_1 \cdot
  \Delta \vy_q \Big) \nonumber \\ &-2 i k_1 \frac{V \Delta s_j}{L_1} \Big(\vt_1 \cdot
  \mathbb{P}_1 \Delta \vy_q - \frac{L_1}{R}\frac{V}{c} \vn_1 \cdot \Delta \vy_q\Big) \Big].
\end{align}

Thus, the problem can be solved with the MMV approach, as described in
section \ref{sect:MMVAlg}. The Doppler correction has two effects: It
gives an extra rotation in the cross-range direction of the imaging
window (the first phase term in (\ref{eq:NewX}), involving $V$), and two extra phase
factors (involving $V$ inside the parentheses) in the reflectivity-to-data matrix $\mathbb{A}$ in (\ref{eq:matrixdp}).

\section{Summary}
\label{sect:summary}
We have introduced and analyzed from first principles a synthetic aperture imaging approach for 
reconstructing direction and frequency dependent reflectivities of localized scatterers.  
It is based on two main ideas: The first one is to segment the data over subsets defined by 
carefully calibrated sub-apertures and frequency sub-bands, and formulate the reflectivity reconstruction for 
each subset as an $\ell_1$ optimization problem.  The direction and frequency dependence of the 
reconstructed reflectivity is frozen for each data subset but varies from one subset to another. 
The second idea is to fuse the sub-aperture and sub-band optimizations
by seeking simultaneously from data subsets those reconstructions of the reflectivity that share the same spatial support in the 
image window. This is done with the multiple measurement vector (MMV) formalism, which leads to a matrix
$\ell_1$ optimization problem. The main result of this paper is showing that synthetic aperture imaging
of direction and frequency dependent reflectivities can be formulated and solved efficiently as an MMV problem.

Data segmentation is a natural idea that has been used before for synthetic aperture imaging of frequency dependent
reflectivities \cite{sotirelis2012study,elachi1990radar}. Here we use it for estimating the direction dependence
of the reflectivity, as well. We analyze how the size of the sub-apertures and frequency 
sub-bands in the data segmentation affects the resolution of the reconstructions as well as the
computational complexity of the inversion.  There is a trade-off in resolution in this approach:
On one hand we want to have large sub-apertures and frequency sub-bands to get good spatial, range and cross-range, 
resolution of the reconstructed reflectivity. But on the other hand we also want to have small sub-apertures and 
frequency sub-bands to resolve well the direction and frequency  dependence of the reflectivity.
Small sub-apertures are also desirable so as to get images efficiently using Fourier transforms.  
The MMV formalism that we have introduced in this paper, and the associated algorithm for its
implementation, deal well with these issues, as indicated by the numerical simulations.

Nearly all synthetic aperture imaging is done with reverse time migration algorithms, without
regard to whether the reflectivities that are to be imaged are direction dependent or not. If the
reflectivities are isotropic, then the spatial resolution of the reconstruction improves as 
the aperture increases. But this is not the case with direction dependent reflectivities as 
only part of the synthetic aperture will sense reflectivities from particular locations. This means
that segmenting the data over sub-apertures is natural. The MMV-based imaging algorithm 
introduced in this paper handles automatically signals received by sub-apertures that are
coming from directional reflectivities located in the image window.



\section*{Acknowledgements}
Borcea's work was partially supported by grant \#339153 from
the Simons Foundation and by AFOSR Grant FA9550-15-1-0118. 
Moscoso's work was partially supported by the Spanish MICINN grant
FIS2013-41802-R. Papanicolaou's work was partially supported by
AFOSR grant FA9550-14-1-0275. Tsogka's work was partially
supported by the ERC Starting Grant Project ADAPTIVES-239959 and the
AFOSR grant FA9550-14-1-0275.  \appendix
\section{Derivation of the reflectivity to data model}
\label{ap:proofs}
Here we show that the expression of $A_{j,q}(\om_l)$ in
(\ref{eq:L1.4}) can be approximated by
$A_{j,q}^{(\alpha,\beta)}(\Delta \om_l)$ given in Lemma \ref{lem.1},
for $\om_l = \om_\beta^\star + \Delta \om_l$ and $s_j = s_\alpha^\star
+ \Delta s_j$. For simplicity of notation we drop the indexes $j$ and
$l$ of the frequency and slow time.

It is easy to see from (\ref{eq:L1.2}) and the assumptions $\om_o \gg
b$ and $L_\alpha \gg a \gtrsim \cY$ that
\begin{equation}
  \frac{k^2 |\hat f(\om)|^2}{\big(4 \pi |\vr(s)-\vy|\big)^2} \approx
  \frac{k_o^2 |\hat f(\om_o)|^2}{(4 \pi L_\alpha)^2},
  \label{eq:AP1}
\end{equation}
for $k = \om/c$ and $k_o = \om_o/c$.  It remains to show the phase
approximation
\begin{align}
  2  \om \big[\tau (s,\vy)-\tau(s,\vy_o)\big] \approx -2  k
  \vm_\alpha \cdot \vy - 2  k_\beta V \Delta s \frac{\vt_\alpha \cdot
    \mathbb{P}_\alpha \Delta \vy}{L_\alpha} + k_\beta \frac{\Delta
    \vy \cdot \mathbb{P}_\alpha \Delta \vy}{L_\alpha},
  \label{eq:AP2}
\end{align}
where $ \om = \om_\beta^\star + \Delta \om$ lies in the frequency
sub-band of width $b$, $s = s_\alpha^\star+\Delta s$ is in the
sub-aperture of size $a$ and $\vy = \vy_o + \Delta \vy$ is in
$\mathcal{Y}$.

We begin by expanding the travel time in $\Delta \vy$,
\begin{align*}
  \Phi &= 2 \om \big[\tau (s,\vy)-\tau(s,\vy_o)\big] \\&= - 2 k
  \vm(s,\vy_o) \cdot \Delta \vy + \frac{k}{|\vr(s)-\vy_o|} \Delta
  \vy \cdot \big[ I - \vm(s,\vy_o) \vm^T(s,\vy_o)\big] \Delta \vy +
  \mathcal{E}_1,
\end{align*}
with small residual
\[
\mathcal{E}_1 = O \Big( \frac{\cYp^2 \cY}{\la_o L_\alpha^2} \Big) \ll 1,
\]
by assumption (\ref{eq:M10}) and $\cYp \lesssim a$, inferred from
(\ref{eq:M5}). Here we used the expression of the gradient
\[
\nabla_\vy |\vr(s)-\vy| = - \frac{\vr(s)-\vy}{|\vr(s)-\vy|} = - \vm(s,\vy),
\]
the Hessian
\[
\nabla_\vy \otimes \nabla_\vy |\vr(s)-\vy| = \frac{1}{|\vr(s)-\vy|}
\Big[ I - \vm(s,\vy)\vm^T(s,\vy)\Big],
\]
and
\begin{align*}
\sum_{i,j,q=1}^3 \Delta y_i \Delta y_j \Delta y_q
\partial^3_{y_i,y_j,y_q} |\vr(s)-\vy| &=  \frac{3 \vm_\alpha \cdot
  \Delta \vy}{|\vr(s)-\vy|^2} \big[ |\Delta \vy|^2 - \big(\vm_\alpha
  \cdot \Delta \vy \big)^2 \big].
\end{align*}
Next, we expand in
$\Delta \om = \om - \om_\beta^\star$ and obtain
\begin{align*}
  \Phi &= - 2(k_\beta + \Delta k) \vm(s,\vy_o) \cdot \Delta \vy +
  \frac{k_\beta}{|\vr(s)-\vy_o|} \Delta
  \vy \cdot \big[ I - \vm(s,\vy_o) \vm^T(s,\vy_o)\big] \Delta \vy +
  \mathcal{E}_2,
\end{align*}
where $\Delta k = \Delta \om/c$ and 
\[
\mathcal{E}_2 = \mathcal{E}_1 + O \Big(\frac{b}{\om_o}
\frac{\cYp^2}{\la_o L_\alpha}\Big) \ll 1.
\]
The last estimate is by assumption (\ref{eq:M8}).  Finally, we expand
in $\Delta s = s - s_\alpha^\star$, and recalling the notation in
section \ref{sect:setup}, we get
\begin{align}
  \Phi = -2(k_\beta + \Delta k) \vm_\alpha \cdot \Delta \vy - 2
  k_\beta \frac{V \Delta s}{L_\alpha} \vt_\alpha \cdot
  \mathbb{P}_\alpha \Delta \vy + \nonumber \\ k_\beta \frac{\Delta \vy
    \cdot \mathbb{P}_\alpha \Delta \vy}{L_\alpha} +
  \mathcal{E}. \label{eq:Phi}
\end{align}
The residual is the sum of four terms
\[
\mathcal{E} = \mathcal{E}_2 + \mathcal{E}_3 + \mathcal{E}_4 +
\mathcal{E}_5,
\]
with $\mathcal{E}_2$ given above. The term $\mathcal{E}_3$ comes from
the quadratic part of the expansion of $k_\beta \vm(s,\vy_o) \cdot
\Delta \vy$,
\begin{align*}
\mathcal{E}_3 \sim k_\beta (V \Delta s)^2 \Big[ \frac{\vn_\alpha \cdot
    \mathbb{P}_\alpha \Delta \vy }{R L_\alpha} + \frac{ \vt_\alpha
    \cdot \big[ \vm'_\alpha \vm_\alpha^T + \vm_\alpha
      (\vm'_\alpha)^T\big] \Delta \vy}{V L_\alpha} +
  \frac{\big(\vt_\alpha\cdot\mathbb{P}_\alpha \Delta \vy\big)
    \big(\vt_\alpha \cdot \vm_\alpha\big)}{L_\alpha^3}
\Big].
\end{align*}
Here $\sim$ denotes order of magnitude, and the primes denote
derivative with respect to $s$. The unit vector $\vn_\alpha $ is
normal to $\vt_\alpha$,  in the plane defined by $\vt_\alpha$ and the center of curvature 
of the trajectory of the platform. It enters the definition 
\begin{equation}
\label{eq:defTp}
\vt'_\alpha = -\frac{V \vn_\alpha}{R},
\end{equation}
 where $R \sim L_\alpha$ is the radius of curvature. Moreover
\begin{equation}
\vm'_\alpha = \frac{V}{L_\alpha} \mathbb{P}_\alpha \vt_\alpha.
\end{equation}
We conclude that
\[
\mathcal{E}_3 = O \Big( \frac{a^2 \cYp}{\la_o L_\alpha^2} \Big) +
O \Big( \frac{a^2 \cY}{\la_o L_\alpha^2} \Big) \ll 1,
\]
where the inequality is by assumption (\ref{eq:M10}).

The term $\mathcal{E}_4$ in the residual is
\begin{align*}
  \mathcal{E}_4 \sim \frac{\Delta \om}{c} \frac{V \Delta s}{L_\alpha}
  \vt_\alpha \cdot \mathbb{P}_\alpha \Delta \vy = O
  \Big(\frac{b}{\om_o} \frac{a \cYp}{\la_o L_\alpha} \Big) \ll 1, 
\end{align*}
by assumption (\ref{eq:M8}), and the last term $\mathcal{E}_5$ comes
from the expansion of the quadratic term in $\Delta \vy$ in the
expression of $\Phi$. We estimate it as
\begin{align*}
  \mathcal{E}_5 = O \Big(\frac{a \cY \cYp}{\la_o L_\alpha^2} \Big) \ll 1,
\end{align*}
where we used assumption (\ref{eq:M10}). The statement of Lemma
\ref{lem.1} follows from (\ref{eq:AP1}) and (\ref{eq:Phi}). $\Box$

Proposition \ref{lem.2} follows easily from the expression
(\ref{lem.1}) of $A_{j,q}^{(\alpha,\beta)}$ and assumptions
(\ref{eq:lem2.1}) and (\ref{eq:lem2.2}). Writing the linear system
(\ref{eq:lem.1p}) component-wise we get
\begin{align*}
  \sum_{q=1}^Q X_q^{(\alpha,\beta)} \exp \Big[ -2 i
    \frac{\Delta \om_l}{c} \vm_\alpha \cdot \vy_q - 2 i k_\beta V
    \Delta s_j \frac{ \vt_\alpha \cdot \mathbb{P}_\alpha \Delta
      \vy_q}{L_\alpha} \Big] = D_j^{(\alpha,\beta)}(\Delta \om_l),
\end{align*}
with $X_q^{(\alpha,\beta)}$ given in (\ref{eq:lem2.4}) and
$D_j^{(\alpha,\beta)}$ defined in (\ref{eq:lem2.6}). The result
(\ref{eq:lem2.8}) follows from this equation and assumptions
(\ref{eq:lem2.1}) and (\ref{eq:lem2.2}). $\Box$

\section{Inner products for rows and columns of the reflectivity-to-data matrix}
\label{ap:discr}

Here we analyze the relation between the discretization of the imaging
window $\mathcal{Y}$ and the linear independence of the columns of the
reflectivity to data matrix. This is done by computing inner products of
of normalized rows and columns of the reflectivity-to-data matrix. If the column inner
products multiplied by the number of elements in the support of the reflectivities are below a threshold 
then the MMV algorithm will give an exact reconstruction, in the noiseless case \cite{chai2014imaging}.

We consider the restriction to a data
subset, defined by a sub-aperture and frequency sub-band satisfying
the assumptions in section \ref{sect:MMVRed}. Thus, we work with matrices
$\bA^{(\alpha,\beta)}$, but to simplify notation we drop the indexes
$(\alpha,\beta)$.

Let us denote by ${\bf a}_q$ the $q-$th column of matrix $\bA$ and
calculate the inner product
\[
\left< {\bf a}_q, {\bf a}_{q'} \right> = \Big(\bA^\star
\bA\Big)_{q,q'} = \sum_{j = 1}^{n_s} \sum_{l=1}^{n_\om}
\overline{A_{j,q}(\Delta \om_l)} A_{j,q'}(\Delta \om_l).
\]
Using Lemma \ref{lem.1} we get
\begin{align*}
\left< \frac{{\bf a}_q}{\|{\bf a}_q\|}, \frac{{\bf a}_{q'}}{\|{\bf
    a}_{q'}\|} \right> = \exp \Big[ -2 i k_\beta \vm_\alpha \cdot
  (\vy_{q'}-\vy_q) + \frac{i k_\beta\Big( \Delta \vy_{q'}
    \mathbb{P}_\alpha \Delta \vy_{q'} - \Delta \vy_{q}
    \mathbb{P}_\alpha \Delta \vy_{q}\Big)}{L_\alpha} \Big] \times
\\\frac{1}{n_s n_\om}\sum_{j = 1}^{n_s} \sum_{l=1}^{n_\om} \exp
\Big[-\frac{2 i \Delta \om_l}{c} \vm_\alpha \cdot (\vy_{q'}-\vy_q) -
  \frac{2 i k_\beta V \Delta s_j}{L_\alpha} \vt_\alpha \cdot
  \mathbb{P}_\alpha (\vy_{q'}-\vy_q)\Big],
\end{align*}
where we normalized the columns by their Euclidian norm. The sums can
be approximated by integrals over the frequency band and aperture,
as long as they are sampled at intervals $h_\om$ and $h_s$ satisfying
\[
\frac{h_\om}{b} \frac{|\vm_\alpha \cdot (\vy_q-\vy_{q'})|}{c/b} \ll 1,
\qquad \frac{V h_s}{a} \frac{|\vt_\alpha \cdot \mathbb{P}_\alpha
  (\vy_{q'}-\vy_q)|}{\la_o L/a} \ll 1.
\]
We obtain after taking absolute values that
\begin{align}
  \left|\left< \frac{{\bf a}_q}{\|{\bf a}_q\|}, \frac{{\bf
      a}_{q'}}{\|{\bf a}_{q'}\|} \right>\right| \approx
  \left|\mbox{sinc} \Big( \frac{b}{c} \vm_\alpha \cdot
  (\vy_{q'}-\vy_q) \Big) \mbox{sinc}\Big( \frac{k_o a}{L_\alpha}
  \vt_\alpha \cdot \mathbb{P}_\alpha (\vy_{q'}-\vy_q)\Big)\right|.
\end{align}
This is small for $q \ne q'$ when we sample the imaging window
$\mathcal{Y}$ in steps that are larger than the resolution limits
$c/b$ in range and $\la_o L/a$ in cross-range.

A similar calculation can be done for the rows of $\bA$, denoted
by ${\bf a}_{(j,l)}$. We have
\[
  \left< {\bf a}_{(j',l')},{\bf a}_{(j,l)} \right> = \Big(\bA \bA^*
  \Big)_{(j,l),(j',l')} = \sum_{q=1}^Q \overline{A_{j',q}(\Delta
    \om_{l'})} A_{j,q}(\Delta \om_l),
\]
and using Lemma \ref{lem.1} we get
\begin{align*}
\left|\left< \frac{ {\bf a}_{(j',l')}}{\|{\bf a}_{(j',l')}\|},
  \frac{{\bf a}_{(j',l')}}{\|{\bf a}_{(j',l')}\|} \right> \right|
  = \frac{1}{Q} \sum_{q=1}^Q \exp \Big[ \frac{2i
      (\om_{l'}-\om_{l})\vm_\alpha \cdot \Delta \vy_q}{c} + \\\frac{2 i
      k_\beta V (s_{j'}-s_{j}) \vt_\alpha \cdot \mathbb{P}_\alpha \Delta
      \vy_q}{L_\alpha} \Big].
\end{align*}
Furthermore, for discretizations of the imaging window in steps
$h$ in range and $h^\perp$ in cross-range, satisfying
\[
\frac{|\om_{l'}-\om_l|}{b} \frac{h}{c/b} \ll 1, \qquad 
\frac{V |s_{j}-s_{j'}|}{a} \frac{h^\perp}{\la_o L_\alpha/a} \ll 1,
\]
we can approximate the sum over $q$ by an integral over the imaging
window and obtain
\begin{align*}
\left|\left< \frac{ {\bf a}_{(j',l')}}{\|{\bf a}_{(j',l')}\|},
\frac{{\bf a}_{(j',l')}}{\|{\bf a}_{(j',l')}\|} \right> \right|
\approx \left| \mbox{sinc} \Big(\frac{(\om_{l'}-\om_l) \cY}{c}\Big)
\mbox{sinc} \Big(\frac{k_o V (s_{j'}-s_j) \cYp}{L_\alpha}\Big) \right|.
\end{align*}
This result shows that the inner product of the rows is small when the
frequency is sampled in steps larger than $\cY/c$ and the slow time is
sampled in steps larger than $(1/V)/(\la_o \cYp/L_\alpha)$.


\bibliographystyle{siam} \bibliography{SPARSE}

\end{document}